\documentclass{article}

\begin{document}

\title{A journey through the  mathematical world of Karen Uhlenbeck}
\author{Simon Donaldson}
\date{\today}
\maketitle


\tableofcontents
\newcommand{\bC}{{\bf C}}
\newcommand{\bR}{{\bf R}}
\newcommand{\bZ}{{\bf Z}}
\newcommand{\bA}{{\bf A}}
\newcommand{\bP}{{\bf P}}
\newcommand{\tA}{\tilde{A}}
\newcommand{\tu}{\tilde{u}}
\newcommand{\tzeta}{\tilde{\zeta}}
\newcommand{\tphi}{\tilde{\phi}}
\newcommand{\tL}{\tilde{L}}
\newcommand{\tR}{\tilde{R}}
\newcommand{\tG}{\tilde{G}}
\newcommand{\uA}{\underline{A}}
\newcommand{\cE}{{\cal E}}
\newcommand{\olambda}{\overline{\lambda}}
\newcommand{\of}{\overline{f}}
\newcommand{\oz}{\overline{z}}
\newcommand{\ow}{\overline{w}}
\newcommand{\db}{\overline{\partial}}
\newcommand{\hatF}{\widehat{F}}
\newcommand{\hatE}{\widehat{E}}
\newcommand{\rad}{{\rm ad}}
\newcommand{\cH}{{\cal H}}
\newcommand{\cM}{{\cal M}}
\newcommand{\cZ}{{\cal Z}}
\newcommand{\cD}{{\cal D}}
\newcommand{\cF}{{\cal F}}
\newcommand{\thalf}{{\small \frac{1}{2}}}
\newcommand{\tthreetwo}{{\small \frac{3}{2}}}
\newcommand{\teighth}{{\small \frac{1}{8}}}
\newcommand{\tquart}{{\small \frac{1}{4}}}
\newcommand{\tthreequart}{{\small \frac{3}{4}}}
\newcommand{\tfourthree}{{\small \frac{4}{3}}}
\newcommand{\ttwothree}{{\small \frac{2}{3}}}
\newtheorem{thm}{Theorem}
\newtheorem{prop}{Proposition}
\newtheorem{lem}{Lemma}
\section{Introduction}

In this article\footnote{To appear in H. Holden and R. Piene (eds.): \lq\lq The Abel Prize 2018-2022'', Springer 2023}
 we discuss some of Karen Uhlenbeck's most prominent mathematical results.  Uhlenbeck's publications range across many mathematical areas,
 including differential geometry and geometric analysis, elliptic and hyperbolic partial differential equations and integrable systems. In this article we only attempt to describe some part of this range. The main omissions are that we say nothing about her work on wave and Schroedinger maps,  and very little about integrable systems. The core of the article is contained in  Sections 3, 4, 5 and 6 which give an account of some highlights of Uhlenbeck's work on the analysis of harmonic maps and Yang-Mills connections. In each case we begin with the theory for the \lq\lq critical dimension''---in Sections 3 and 5---before going on to higher dimensions (in Sections 4 and 6). This body of work has been absolutely fundamental in the developments of geometric analysis over the past 40 years and with an impact
 that extends to many fields, from  symplectic geometry and low-dimensional topology to Quantum Field Theory and the mathematics of liquid crystals.
  At the beginning and end of the article we discuss two other contributions of Uhlenbeck which take somewhat different directions to that in the core; each very influential and highly-cited. The first (in Section 2) is a paper on   nonlinear elliptic PDE theory and the other (in Section 7) is on  integrable systems aspects of  harmonic maps from surfaces to  Lie groups. 

The author has written another review \cite{kn:SKDNotices} of some of  Karen Uhlenbeck's mathematical work, which focused on  variational methods. While there is  overlap with the current article   we have made the focus here different and sought to avoid duplication. At some points in this article we refer to \cite{kn:SKDNotices} for further discussion of  literature and background.

\section{Nonlinear systems and p-harmonic functions}
\subsection{A regularity theorem}
We begin our tour by discussing the 1977 {\it Acta Mathematica} paper \cite{kn:UActa} of Uhlenbeck which was one of her  first papers with a focus on \lq\lq hard'' PDE theory.  To set the scene for this, recall that the Laplace operator $\Delta$ on functions on $\bR^{n}$ is the Euler-Lagrange operator associated to the Dirichlet energy, the integral of $\vert du\vert^{2}$.
  The solution of the boundary value problem for
   a harmonic function on a domain with prescribed boundary values minimises the Dirichlet energy over the set of all functions with those boundary values. A generalisation is to  take any $p>1$ and the functional defined by  the integral of $\vert du\vert^{p}$.  
The associated Euler-Lagrange equation is the nonlinear {\it p-Laplace} equation
     \begin{equation} d^{*} ( \vert d u\vert^{p-2} d u) =0 . \end{equation}
     The existence of weak solutions to this equation, lying in the Sobolev space $L^{p}_{1}$ and with prescribed boundary values,  is relatively straightforward but the question of the regularity of these weak solutions is very subtle. The equation
(1) is a degenerate elliptic equation at points where the derivative of $u$ vanishes. We can write the equation as
  $$  \Delta u  + (p-2) \sum \nu_{i} \nu_{j} \frac{\partial^{2} u}{\partial x_{i}\partial x_{j}}=0  , $$
where $\nu$ is the unit vector field 
$$    \nu= \frac{1}{\vert d u\vert} du, $$
and $\nu$ will usually be discontinuous at zeros of $du$. This means that one cannot expect solutions of the $p$-Laplace equation to be smooth at such zeros. For example the function $u(x)= \vert x\vert^{\beta}$ with $\beta=(p-n)/(p-1)$ is a solution.

 Uhlenbeck's Acta paper established a central result on this regularity question, as a particular case of a more general theory,   showing that  the derivative of a $p$-harmonic function satisfies a $C^{,\alpha}$ Holder estimate for some $\alpha$ depending on $p,n$. This particular result had been obtained before by Ural'ceva, appearing in Russian \cite{kn:Ural}. But the theory developed in Uhlenbeck's paper covers much more than this model case, as we will now explain.

 Let $M$ be a smooth manifold, with its complex of differential forms
$$   \Omega^{0}\stackrel{d}{\rightarrow} \Omega^{1}\stackrel{d}{\rightarrow}
\Omega^{2} \dots  $$
If $M$ is compact and Riemannian then in Hodge theory the harmonic representative $\omega$ of a $k$-dimensional de Rham cohomology class is characterised as the minimiser of the $L^{2}$ norm over all representatives of that class.
 It satisfies the equations $d\omega=0, d^{*}\omega=0$. In the spirit of the discussion above, it is natural to consider the generalisation of this where one takes a positive function $g$ on $\bR$ and minimises
$$  \int_{M} g(\vert \omega\vert). $$
For example we could take $g(\vert \omega\vert)= \vert\omega\vert^{p}$. 
For a small variation $\omega + d\alpha$ in the fixed cohomology class
$$  g(\vert \omega +d\alpha\vert)= g(\vert \omega\vert)+ ( d\alpha, \rho(\vert \omega\vert) \omega), +O(\alpha^{2})$$
where $\rho$ is the function $\rho(t)= g'(t)/t$. So a minimiser satisfies the Euler-Lagrange equation
\begin{equation} d^{*}( \rho(\vert \omega\vert) \omega)=0, \end{equation}
in addition to the closed condition $d\omega=0$.

 This nonlinear generalisation of Hodge Theory was studied by the  Sibners \cite{kn:Sib} who established that for a large class of functions $g$ there is indeed a unique minimiser, giving a weak solution of the equation (2). When $k=1$ the closed $1$-form $\omega$ can be written locally as the derivative of a function $u$ and, when $g(\vert\omega\vert)=\vert\omega\vert^{p}$ we get back to the  $p$-harmonic equation. The equations derived from other functions $g$ arise in the theory of gas dynamics, as explained in \cite{kn:Sib}.

Uhlenbeck's main theorem in \cite{kn:UActa} asserts that, for a large class of functions $g$, these weak solutions are H\"older continuous. In fact her result is formulated for more general elliptic complexes, such as the   $\db$-complex. The conditions imposed on the function $g$   are, roughly speaking, that it should have the  character of  $(\vert \omega\vert^{2} +c)^{p/2}$ for some $c\geq 0$. In the case when $c=0$ the equation becomes degenerate
 at the zeros of $\omega$, just as we saw for the $p$-harmonic equation. But even in easier case when $c>0$ the result was new.  The force of the result in that case is that it applies to {\it systems} of PDE rather than to an equation for a single function. Problems 19 and 20 in  Hilbert's 1900 problem list  asked about the existence and regularity of solutions
 to   variational problems. In the 1950's, De Giorgi  and Nash   obtained very general results on the regularity of  weak solutions to elliptic variational problems  for a single function,  but examples show that these results do not extend to systems: we refer to the discussion in \cite{kn:Giaquinta}, Chapter II. Such regularity questions form a theme running through much of Uhlenbeck's work discussed in this article.

In the remainder of this Section we  sketch some of the main parts of Uhlenbeck's arguments in \cite{kn:UActa}.  To simplify our presentation we will consider only the case of the de Rham complex and the function $g(\vert \omega\vert)=\vert \omega\vert^{p}$. (Uhlenbeck's results are stated for domains in $\bR^{n}$ but the proofs should extend to general Riemannian manifolds.) The theorem we are discussing then is:
\begin{thm}

Let $\omega$ be a $k$-form on the domain $U\subset \bR^{n}$ with coefficients in $L^{p}$ which is a weak solution of the equations
\begin{equation}   d\omega=0\ \ \ \  d^{*}(\vert \omega\vert^{p-2}
 \omega)=0. \end{equation} Then $\omega$ is H\"older continuous on compact subsets of $U$.
\end{thm}
\subsection{A differential inequality}
The foundation of Uhlenbeck's proof is an idea which we will meet many other times below: artful use of differential inequalities for {\it functions} can produce important information results about solutions of complicated {\it systems} of PDE. To set things up, given an exterior $k$-form $\nu$ with $\vert \nu\vert=1$ define a symmetric matrix $(a_{ij})$ by the inner products
\begin{equation}  a_{ij}= (dx_{i}\wedge \nu, dx_{j}\wedge \nu) \end{equation}
This is clearly a positive symmetric matrix $(a_{ij})\geq 0$ and we also have an upper bound $(a_{ij})\leq (\delta_{ij})$. Indeed if we define $(b_{ij})$ by
$$    b_{ij}= ( I_{i} \nu, I_{j} \nu), $$
where $I_{i}$ is the operation of contraction with $\frac{\partial}{\partial x_{i}}$ then $(b_{ij})\geq 0$ and it is a basic fact of exterior algebra that $a_{ij}+ b_{ij}=\delta_{ij}$. Now, given a $k$-form $\omega$ on $U\subset \bR^{n}$ we apply this at each point where $\omega\neq 0$, taking $\nu=\omega/\vert \omega\vert$, so we get functions $a_{ij}$, defined away from these zeros.
Let $\Lambda$ be the linear differential operator, depending on $\omega$,
$$   \Lambda(f)= \sum_{ij} \frac{\partial}{\partial x_{i}} \left( a_{ij} \frac{\partial f}{\partial x_{j}}\right), $$
and define an operator $L$ by $$L= \Delta + \frac{p-2}{p-1} \Lambda. $$
 (In this article we use the \lq\lq analysts'' sign convention for the Lapacian: $\Delta =\sum\frac{\partial^{2}}{\partial^{2} x_{i}^{2}}$.)

\begin{prop}

Let $\omega$ be an $L^{p}_{\rm loc}$ solution of the equations (3)
 and define $H=\vert \omega\vert^{p}$ and $\theta= \vert \omega\vert^{(p-2)/2} \omega$. Then, with the operator $L$ as defined above,
$$  L(H) \geq c_{p} \vert \nabla \theta \vert^{2}, $$
where $c_{p}= 4p/(p-1)(p+2)$. 
\end{prop}

(In what follows we calculate  as though all derivatives are defined in the elementary sense. Of course  one has to make precise the meaning of the formula when $\omega$ is {\it a priori} only in $L^{p}$, but we will ignore such technicalities here.) 

\

To establish the inequality we begin with the formula
$$  \nabla_{i} \vert \omega\vert^{p}= p \vert \omega\vert^{p-2} (\nabla_{i}\omega, \omega). $$
Replacing $p$ by $(p-2)$ we get
 $$ \nabla_{i} \left( \vert\omega\vert^{p-2}\omega\right) = 
\vert\omega\vert^{p-2} \nabla_{i}\omega + (p-2)\vert \omega\vert^{p-4} (\nabla_{i}\omega, \omega) \omega, $$
so $$
  (\ \omega, \nabla_{i} ( \vert\omega\vert^{p-2}\omega)\ )= (p-1)\vert\omega\vert^{p-2} (\omega, \nabla_{i}\omega)= \frac{p-1}{p} \nabla_{i}\vert \omega\vert^{p}. $$
Thus \begin{equation}\Delta \vert\omega\vert^{p}= \sum \nabla_{i}\nabla_{i} \vert \omega\vert^{p}=
=\frac{p}{p-1} \sum \nabla_{i}\left( \omega, \nabla_{i}( \vert\omega\vert^{p-2}
\omega)\right). \end{equation}
Now recall the basic fact of Hodge theory, that the Laplacian on $k$-forms on $\bR^{n}$ has two expressions $\Delta= -(d^{*} d+ d d^{*})= \sum \nabla_{i}\nabla_{i}$. (This is the same as the statement that $a_{ij}+ b_{ij}=\delta_{ij}$ in the preceding discussion.) From (5) we get
\begin{equation}  \Delta \vert \omega\vert^{p} = \frac{p}{p-1}(
 P+Q) \end{equation}
where
$$  P= (\omega,\Delta (\vert \omega\vert ^{p-2}\omega))$$
and
$$ Q= (\nabla_{i}\omega, \nabla_{i} (\vert \omega\vert^{p-2} \omega)). $$

To understand the term  $P$ it is convenient to consider a compactly supported test function $f$ and the $L^{2}$ inner product
$ \langle P, f\rangle_{L^{2}}$. This is
$$   \langle f\omega,  \Delta \vert \omega\vert ^{p-2}\omega\rangle_{L^{2}}. $$
Since $d^{*}(\omega\vert ^{p-2}\omega)=0$ we can write this as
$$ \langle d(f\omega),  d( \vert \omega\vert ^{p-2}\omega\rangle_{L^{2}}, $$
and since $d\omega=0$ this becomes
$$  \langle df\wedge \omega, d(\vert \omega\vert^{p-2})\wedge \omega\rangle_{L^{2}}.$$

Now $$d\vert \omega\vert^{p-2}= \frac{p-2}{p} \vert \omega\vert^{-2} d\vert\omega\vert^{p}=\frac{p-2}{p} dH$$
so we can write this as
$$  \langle P, f\rangle_{L^{2}}= \frac{p-2}{p}\langle df\wedge \frac{\omega}{\vert \omega\vert}, dH\wedge\frac{\omega}{\vert \omega\vert}
\rangle_{L^{2}}. $$
By the definition of the $(a_{ij})$ this equation is
$$  \langle P, f\rangle_{L^{2}}= \int \sum a_{ij} \nabla_{i} f \nabla_{j} H. $$
Since this is true for all $f$ we have
$$  P= -\frac{p-2}{p} \sum \nabla_{i}(a_{ij}\nabla_{j} H) =\frac{p-2}{p} \Lambda ( H). $$
Now (6) becomes
$   L (H)= Q$. Turning attention to the term $Q$, we have
$$  Q= \vert \omega\vert^{p-2}\left(\vert \nabla \omega\vert^{2}+ (p-2) \sum_{i}(\frac{\omega}{\vert \omega\vert}, \nabla_{i}\omega)^{2}\right). $$
While, for $\theta=\vert\omega\vert^{p/2-1}\omega$,
$$  \vert \nabla \theta\vert^{2}= \vert \omega\vert^{p-2}
 \left(\vert \nabla\omega\vert^{2} +
 ((p-2)+\left( \frac{p-2}{2}\right)^{2})  \sum_{i}(\frac{\omega}{\vert
\omega\vert}, \nabla_{i}\omega)^{2}\right), $$
and, comparing the two, we see that $  Q\geq \frac{4}{p+2} \vert \nabla\theta\vert^{2}$, completing the verification of Proposition 1.

\

The significance of this Proposition 1 is that the divergence-form operator $L$ is uniformly elliptic  with bounded measurable coefficients (we assume that $(a_{ij})$ is defined almost everywhere). That is, the eigenvalues of the coefficient matrix $\delta_{ij}+ \frac{p-2}{p-1}a_{ij}$ of $L$ are bounded between $1- \vert (p-2)/(p-1)\vert $ and $1+ \vert (p-2)/(p-1)\vert $.
This opens the way to apply  the deep theory from the 1950's of elliptic operators with measurable coefficients, which were the foundation for the results of 
di Giorgi and Nash mentioned above.

\

\subsection{Outline of  proof of Theorem 1}

One first issue is to show that the form $\omega$ in Theorem 1 is bounded
on compact subsets of the domain $U$, but we will pass over this to focus
on Uhlenbeck's proof of H\"older continuity. (Her proof of boundedness uses
related arguments.)
\

For background, we review some relatively elementary results for the standard Laplace operator.
Let $B'\subset B$ be  balls in $\bR^{n}$, for example the unit ball and the concentric ball of half the radius. Let  $h$ be a positive function on $B$ with $\Delta h\geq 0$. Let $M, M'$ be the suprema of $h$ on the balls $B,B'$ respectively. Then by the definition $M\geq M'$ and the maximum principle implies that $M=M'$ if and only if   $h$ is a constant, in which case $\Delta h=0$. The next Proposition gives two quantitative versions of this.
\begin{prop} There are constants $C_{1}, C_{2}$ such that if $\Delta h\geq \rho\geq 0$ on $B$ then
\begin{enumerate}\item $$   \int_{B'} \rho \leq C_{1} (M-M') ; $$
\item $$  \int_{B'} (M'-h) \leq C_{2} (M-M'). $$
\end{enumerate}
\end{prop}

We give a  proof of the first item of Proposition 2.
  Let $f$ be the solution of $\Delta f =-\rho$ in $B$ with $f=0$ on the boundary of $B$ and let $g= h+f$. Then $\Delta g\geq 0$ and so the maximum principle implies that that the supremum of $g$ on $B'$ is at most that on $\partial B$, which is $M$.
For $x\in B'$ we have 
$$  f(x)= \int_{B} G(x,y)\rho(y) dy, $$
where $G$ is the Green's function, which is positive in the interior of $B$. So there is an $\epsilon>0$ such that for $x,y\in B'$ we have $G(x,y)\geq  \epsilon$ which implies that 
$$  f(x) \geq \epsilon \int_{B'} \rho. $$ So for $x\in B'$
$$h(x)=g(x)-f(x) \leq M- \epsilon\int_{B'} \rho, $$
so
$$  M'\leq M- \epsilon\int_{B'} \rho, $$
which is the desired inequality with  $C_{1}=\epsilon^{-1}$.

It is also  easy to deduce the first item  from the second,  with a slightly different choice of balls.  Let $\chi$ be a cut-off function supported in $B'$, equal to $1$ on some smaller ball $B''\subset B'$.  Then
$$  \int_{B''}\rho \leq \int_{B'} \chi \Delta(h-M')= \int_{B'} (\Delta \chi) (h-M')\leq c \int_{B'} (M'-h) $$
where $c={\rm max} \vert \Delta \chi\vert$.

\

Now suppose that $\omega$ is a $k$-form on the ball $B$ satisfying the equation (3) in Theorem 1 and set $H= \vert \omega\vert^{p},  \theta=\vert\omega\vert^{(p-2)/2}\omega$ as above.  So Proposition 1 gives
$L(H)\geq \rho$ with $\rho= c_{p} \vert\nabla \theta\vert^{2}$.
Let $M, M'$  be the suprema of $H$ on $B,B'$. Uhlenbeck shows that an inequality of the same nature  as the first item in Proposition 2 holds in this situation, so  that, for a suitable constant $C$,
\begin{equation}  \int_{B'}\vert \nabla \theta \vert^{2}   \leq C (M-M').\end{equation}
 The proof uses results of Moser---part of the theory of operators with bounded coefficients mentioned above---and many substantial additional arguments. In fact Moser's result gives the analogue of the first item in Proposition 2 for the operator $L$ and Uhlenbeck obtains the analogue of the first item
  in the manner indicated above, but  additional arguments are required to carry this through because $L$ depends on $\omega$. So we do not  have the same control of
 $\vert L(\chi)\vert $ for the cut-off function $\chi$.  Of course, in all this the techniques required to treat the operator $L$ are quite different from the elementary techniques which suffice for the Laplace operator.

The conclusion is that, in going from the ball $B$ to the smaller ball $B'$, either the supremum of $H$ decreases substantially or $\theta$  is  approximately constant on $B'$, in the sense that the $L^{2}$ norm of $\nabla \theta$ is small.

The other main component in Uhlenbeck's proof is a \lq\lq perturbation theorem'' for solutions which are close to a constant. We state this over a fixed pair of balls $B''\subset B'$.
\begin{prop}
 There are $\epsilon, \kappa>0$ such that if $\omega_{0}$ is a constant form with norm $1$ and $\omega$ is a solution of (3) over  $B'$ with $M'\leq 2$ such that
\begin{equation}  \int_{B'} \vert \omega-\omega_{0}\vert^{2} \leq \epsilon, \end{equation}
then $\omega$ satisfies a H\"older estimate 
$$\vert \omega(x)-\omega(y)\vert \leq \kappa \vert x-y\vert^{1/2}$$
for $x,y\in B''$.
\end{prop}

 This is of the same flavour as \lq\lq small energy'' results which we will encounter throughout this article.  In fact one can go on from this, with a sufficiently small $\epsilon$,  to obtain estimates on all derivatives of $\omega$ over $B''$. The general idea is that the constraint (8) keeps the solution in the regime where the nonlinear equation is well approximated by its linearisation.

\

We now outline how Uhlenbeck puts these components together to prove Theorem 1. 
Notice that our equations (3) are preserved by translations and dilations of $\bR^{n}$ and also by multiplying the solution $\omega$ by a non-zero constant. So the statements above for fixed pairs of balls and---in Proposition 3--- for
$\omega_{0}$ of unit norm,  scale to  corresponding result on balls of arbitrary size and  
for any non-zero constant form $\omega_{0}$.

Suppose again that $\omega$ is a solution over $B$ and normalise so that  $M=1$. Suppose that $M'$ is close to $M$, so that the $L^{2}$ norm of
 $\nabla \theta$ is small by (7). Let $\theta_{0}$ be the average  of $\theta$ over $B'$. The Poincar\'e inequality implies that
 the $L^{2}$ norm of $\theta-\theta_{0}$ is small and it also follows from the hypotheses of $M,M'$ that $\vert \theta_{0}\vert$ will be close to $1$. This is not immediately what is needed to apply Proposition 3, because that needs control of the $L^{2}$ norm of $\vert \omega-\omega_{0}\vert$ for a constant form $\omega_{0}$. But, with additional arguments, Uhlenbeck achieves this control, for $\omega_{0}= \vert \theta_{0}\vert^{2/p-1}\theta_{0}$. 

 The conclusion is that there is some fixed small $\lambda>0$ such that if
$M'>(1-\lambda) M$ then $\omega$ satisfies the hypotheses of the perturbation theorem (after rescaling) over $B'$ and hence a $\thalf$-H\"older estimate over the interior ball $B''$. This number $\lambda$ will determine the H\"older exponent achieved in Theorem 1. (Uhlenbeck remarks on page 238 of \cite{kn:UActa}
{\it \lq\lq it looks like $\lambda$ will be rather small !''})

\

 To prove Theorem 1 we suppose that the domain $U$ contains the unit ball $B=B_{1}$.  It suffices to estimate  $\vert \omega(x)-\omega(0)\vert$ for small $x$, say $\vert x\vert\leq \tquart$. (Strictly, $\omega$ is {\it a priori} only defined almost everywhere, so some some extra words are needed to make sense of pointwise values, but we are ignoring such technicalities here.)
Let $M_{j}$ be the supremum of $H$ on the $2^{-j}$ ball centred at the origin.
If $M_{1}\geq (1-\lambda)M_{0}$ we get a H\"older estimate on $\omega$ over the $\tquart$-ball $B_{\tquart}$ and we are done. If not, we have some definite decrease in the supremum: $M_{1}\leq (1-\lambda)M_{0}$. Now we consider the same alternative for $M_{1}$ and $M_{2}$. If $M_{2}\geq (1-\lambda) M_{1}$ we have our H\"older estimate over $B_{\teighth}$ and in addition we know that for $x$ in $B_{\tquart}$
$$ 
 \vert \omega(0)-\omega(x)\vert 
 \leq 2 M_{2}^{1/p}\leq 2(1-\lambda)^{2/p}M_{0}^{1/p}, $$
 where we have used the facts that $M_{2}$ is the supremum of $
  \vert \omega\vert^{p}$ over the $\tquart$-ball and $\vert \omega(0)-\omega(x)\vert \leq \vert \omega(0)\vert + \vert \omega(x)\vert$. Continue in the same way: either $M_{j+1}\leq (1-\lambda)M_{j}$ for all $j$ or there is some $k$ such that $M_{j+1}\leq (1-\lambda)M_{j}$ for all $j<k$ but $M_{k+1}\geq (1-\lambda)M_{k}$. In either case we get an estimate on $\vert \omega(x)-\omega(0)\vert$ in the manner above and a little bookkeeping shows that this yields the desired H\"older estimate. For example,  consider the first situation when
$M_{j}\leq (1-\lambda)^{j} M_{0}$ for all $j$.
 Then we must have $\omega(0)=0$ and for
$x$ with $2^{-k-1}\leq \vert x\vert\leq 2^{-k}$ we have 
$\vert \omega(x)\vert\leq M_{0}^{1/p} (1-\lambda)^{k/p}$. This gives
$$  \vert \omega(x)\vert \leq K \vert x\vert^{\alpha}, $$
with $\alpha= -p^{-1}\log_{2} (1-\lambda)$ and $K= M_{0}^{1/p}2^{\alpha}$.

\

\

There is an enormous literature in this area, especially on the $p$-harmonic equation.(At the time of writing, the  {\it Acta} paper of Uhlenbeck has 316 citations on MathSciNet.) One subtle question  is the  optimal H\"older exponent. For example, in the case of $p$-harmonic functions in  dimension $n=2$,  Iwaniec and Manfredi show in \cite{kn:IM} that the optimal exponent (for $p\neq 2$) is
$$  \frac{1}{6} \left( \frac{p}{p-1} + \sqrt{ 1+ \frac{14}{p-1} + \frac{1}{(p-1)^{2}}}\right). $$
There are also many papers on the limiting cases $p=1,\infty$. One recent paper of Daskalopoulos and Uhlenbeck \cite{kn:DU2} makes connections between $\infty$-harmonic functions and  Thurston's theory of homeomorphisms between hyperbolic surfaces
 minimising the Lipschitz constant.

\section{Harmonic maps of surfaces}

\subsection{Background}
Let $(M,g)$ and $(N,h)$ be Riemannian manifolds. The harmonic mapping equation for a map $f:M\rightarrow N$ is the Euler-Lagrange equation associated to the energy functional
$$ E(f)= \int_{M} \vert df\vert^{2} , $$
where the norm $\vert df\vert$ is the standard one defined by $g,h$. Familiar cases are when $M$ is $1$-dimensional, where we get geodesics in $N$, and when $N$ is $1$-dimensional where we get harmonic functions on $M$.  Written explicitly in local coordinates $x^{i}$ on $M$ and $y^{\alpha}$ on $N$ the equations are
$$  \Delta_{M} y^{\alpha} + \Gamma^{\alpha}_{\beta \gamma} y^{\beta}_{,i} y^{\gamma}_{,j} g^{ij}, $$
where $\Gamma^{\alpha}_{\beta \gamma}$ are the Christoffel symbols on $N$.
For analysis, it is often convenient to take $N$ to be isometrically embedded in some large Euclidean space $V$, which is possible by Nash's embedding theorem.
 This is never essential but we will use that set-up in this article. Thus
 $f$ can be thought of as  a vector-valued function on $M$, constrained to lie in $N\subset V$. The harmonic mapping condition is that the projection
of $\Delta_{M} f$ to the tangent bundle of $N$ in $V$ is zero. At each point $y$ of $N$ we have a second fundamental form $B_{y}$ which is a symmetric bilinear map $TN_{y}\times TN_{y}\rightarrow\nu_{y}$, where $\nu$ is the normal bundle. The harmonic mapping equation  can be written as
\begin{equation}  \Delta_{M} f + A_{f}(df,df)=0, \end{equation}
where $A_{f}$ is the symmetric bilinear map, at a point $x\in M$,
$$  A_{f}: {\rm Hom}(TM_{x}, TN_{f(x)})\times {\rm Hom}(TM_{x}, TN_{f(x)})\rightarrow V, $$
obtained from $B_{f(x)}$, the metric on $TM_{x}$ and the inclusion $\nu\subset V$. It will sometimes be convenient to extend $A_{f}$, using orthogonal projection from $V$ to $TN$, 
to a bilinear map ${\rm Hom}(TM_{x},V)\times {\rm Hom}(TM_{x},V)\rightarrow V$, depending on $x\in M$ and $f(x)\in N$.

The dimension of $M$ plays a crucial role in the theory of harmonic maps.
In the words of Eells and Lemaire in the Introduction to \cite{kn:EL}{\it  \lq\lq we imagine $M$ made of rubber and $N$ made of marble\dots the map is harmonic if it constrains $M$ to lie on $N$ in a position of elastic equilibrium''}. In that picture we could say that  higher dimensional rubber is weaker and is inclined to tear when searching for an equilibrum position.  When $M$ is the $1$-dimensional circle there is a geodesic in each homotopy class (a rubber band)  minimising energy, but the analogue is not true in higher dimensions. For example it is easy to show that if $M$ is a sphere of dimension $3$ or more then the infimum of energy in any homotopy class of maps from $M$ to $N$ is zero. The critical dimension in the theory is ${\rm dim}\ M=2$. This is bound up with Sobolev inequalities.  Regarding $N$ as isometrically embedded in the Euclidean space $V$,  the energy functional is just the square of the usual $L^{2}$ norm of the derivative of a map $f:M\rightarrow V $, but restricted to maps with image in $N$.
  If ${\rm dim}\ M=1$ the Sobolev inequalities state that maps with derivative in $L^{2}$ are continuous, in fact H\"older continuous with exponent $\thalf$.
A sequence of maps convergent in the Sobolev space $L^{2}_{1}$ converges pointwise and the  constraint that the map takes values in $N$ is preserved in the limit. In higher dimensions this is not true: dimension $2$ is the borderline where  a map with derivative in $L^{2}$ is in $L^{p}$ for all $p$ but not necessarily continuous: evaluation at a point is not well-defined for such a map---it is only defined up to sets of measure zero. The critical nature of dimension two is related to conformal invariance of the energy.  In a general
dimension ${\rm dim}M=n$, if we multiply the metric $g$ by a conformal   factor
$\lambda$ we change $\vert df\vert^{2}$ by $\lambda^{-1}$ and the volume
element by $\lambda^{n/2}$, so when $n=2$ these factors cancel.

\

The main topic of this  Section 3 is the  paper \cite{kn:SU} of Sacks and  which opened up the theory of harmonic maps in the critical dimension 2. 

\subsection{Bubbling}

We begin with an illuminating example of maps from a   flat $2$-torus $M$ to  the standard round $2$-sphere $S^{2}$, both oriented. Consider the homotopy class of maps of degree $1$ from $M$ to $S^{2}$. There is a simple lower bound on the energy of such maps. We can consider $M$ and $S^{2}$ as Riemann surfaces with area forms 
$\omega_{M},\omega_{S^{2}}$. Then we have at each point of $M$:
\begin{equation}  f^{*}(\omega_{S^{2}})\leq \thalf\vert df\vert^{2} \omega_{M}\end{equation} with equality if and only if $df$ is complex linear. This is a simple calculation with $2\times 2$ matrices. If $f$ has degree $1$ then the integral of $f^{*}(\omega_{S^{2}})$ is the area of $S^{2}$. So we get the lower bound
$$  E(f)\geq 2\ {\rm Area}(S^{2}).$$
This lower bound is not achieved, because if it were the map would be holomorphic and by elementary Riemann surface theory there is no degree $1$ holomorphic map from a torus to the Riemann sphere. On the other hand we can construct maps with energy arbitrarily close to this lower bound. Let $D_{r}$ be a small disc of radius $r$ in $M$ centred at a point $x_{0}$ and identify it isometrically with the standard $r$-disc in $\bC$. Now take a very large disc $D_{R}\subset \bC$ and consider it as a subset of $S^{2}$ via the usual description $S^{2}=\bC\cup\{\infty\}$. So the complement of $D_{R}$ in $S^{2}$ is a small disc centred at the point at infinity. Let  $F:S^{2}\rightarrow S^{2}$ be a map which is the identity on most of $D_{R}$ but which collapses the boundary of $D_{R}$ to the point at infinity. It is clear that when $R$ is large we can do this in such a way that the energy of $F$ is as  close as we please to that of the identity map, which is $2 {\rm Area}(S^{2})$. Finally, define a map $f:M\rightarrow S^{2}$ which sends the complement $M\setminus D_{r}$ of $D_{r}$ to  the point $\infty\in S^{2}$ and on $D_{r}$ is the composite
$F\circ \underline{\lambda}$ where $\underline{\lambda}:D_{r}\rightarrow D_{R}$ is multiplication by $\lambda=R/r$. 

The energy of this map $f$ is exactly the same as that of $F$. This follows immediately from the fact that the energy is a conformal invariant for $2$-dimensional domains, it only depends on the conformal class of the metric. Thus  we get a \lq\lq minimising sequence'' $f_{i}$ of degree-$1$ maps from $M$ to $S^{2}$ whose energy tends to the infimum  $2 {\rm Area}\ S^{2}$ by making the construction above with a sequence
$R_{i}\rightarrow \infty$. For large $i$ the image of a small disc in $M$ covers most of $S^{2}$ and away from $x_{0}$ the maps approach the constant harmonic map.

\

What Sacks and Uhlenbeck established is, roughly speaking, that this is the only way that things can go wrong when trying to apply variational arguments to the energy functional on surfaces. More precisely, they consider a $1$-paramter family of deformations of the functional, with parameter $\alpha\geq 1$:
$$  E_{\alpha}(f)= \int_{M} (1+ \vert df\vert^{2})^{\alpha}. $$
(Some formulae would be neater is one used the integral of $\vert df\vert^{2\alpha}$ but this would lead to a degenerate equation and extra difficulties of the kind discussed in Section 2.) When $\alpha=1$ the functional $E_{\alpha}$ is equal to $E(f)$, up to a constant. For $\alpha>1$ the functional controls the $L^{2\alpha}$ norm of the derivative and one has the favourable Sobolev embedding $L^{2\alpha}_{1}\rightarrow C^{0}$. This means that there is a complete \lq\lq Palais-Smale'' variational theory, something which was worked out in the earlier paper \cite{kn:U1} of Uhlenbeck. So the  functional $E_{\alpha}$ attains its minimum in each homotopy class and more generally there must be sufficient critical points to account for the topology of the mapping space, by minimax and Morse theory arguments.
The Sacks and Uhlenbeck strategy is to seek critical points of $E_{1}$ as limits of critical points of the $E_{\alpha}$ as $\alpha\rightarrow 1$. 
The advantage of this approach, compared with studying minimising or minimax sequences for $E$ directly, is that the critical points of $E_{\alpha}$ satisfy an elliptic equation and this improves the convergence properties, as we will see. Even if a minimiser for $E_{1}$ exists there will always be minimising sequences which only converge in a weak sense, not in $C^{\infty}$. 

\

We can now state more precisely one of  the main results of Sacks and Uhlenbeck for maps between compact manifolds $M,N$ with ${\rm dim} \  M=2$.

\begin{thm}
Let $\alpha_{i}\geq 1$ with $\alpha_{i}\rightarrow 1$ as $i\rightarrow \infty$ and let $f_{i}: M\rightarrow N$ be critical points of $E_{\alpha_{i}}$ with $E_{\alpha_{i}}(f_{i})\leq E_{\rm max} $ for some fixed $ E_{\rm max}$. Then, after perhaps passing to a subsequence $\{i'\}$, there is a finite set $S=\{q_{1}, \dots q_{d}\}\subset M$ and a harmonic map $f:M\rightarrow  N$ such that $f_{i'}$ converge to $f$ in $C^{\infty}$ on compact subsets of $M\setminus S$. In addition there are harmonic maps $F_{1}, \dots ,F_{d}: S^{2}\rightarrow N$ such that for each $j\in \{1,\dots,d\}$ a suitable sequence of rescaling of $f_{i'}$ near $q_{j}$ converge on compact subset of $\bC= S^{2}\setminus \{\infty\}$ to $F_{j}$.
\end{thm}

To explain the last statement; we mean that there are points $p_{i'j}$ converging to $q_{j}$ and scale factors $\lambda_{i'j}$ tending to $\infty$ with $i'$ so that if we identify a small disc centred at $p_{ij}$ with a small disc in $\bC$ and compose with a scaling map $\underline{\lambda}_{ij}$ of the kind discussed above the resulting maps converge to $F_{j}$. One says that the sequence of maps $f_{i'}$ is \lq\lq bubbling'' at the points $q_{j}$.

 The statement of our Theorem 2 here does not capture all that Sacks and Uhlenbeck established. For example, they show that the homotopy classes of $f_{i}$ and $f$ in $[M,N]$ differ by a class in $[S^{2}, N]$. But the statement of Theorem 2 gives the general idea. A complete discussion involves the notion of a \lq\lq bubble tree'' of maps, which was worked out later; see for example \cite{kn:Parker}.

This theorem of Sacks and Uhlenbeck implies the existence of harmonic maps in many specific situations. For example, if it is known that there is no harmonic map from $S^{2}$ to $N$ then for any surface $M$ there is a minimising harmonic map in any homotopy class $[M,N]$. One early and famous application came in the proof by Siu and Yau of the Frankel conjecture \cite{kn:SY}. The conjecture was that  projective spaces are the only  compact  complex manifolds admitting  K\"ahler metric with positive biholomorphic sectional curvature.  The result of Sacks and Uhlenbeck shows that there is a nonconstant minimising harmonic map from $S^{2}$ to such a manifold. Siu and Yau proved, by studying the second variation formula and using the curvature condition, that this map is holomorphic  and then the geometry of the resulting family of holomorphic curves shows that the manifold is a projective space. Another important application of the Sacks-Uhlenbeck result, this time in Riemannian geometry, was the \lq\lq sphere theorem'' of Micallef and Moore \cite{kn:MM}, discussed in \cite{kn:SKDNotices}.

By far the greatest impact of the phenomena uncovered by Sacks and Uhlenbeck came in the special case of the holomorphic maps introduced as a tool in symplectic topology by Gromov in 1987.  Here we consider a symplectic manifold $(N,\omega)$ with a compatible almost-complex structure $J$ and resulting Riemannian metric $\vert \xi\vert^{2}= \omega(\xi,J\xi)$. For any oriented Riemannian surface $f:M\rightarrow N$ there is an inequality
$$  2  \int_{M} f^{*}(\omega) \leq  E(f), $$
with equality if and only if $f$ is  holomorphic ({\it i.e.} the derivative at each point is complex linear with respect to the Riemann surface structure on $M$ and  the almost-complex structure $J$). This is a generalisation of (10), in the case when $N$ is a surface. This inequality implies that a holomorphic map minimises  energy in its homotopy class,  so in particular is harmonic. Theorem 2, with all $\alpha_{i}=1$, describes the convergence behaviour of sequences of these holomorphic maps and becomes  the foundation for all of the  applications to symplectic topology such as Gromov-Witten invariants, Lagrangian Floer homology, Fukaya categories\dots. There are many expositions of the theory in this restricted context of holomorphic maps, for example \cite{kn:MDS}, \cite{kn:W}.

\subsection{Small energy}

For simplicity we just discuss the proof of  Theorem 2  in the case when all $\alpha_{i}$ are $1$: the general case does not involve major extra difficulties. The proof has two main components. The first is a \lq\lq small energy'' estimate.
\begin{thm}
Let $N$ be a compact Riemannian manifold and $D$ be the unit disc in $\bC$.
There are $\epsilon, C $ such that if $f:D\rightarrow N$ is harmonic with $E(f)\leq\epsilon$ then
$$  \vert df(0)\vert^{2}\leq C E(f). $$
\end{thm}
More generally, we can choose $\epsilon$ so that if $E(f)\leq \epsilon$ then $E(f)$ controls all derivatives of $f$ on a fixed interior disc, say the $\thalf$-sized disc.  The conformal invariance of the energy implies that these estimates apply {\it with the same small energy threshold} to discs of any size: if $f$ has energy less than $\epsilon$ on the disc $D_{r}$ we get 
$$   \vert \nabla^{k} f\vert\leq c_{k} r^{-k} \sqrt{E(f)}$$
on $D_{r/2}$.
\

The second component is the removability of point singularities.
\begin{thm}
If $f:D\setminus\{0\}\rightarrow N$ is a harmonic map with $E(f)<\infty$ then  $f$ extends smoothly to $D$.
\end{thm}

\

Given these local statements the proof of Theorem 2 is relatively straightforward, using arguments of a kind which we will see several other times in this article. 
 We have $f_{i}:M\rightarrow N$ harmonic with energy less than a fixed number $H$. After passing to a subsequence we can suppose that the energy densities $\vert df_{i}\vert^{2}$ converge as Radon measures: that is, for any continuous function $\phi$ on $M$ the integrals of $\phi \vert df_{i}\vert^{2}$ over $M$ have a limit as $i\rightarrow \infty$. Fix a non-increasing cut-off function $\sigma$ on $[0,\infty)$, equal to $1$ on $[0,1]$ and supported in $[0,2]$. For $x\in M$ and $r>0$, let $\chi_{r,x}$ be the function on $M$
$$  \chi_{r,x}(y)= \sigma ( r^{-1}d(x,y)), $$
where $d(\ ,\ )$ is the Riemannian distance.
So $\chi_{r,x}$ is a smoothing of the characteristic function of the
 $r$-disc $D_{r,x}$ about $x$. Define $\mu(x,r)$ by
$$  \mu(x,r)= {\rm lim}_{i\rightarrow \infty}\int_{M} \chi_{r,x}\  \vert df_{i}\vert^{2} .$$
Then $\mu(x,r)$ is an increasing function of $r$ and has a limit $\mu(x)$ as $r\rightarrow 0$. By construction
\begin{equation}   \int_{D_{x,r}} \vert df_{i}\vert^{2}\leq \int \chi_{r,x} \vert df_{i}\vert^{2}\leq \int_{D_{x,2r}}  \vert df_{i}\vert^{2}. \end{equation}
Let $S$ be the set of points $x$ in $M$ where $\mu(x)>\epsilon/2$. The right hand inequality in (11) implies that there are at most $2 E_{\rm max}/\epsilon$ points $x$ in $S$ (by taking $2r$ less than half the distance between the points, so that the $2r$ discs with these centres are disjoint). On the other hand, taking $\epsilon$ as in Theorem 3, if $\mu(x)<\epsilon/2$ then for some sufficiently small $r=r(x)>0$ and all large enough $i$ the left hand inequality in (11) gives
$$  \int_{D_{x,r}} \vert df_{i}\vert^{2} \leq \epsilon/2, $$
and the small energy theorem  gives estimates on all derivatives of the $f_{i}$ in $D_{x,r/2}$.

From these arguments we get a finite set $S=\{q_{1}, \dots q_{d}\}$ in $M$ such that all derivatives of the  $f_{i}$ are bounded on compact subsets of the complement $M\setminus S$. Taking a subsequence we can assume that the maps converge on the complement to a harmonic map on the punctured manifold with energy at most $C$, and the removal of singularities theorem implies that this extends to a smooth harmonic map $f$ from $M$ to $N$, as stated in the first part of Theorem 2.
 
The second part of Theorem 2 involves the rescaling construction. Fix a point  $q_{j}$. For all large $i$ there must be points near $q_{j}$ where the derivative of $f_{i}$ is large.  Let $p_{i j}$ be a point near $q_{j}$ where
$\vert df_{i}\vert$ is maximal. Define  $\lambda_{ij}$ to be these local  maximal values. Then after rescaling by these factors with centre $p_{ij}$ we get a sequence of harmonic maps $F_{ij}$ defined on  a sequence of large discs in $\bC$ which exhaust $\bC$ as $i\rightarrow \infty$. These maps have  bounded derivative and energy so by the same arguments, perhaps passing to a suitable subsequence, they converge to a harmonic map from $S^{2}$ to $N$. That is, we first get a harmonic map from $\bC$ to $N$ and then apply the removal of singularities theorem at the point at infinity in $S^{2}$.
The limiting  map is not constant since by construction the derivative of the map $F_{ij}$ at the origin has size $1$.

\

We proceed to discuss the proof of the small energy result Theorem 3, leaving that of Theorem 4 to the next subsection.  To recap, we have a  $V$-valued  function $f$ on the disc $D$ which satisfies the PDE $\Delta f= A_{f}(df,df)$ and takes values in the compact submanifold $N$. The discussion below applies to any PDE of this shape, for a smooth map $A$ from $D\times V$ to 
symmetric bilinear maps ${\rm Hom}(T,V)\times {\rm Hom}(T,V)\rightarrow V$,
where $T$ denotes the tangent space of the disc. (In the case at hand we could always extend $A$ in some way to fit into this framework. The fact that, in the case at hand, $f$ maps into the submanifold $N$ is only used in that it gives a bound on $\vert f\vert$.) 

A basic fact of elliptic PDE theory is that for any $q>1$ there is a constant
$K_{q}$ such that for all compactly supported functions $\phi$ on $D$ we have
\begin{equation}   \Vert \phi\Vert_{L^{q}_{2}}\leq K_{q} \Vert \Delta \phi\Vert_{L^{q}}. \end{equation}

(Here, and throughout this article we write $L^{q}_{k}$ for the Sobolev space based on the $L^{q}$ norm of all derivatives of order $\leq k$. We write the norm as $\Vert\ \Vert_{L^{q}_{k}}$ or sometimes $\Vert\ \Vert_{q,k}$ to improve readability.)

Before going on to the proof of Theorem 3 we consider a different situation where we suppose given a solution of the equation with  bound on the $L^{2\alpha}$ norm of $df$ for some $\alpha$ with $\alpha>1$. Then we can do a straightforward \lq\lq bootstrapping'' argument. (Here, and in various other parts of this article, we use the convention that $c$ is a constant that can change from line to line.)

 By choice of the origin in $V$ we may suppose that the integral of $f$ over the disc is $0$. Let $\chi$ be a compactly supported function on the disc, equal to $1$ on an interior disc $D'$. Then we have
$$  \Delta(\chi f)= \chi \Delta f +2 d\chi.df+ \Delta \chi\  f, $$
and
$$   \chi A(df,df)= A(\chi df, df) - A(d\chi \otimes f, df). $$
Thus we have a pointwise bound
\begin{equation}  \vert \Delta(\chi f)\vert \leq c \left( \vert d(\chi f)\vert\  \ \vert df\vert + 
   \vert f\vert +\vert df\vert + \vert f\vert\  \vert df\vert\right), \end{equation}
where $c$ depends only on $\chi$ and the given map $A$. From this we readily obtain, applying the Cauchy-Schwartz inequality,
$$\Vert \Delta(\chi f)\Vert_{{\alpha}}\leq c \left(
 \Vert f \Vert^{2}_{2\alpha,1} + 
  \Vert f \Vert_{2\alpha,1} \ \Vert f \Vert_{{2\alpha}} + \Vert f\Vert_{{\alpha}}\right)
  $$
Since the disc has finite area the $L^{2\alpha}$ norm of $f$ controls the $L^{\alpha}$ norm and the fact that the 
 integral of $f$ vanishes means that the $L^{2\alpha}$ norm of $df$ controls that of $f$. So we get an inequality
$$   \Vert \Delta(\chi f)\Vert_{{\alpha}}\leq c
 \left( \Vert f\Vert^{2}_{2\alpha,1} + \Vert f\Vert_{2\alpha,1}\right), $$
 and hence by the elliptic inequality (12),
$$  \Vert \chi f\Vert_{\alpha,2} \leq c \left( \Vert f\Vert^{2}_{2\alpha,1} + \Vert f\Vert_{2\alpha,1}\right) . $$

If $\alpha<2$ we have a Sobolev embedding $L^{\alpha}_{2}\rightarrow L^{r}_{1}$ with $r=2\alpha/(2-\alpha)$. Thus the inequality above gives an $L^{r}$ bound on $d (\chi f)$ and so an $L^{r}$ bound on $df$ over the interior disc $D'$. Since $r>2\alpha$ this is an improvement on the $L^{2\alpha}$ bound that we started with. If $\alpha>2$ we have a Sobolev embedding $L^{\alpha}_{2}\rightarrow C^{1,\mu}$ for $\mu= 1-2/\alpha$. Starting with our $L^{2\alpha}$ bound on $df$, for any $\alpha>1$ we can iterate such arguments, working on a decreasing sequence of discs, to  get interior  bounds on all derivatives of $f$ in terms of the $L^{2\alpha}$ norm of $df$ over $D$.
 In the same fashion we get a {\it regularity} statement: if we only know at the outset that $f$ is in $L^{2\alpha}_{1}$ we show that in fact it is smooth in the interior of the disc.

\

Now we go on to the proof of Theorem 3. In the situation above {\it any} bound on $\Vert df\Vert_{L^{2\alpha}}$ gives bounds on higher derivatives in the interior. The difference in Theorem 3 is that we only get such a bound when $\Vert df \Vert_{L^{2}}$ is small. Taking  $\alpha\in (1,2)$ and 
$r=2\alpha/(2-\alpha)$ we observe that $1/\alpha= 1/r+1/2$. Thus we can apply Holder's inequality to (13) to get 
$$  \Vert \Delta(\chi f)\Vert_{{\alpha}} \leq c_{1}\Vert d(\chi f)\Vert_{{r}} \Vert df\Vert_{{2}} + c_{2}\left( \Vert df \Vert_{{2}}\Vert f \Vert_{{r}} + \Vert f \Vert_{{\alpha}}\right). $$

By the elliptic inequality and the Sobolev embedding, $\Vert \Delta(\chi f)\Vert_{L^{\alpha}}$ controls
$\Vert d(\chi f)\Vert_{L^{r}}$ and the assumption that $f$ has integral zero means that $\Vert df\Vert_{L^{2}}$ controls both $\Vert f\Vert_{L^{\alpha}}$ and
$\Vert f\Vert_{L^{r}}$. So we have 
$$  \Vert d(\chi f)\Vert_{{r}}\leq c_{3}  \Vert df\Vert_{{2}} \Vert d (\chi f)\Vert_{r} + c_{4}( \Vert df \Vert_{2} + \Vert df\Vert^{2}_{2}). $$
Take $\sqrt{\epsilon}= 1/(2c_{3})$. Then if $\Vert df\Vert_{L^{2}}\leq \sqrt{ \epsilon}$ we re-arrange  to get
$$   \Vert d(\chi f)\Vert_{{r}}\leq 2 c_{4} (\Vert df \Vert_{{2}} + \Vert df\Vert_{{2}}^{2}). $$
We have $r>2$ so, replacing $D$ by the smaller disc $D'$, we are in the position considered before  and we can go on to estimate all higher derivatives in the interior.

\

Arguments with the same structure as this will appear often in this article so we give them a name: \lq\lq {\it critical quadratic re-arrangement}''.
The crux is that we get the same exponent in H\"older's inequality
$L^{2}\times L^{r}\rightarrow L^{\alpha}$  and in the Sobolev embedding
$L^{\alpha}_{1}\rightarrow L^{r}$. This is not a coincidence: it can be traced back to the scaling behaviour of the norms. 
The small-energy threshold $\epsilon$ produced by this  argument is computable.
That number might be  rather small but it follows from the further development of the theory, as in the proof of Theorem 2 above, that in Theorem 3 $\epsilon$ can be taken to be any number less than 
 the least energy of a harmonic map from $S^{2}$ to $N$. 

\subsection{The stress energy tensor and removal of point singularities}

Before beginning the proof of Theorem 4 we make a digression to review some background which will be used in the proof and also later in this article.
 
 \

Suppose that we have some functional ${\cal F}$ which depends on a Riemannian metric $g$ on a manifold $M$ and other \lq\lq fields'' (in the case at hand the fields are maps from $M$ to the fixed Riemannian manifold $N$ and the functional is the energy). The variation of ${\cal F}$ with respect to the fields, holding the metric $g$ fixed, produces Euler-Lagrange equations like the harmonic map equation. But we can also consider variations of the metric $g$, holding the fields fixed. By general principles the first variation of ${\cal F}$ can be written as
$$  \delta_{g} {\cal F}= \int_{M} (T, \delta g). $$
where the tensor $T$ is a section of $s^{2}T^{*}M$ called the stress-energy tensor, which depends on the fields and the metric. Any natural functional arising in differential geometry will be diffeomorphism invariant. It follows that if the fields satisfy the Euler-Lagrange equations generated by ${\cal F}$ then if $\delta g$ is defined by an infinitesimal diffeomorphism---i.e. $\delta g$ is the Lie derivative $L_{v}g $ of $g$ along a vector field $v$ on the manifold---then $\delta_{g}{\cal F}=0$. This is the identity ${\rm div}\ T=0$ or in index notation 
     \begin{equation}  T^{ij}_{;j}=0. \end{equation}
   When the functional ${\cal F}$ is conformally invariant the tensor $T$ is trace-free. If $v$ is a conformal Killing field on the manifold $(M,g)$ (that is,  $L_{v} g= \mu g$ for some function $\mu$ on $M$) then the contraction of $T$ by $v$ is a co-closed $1$-form. In index notation
$$   (T^{ij}v_{i})_{;j}= T^{ij}_{;j} v_{i} + T^{ij} v_{i;j}. $$
The first term on the right hand side vanishes by (14) and the second vanishes because $T^{ij}$ is symmetric and trace-free and the Lie derivative $L_{v}g$ is the symmetrisation $v_{i;j}+v_{j;i}$. 

In this section we will apply this discussion in the case of the harmonic maps energy with 2-dimensional oriented domain $M$, which can also be viewed as a Riemann surface. Taking the real part gives an isomorphism  between the  tensor square
 of $T^{*}M$, regarded
as a complex line bundle, and the trace-free symmetric tensors.  So we have a quadratic differential $\tau$ with  $T= {\rm Re} \tau$.  The equation (14) goes over to the condition that $\tau$ be a holomorphic quadratic differential; this is the 
{\it Hopf differential} defined by a harmonic map from a Riemann surface.
In a local complex coordinate $z=x+i y$
$$  \tau= \left( (f_{x}, f_{x})- (f_{y},f_{y}) + 2i (f_{x},f_{y}) \right) dz^{2} . $$
Here we are writing $f_{x}= \frac{\partial f}{\partial x}$ etc.
  Similarly,  a conformal Killing  field $v$ can be viewed as a holomorphic vector field $v$ and the contraction of $v$ with $\tau$ is a holomorphic---hence closed and co-closed---$1$-form. This completes our digression.

\

The Sacks and Uhlenbeck proof of the removal of singularities theorem goes
through a differential inequality for the energy on small discs. Let
$$  E(r)= \int_{D_{r}} \vert df\vert^{2}, $$
so clearly $E(r)$ is an increasing function of $r$ and tends to zero as $r\rightarrow
0$. We may suppose that $E(\tthreetwo)$ is less than the small energy value $\epsilon$
of Theorem 3. Applying that result to the disc of radius $\vert z\vert/2$,
say, centred at a point $z$ we get, 
\begin{equation}  \vert df (z)\vert^{2} \leq C E(3\vert z\vert/2)\  \vert z\vert^{-2},\end{equation}
so $\vert df\vert$ is $o(\vert z\vert^{-1})$. The removal of singularities theorem is proved by showing that $\vert df\vert$ is 
$O(\vert z\vert^{\delta-1})$ for some $\delta>0$. If we know this then 
$df$ is in $L^{2\alpha}$ for some $\alpha>1$ and we can apply the regularity theory discussed in the previous section to see that $f$ is smooth across the origin.

 Using (15), we see then that it suffices to show  that $E(r)$ is $O(r^{\kappa})$ for some $\kappa>0$.
The differential inequality to be established is that
\begin{equation}  \kappa E(r)\leq   r \frac{d}{dr} E(r) . \end{equation}
 If we know this then it follows by a simple comparison
argument that 
$$  E(r)\leq r^{\kappa} E(1), $$
as required.

It is convenient to exploit the conformal invariance of the problem and to work on the cylinder $(-\infty, 0]\times S^{1}$ with coordinates $(s,\theta)$.
So $r=e^{-s}$ and we now write 
$$   E(S)= \int_{s\leq S}\int f_{s}^{2}+ f_{\theta}^{2}\  d\theta ds, $$
where subscripts denote partial derivatives and we are writing $f_{s}^{2}$ for $(f_{s},f_{s})$. We know that $f$ is bounded and that the derivatives $f_{s}, f_{\theta}$ tend to zero as $s\rightarrow -\infty$. We want to show that for some $\kappa>0$
$$    \kappa E\leq \frac{dE}{dS} . $$
By translation invariance it suffices to prove this when $S=0$, that is:
\begin{equation}  \kappa \int_{s\leq 0} \int f_{s}^{2}+ f_{\theta}^{2}\  d\theta ds\leq \int_{s=0} f_{2}^{2}+ f_{\theta}^{2}\  d\theta. \end{equation}

By the  general theory reviewed above, the contraction of the Hopf differential with the Killimg field $\frac{\partial}{\partial s}$ gives the closed $1$-form
$ (f_{s}^{2}-f_{\theta}^{2})d\theta$. It follows that the integral
$$  \int_{0}^{2\pi} f_{s}^{2}-f_{\theta}^{2}\  d\theta $$
is independent of $s$ and since the integrand tends to zero as $s\rightarrow -\infty$ the integral vanishes. In other words, for each fixed $s$,
\begin{equation}  \int f_{s}^{2}\  d\theta = \int f_{\theta}^{2}\  d\theta. \end{equation}

For purposes of exposition, let us consider for a moment the case when
 $A=0$, so $f$ is an ordinary harmonic function: $\Delta f=0$. Then we have the usual integration-by-parts formula over a finite cylinder:
$$  \int_{S_{0}\leq s\leq 0} \int f_{s}^{2}+f_{\theta}^{2} \ d\theta ds= \int_{s=0} (f,  f_{s})\ d\theta -  \int_{s=S_{0}} (f, f_{s})\  d\theta .$$
Since $f_{s}$ tends to $0$ as $s\rightarrow -\infty$ and $f$ is bounded the boundary term at $s=S_{0}$ tends to zero as $S_{0}\rightarrow -\infty$ and we get
\begin{equation}  \int_{s\leq 0} \int f_{s}^{2}+f_{\theta}^{2}\  d\theta ds= \int_{s=0}
(f,  f_{s})\ d\theta. \end{equation}

 The left hand side of this formula is unchanged if we add a constant to $f$,  so we can suppose that the integral of $f$ over the boundary $\{s=0\}$ vanishes. Now for any function $g$ on the circle of integral zero we have an inequality 
\begin{equation}  \int g^{2}\  d\theta \leq \int g_{\theta}^{2}\  d\theta. \end{equation}
This is clear from the Fourier series. Thus, combining with Cauchy-Schwarz,
\begin{equation}  \left( \int_{s=0} (f,  f_{s})\  d\theta \right)^{2}\leq  \int_{s=0} f_{s}^{2}d\theta \ \int_{s=0} f_{\theta}^{2} \ d\theta . \end{equation}
Combining (19) and (21) we have

$$ \int_{s\leq 0} \int  f_{s}^{2}+f_{\theta}^{2}\  d\theta ds \leq 
\thalf \int_{s=0} f_{s}^{2}+ f_{\theta}^{2}\  d\theta, $$
which is the required inequality with $\kappa=2$.  This gives the growth rate $E(r)=O(r^{2})$ which is indeed what will occur for  smooth maps.

\

The idea now is to modify this discussion to take account of the nonlinear term $A(df,df)$ at the cost of changing the constant $\kappa$. So suppose again that $f$ is our harmonic map with $\Delta f= A(df,df)$.
 Taking the inner product of this equation with $f$ and integrating over the cylinder we  get an extra term
$$  \int_{s\leq 0}\int  ( f, A(df,df ))\  d\theta ds $$ which is bounded in modulus by the integral of $c_{1} \vert f \vert \ \vert df\vert^{2}$ for some $c_{1}$. If we knew that, over the cylinder, we have $\vert f\vert \leq \sigma c_{1}^{-1}$ for some $\sigma<1$ this would give the desired inequality 
with $\kappa= 2/(1-\sigma)$. The problem is that at this stage we know that $f$ is bounded on the cylinder but we do not know how to make this bound arbitrarily small. 
 To overcome this, take the average values
$$ F(S)=\frac{1}{2\pi} \int_{s=S} f(s,\theta) d\theta, $$ and write $g(S,\theta)= f(S,\theta)-F(S,\theta)$. 
Then we have
\begin{equation}  \int_{s\leq 0} \int (g, \Delta f)\  d\theta ds= \int_{s\leq 0}\int (dg,df)\ d\theta ds + \int_{s=0} (g, f_{s})\  d\theta , \end{equation}
since one sees as before that the other boundary term for a finite cylinder tends to $0$. Exactly the same argument as before gives
$$  
 \int_{s=0} (g,f_{s})\leq (1/2) \int_{s=0} f_{s}^{2}+f_{\theta}^{2} \ d\theta $$
We have $(dg,df)= \vert df\vert^{2}- (f_{s}, F_{s})$ so we get
$$  \int_{s\leq 0}\int f_{s}^{2}+ f_{\theta}^{2} \  d\theta ds \leq
(1/2)  \int_{s=0} f_{s}^{2}+f_{\theta}^{2}\  d\theta + I + II, $$
where
$$ I= \int_{s\leq 0}\int (g,\Delta f)\  d\theta ds, $$
and
$$ II= \int_{s\leq 0}\int (f_{s}, F_{s})\  d\theta ds. $$

The integrand in I is bounded by $c_{1} \vert g \vert \vert df\vert^{2}$.
For each fixed $s$ the integral of $g$ is zero and there is an inequality in the same vein as (20)
$$  \vert g\vert^{2} \leq c_{2} \int g_{\theta}^{2}\  d\theta = c_{2} \int f_{\theta}^{2}\  d\theta. $$
(This is the Sobolev embedding $L^{2}_{1}\rightarrow C^{0}$ in dimension $1$.)
So we can suppose that $\vert g\vert$ is as small as we please; say $\vert g\vert \leq \sigma c_{1}^{-1}$ for $\sigma<\thalf$. Then the term I is bounded by
$\sigma$ times the energy.

Turning to the term II: consider the integral over $\theta$ for fixed $s=S$. This is
$$  (2\pi)^{-1} \vert \int_{s=S} f_{s}\ d\theta \vert^{2}, $$ which is bounded by $$  \int_{s=S} f_{s}^{2} d\theta = (1/2) \int_{s=S} f_{s}^{2}+ f_{\theta}^{2} \ d\theta, $$
using (20) again. Putting things together we get
$$  (\thalf-\sigma) \int_{s\leq 0}\int f_{s}^{2}+ f_{\theta}^{2}\   d\theta ds \leq  (1/2) \int_{s=0} f_{s}^{2}+f_{\theta}^{2}d\theta, $$
which is the desired inequality with $\kappa= 1-2\sigma$.

\

In the case of pseudoholomorphic curves the proof of the differential inequality is simpler, using the fact that the energy is twice the integral of $f^{*}(\omega)$, where $\omega$ is the symplectic form. If one can write $\omega=d\alpha$ for a $1$-form $\alpha$ over a neighbourhood of the image of $f$ then Stokes' theorem expresses the energy as a boundary integral. (We will use an argument like this in the proof of Theorem 12 below.)

\section{Harmonic maps in higher dimensions}

\subsection{Monotonicity of normalised energy}

The main topic of this Section 4 is Uhlenbeck's work with Schoen in the 
 paper \cite{kn:SchU1},  on weak solutions to the harmonic map equation, but we begin in this subsection with some background and  results for smooth  maps.

Consider again the variational theory on an $n$-dimensional manifold $M$ of a functional ${\cal F}(g,\Phi)$ given by the integral of an $n$-form $F(g,\Phi)$. Suppose that this functional has a non-zero scaling weight $w$ under conformal change of the metric in that
$$ F(\lambda g, \Phi)= \lambda^{w} F(g,\Phi). $$
It follows that the trace of the energy momentum tensor is $ w L(g,\Phi) {\rm vol}_{g}$. Let $v$ be a conformal vector field, so 
$ v_{i;j}+ v_{j;i}= 2 \mu\  g_{ij}$ for some function $\mu$. Then we have
$$  (v_{i} T^{ij})_{j}= (v_{i;j} T_{ij})= w \mu F(g,\Phi). $$
So we can write
   $$  w \mu F(g,\Phi)=  d \eta, $$
   where $\eta$ is the $(n-1)$-form  $*(v_{i} T^{ij})$.
   Thus if $U$ is a domain in $M$ with compact closure and smooth boundary we have
   \begin{equation} w \int_{U} \mu\ F(g,\Phi) =  \int_{\partial U}  \eta.  \end{equation}

Apply this discussion to the harmonic maps energy functional on  a manifold $M$ of dimension $n>2$, so the field $\Phi$ is $f:M\rightarrow N$ and
$F(g,f)= \vert df\vert^{2} {\rm vol}_{g}$ which has a conformal weight
 $w=n/2-1$.
Suppose that $M=\bR^{n}$, so we have a conformal vector field $r\frac{\partial}{\partial r}$ as before and the function $\mu$ is the constant $1$. The  stress-energy tensor is 
$$  T_{ij}= (\nabla_{i}f,\nabla_{j}f)- \thalf \vert df\vert^{2} \delta_{ij}$$
One sees then that the restriction of the $(n-1)$-form $\eta$ to the unit sphere is
$$  \eta\vert_{S^{n-1}} = \left(\thalf \vert df\vert^{2} - \vert \nabla_{r} f\vert^{2}\right) d{\rm vol}_{S^{n-1}}, $$
where $\nabla_{r}$ is the radial derivative. So the identity (23)  is
\begin{equation} (n-2) \int_{B} \vert df\vert^{2} =  \int_{\partial B}  \vert df\vert^{2}-2  \vert \nabla_{r} f\vert^{2}, \end{equation}
 (which agrees with (18) in the case $n=2$). Thus
$$  (n-2)\int_{B} \vert df\vert^{2} \leq  \int_{\partial B} \vert df\vert^{2}  $$ with equality if and only if $\nabla_{r} f=0$ on $\partial B$.
If we apply the same argument to the ball $B_{r}$ of radius $r$ we get
$$  (n-2) \int_{B_{r}} \vert df\vert^{2} \leq  r\int_{\partial B} \vert df\vert^{2} $$ 
Let $\hatE(r)$ be the {\it normalised energy}
$$  \hatE(r)= \frac{1}{r^{n-2}} \int_{B_{r}} \vert df\vert^{2}. $$
The inequality above is equivalent to the {\it monotonicity condition} $d\hatE/dr\geq 0$ and in fact for $r_{1}<r_{2}$
$$    \hatE(r_{1})= \hatE(r_{2})- 2 \int_{r_{1}<r<r_{2}} r^{2}
\vert \nabla_{r} f\vert^{2} . $$

\

A good way to think about the normalised energy is through rescaling. Map the unit ball $B$ to the  ball $B_{r}$ by $x\mapsto r x$ and compose with the restriction of  $f$ to $B_{r}$ to get $\tilde{f}:B\rightarrow N$. Then the normalised energy of $f$ on $B_{r}$ is the energy of $\tilde{f}$ on the unit ball $B$. The monotonicity condition says that  the map $f$ \lq\lq looks  better''---in the sense of having smaller energy---if we look at it on smaller and smaller scales in this manner. Another useful observation is that if $f$ is the composite of a map $\underline{f}$  from $\bR^{2}$ to $N$ with an orthogonal projection from $\bR^{n}$ to $\bR^{2}$ then the normalised energy for $f$ agrees with the ordinary energy for $\underline{f}$,  up to a factor.

\

One important consequence of 
the monotonicity property is a small energy result. The statement is essentially the same as in the 2-dimensional case of the previous section. A difference is that the equations now depend essentially on the Riemannian metric on $M$ so we formulate the statement a bit differently. 
\begin{prop}

Let $M,N$ be compact Riemannian manifolds. There are $\epsilon, r_{0}, C>0$ such that if $B_{x,r}$ is a metric $r$-ball in $M$ with $r\leq r_{0}$ and the normalised energy $\hatE(B_{r})$ of $f$ on $B_{r}$ is less than $\epsilon$ then on the half-sized ball $B_{x,r/2}$ we have
$$\vert df\vert^{2}\leq C r^{-2} E(B_{r}).$$
\end{prop}
As before, once we have the $L^{\infty}$ bound on $df$ we can go on to get estimates on all higher derivatives.
This small energy statement can be proved in  a manner similar to the proof above of Theorem 3 but using Morrey spaces  (which we will encounter in  6.3 below) in place of $L^{p}$ spaces. We will  discuss here the proof of   Schoen in \cite{kn:Schoen}, using a \lq\lq worst point'' argument,  which could also be used in the 2-dimensional case of Theorem 3.

For simplicity we suppose that $M$ is locally Euclidean, so the discussion above applies to give a monotonicity formula for balls of sufficiently small size.  For a general Riemannian manifold $M$ we do not have exact formulae for the normalised energy but the equations hold with extra error terms which can be made as small as we please, since the geometry is close to Euclidean on small scales, and  the same arguments work with minor modifications. By scale invariance and adjustment  of constants we can suppose that the ball  $B_{x,r}$ is the unit ball $B$ in $\bR^{n}$ and that $f$ has normalised energy  at most $E\leq \epsilon$ on any interior  ball. 

For a point $x\in B$ let $D(x)$ be the distance to the boundary of $B$, i.e.
$D(x)= 1-\vert x\vert$. The idea is to consider the quantity 
$$ M= {\rm max}_{x\in B}  D(x) \vert df (x)\vert .$$ 
 Since the function $D$ vanishes on the boundary, the maximum is achieved at some interior point $x_{0}$. For $\rho\leq 1$, rescale the ball of radius $\rho D(x_{0})/2$ with centre $x_{0}$ to unit size to get a map $\tilde{f}_{\rho}$ on the unit ball $B$ with the properties
\begin{itemize}
\item $\int_{B}\vert d\tilde{f}_{\rho}\vert^{2}\leq E $;
\item $\vert d\tilde{f}_{\rho}\vert\leq 4M\rho$ on $B$;
\item $\vert d\tilde{f}_{\rho}(0)= 2M\rho$;
\end{itemize}
where the first item uses the small energy property of $f$ on the interior ball. Now we have
$$   \vert \Delta \tilde{f}_{\rho}\vert \leq c \vert d\tilde{f}_{\rho}\vert^{2}. $$

Elliptic theory gives an inequality 

$$    \vert d\tilde{f}_{\rho}(0)\vert \leq c\left( \Vert \Delta\tilde{f}_{\rho}\Vert_{L^{\infty}} + \Vert d\tilde{f}_{\rho}\Vert_{L^{2}} \right), $$
so we get
$$  \vert d\tilde{f}_{\rho}(0)\vert \leq c_{1}\rho^{2} M^{2}+ c_{2}
 \sqrt{E}, $$
 and $$\rho M \leq c_{3} \rho^{2} M^{2} + c_{4} \sqrt{\epsilon}. $$
 If we choose $\epsilon$ small enough the equation
    $$  y = c_{3} y^{2} + c_{4} \sqrt{E} $$
    will have two solutions: a small solution $y_{0}$,  approximately  $c_{2}\sqrt{E}$, and a large solution $y_{1}$, approximately $c_{1}^{-1}$. Fix  such an $\epsilon$. Then the inequality implies that either $\rho M\leq y_{0}$ or $\rho M\geq y_{1}$. For very small $\rho$ the first alternative must hold and by continuity it must continue to hold for all $\rho\leq 1$. So we conclude that $M\leq {\rm const.} \sqrt{E}$,  which establishes the Proposition.

Another  result about smooth harmonic maps  that can be proved using a similar approach is:
\begin{prop}

There is a constant $M_{0}$ such that if  the harmonic map $f$ satisfies a H\"older bound on the unit ball $B$:
 $$  \vert f(x)- f(y)\vert \leq  \vert x-y \vert^{\alpha} $$
 then $\vert df(x)\vert \leq M_{0}  D(x)^{-1} $, where $D(x)$ is the distance to the boundary, as above.
 \end{prop}
 
 As usual, we can go on to get estimates on all derivatives of $f$ in the interior depending only on the  H\"older bound.
 To prove this Proposition we define $M$ and $x_{0}$ as above. If $M\leq 2$ we can take $M_{0}=2$. If $M\geq 2$ rescale the ball of radius $M^{-1}$ centered at $x_{0}$ to unit size to get a harmonic map $\tilde{f}$on the unit ball $B$ with 
\begin{enumerate}
\item $\vert \tilde{f}(x)-\tilde{f}(y)\vert\leq M^{-\alpha} \vert x-y\vert^{\alpha}$;
\item $\vert d\tilde{f}\vert \leq 2$ on $B$;
\item $ \vert d\tilde{f}(0)\vert =1$.
\end{enumerate}
The harmonic map equation and item (2) give a bound on
 $\vert \Delta \tilde{f}\vert$ over $B$, which gives a $C^{,\alpha}$ bound on $d\tilde{f}$ in the half-sized ball. It then follows from item (3) that for some computable number $\kappa$ we can choose a ray $\{t \nu\}$ through the origin such that $\vert \tilde{f}(t\nu)-\tilde{f}(0)\vert\geq t/2$ for $t\leq \kappa$. Then item (1) implies that $M\leq 2^{\alpha} \kappa^{1-1/\alpha}$.

\

 The small energy result gives a partial compactness property, extending what we have seen for surfaces in Section 3. Let $M$ and $N$ be compact and let $f_{i}:M\rightarrow N$ be a sequence of harmonic maps with energy bounded by a fixed constant $E_{\rm max}$. As in Section 2 we can suppose that
 the energy densities $\vert df_{i}\vert^{2}$ converge as Radon measures.  Now we define
  $$\mu_{i}(x,r)= r^{2-n} \int \chi_{x,r} \vert df_{i}\vert^{2}, $$
  and $\mu(x,r)=\lim_{i\rightarrow \infty} \mu_{i}(x,r)$.
We assume that the $\mu_{i}(x,r)$ are  increasing functions of $r$. If $M$ is locally Euclidean this follows from the monotonicity property and in general it will be true up to an unimportant error term.   Then $\mu(x,r)$ is
 increasing  and   has a limit $\mu(x)$ as $r$ tends to $0$. As  before, we define the set
$S\subset M$ to be the set of points where $\mu(x)\geq \epsilon/2$, for the constant $\epsilon$ in the small energy result,  Proposition 4. Just as before we get, after passing to a subsequence, a limiting harmonic map $f:M\setminus S \rightarrow N$ and the $f_{i}$ converge to $f$ in $C^{\infty}$ on compact subsets of
$M\setminus S$. 
   Given a small $\delta>0$, choose a maximal collection of  disjoint $\delta/2$ balls $\{ B_{\alpha}\}$ centred at points $x_{\alpha}$ of $S$. Let $A$ be the number of balls. Then the $\delta$-balls with the same centres cover $S$ and
    $$  \int_{B_{\alpha}} \vert df_{i}\vert^{2}
     \geq (\delta/4)^{n-2}\mu_{i}(x_{\alpha},\delta/4). $$
     Since the $B_{\alpha}$ are disjoint we have
     $$ E_{\rm max} \geq \sum_{\alpha} (\delta/4)^{n-2}\mu_{i}(x_{\alpha},\delta/4),  $$
and taking the limit we can replace $\mu_{i}(x_{\alpha,\delta/4})$ by
  $\mu(x_{\alpha,\delta/4})$. But $\mu(x_{\alpha},\delta/4)\geq \mu(x_{\alpha})\geq \epsilon/2$ so we get a bound on the number $A$ of balls $B_{\alpha}$
$$  A \leq  C \delta^{2-n} $$
 with $C= 2^{2n-3} \epsilon^{-1} E_{\rm max} $. So  the set $S$ is covered by at most $C'\delta^{2-n}$ balls of radius $\delta$. This implies that $S$ has Hausdorff dimension  at most $(n-2)$ and  the $(n-2)$-dimensional Hausdorff measure is bounded by $C$.  
(In fact we get a stronger statement,  that the $(n-2)$-dimensional \lq\lq Minkowski content'' is finite. This is because  the balls in our cover have the same  radius $\delta$: in the definition of Hausdorff measure one is allowed to cover by balls of varying radii.) 

\

In sum we have:
\begin{prop}
For $M,N$ compact a sequence of harmonic maps from $M$ to $N$ with a fixed energy bound has a subsequence which converges off a set of Hausdorff codimension at least $2$.
\end{prop}

  \subsection{Minimising maps}
  
In the previous subsection we have considered {\it smooth} harmonic maps. The thrust of the Schoen and Uhlenbeck paper \cite{kn:SchU1} is different because they consider the much more formidable case of a class of weak solutions in $L^{2}_{1}$.
These  can have singularities, and the great achievement of Schoen and Uhlenbeck was to make these singularities
 somewhat tractable.

More precisely, for $N$ embedded in the Euclidean space $V$, let $L^{2}_{1}(M,N)$ be the set of maps $f:M\rightarrow V$ which are in $L^{2}_{1}$ in the usual sense and which map almost all points of $M$ to $N$. (Of course, for ${\rm dim}\  M>1$, such a map $f$ is only  defined almost everywhere.) The energy functional is defined on these maps and $f$ is called a weak solution if the first variation of the energy vanishes, which is equivalent to $f$ being a weak solution of the equation (9). This notion makes sense because the  nonlinear term $A_{f}(df,df)$ is in $L^{1}$. These weak solutions can be bizarre: there are examples which are not continuous at any point of $M$ \cite{kn:Riv}.
Schoen and Uhlenbeck showed that if one   restricts to the class of 
{\it energy minimising} maps the situation is much better. For our discussion,  we could take the definition of energy minimising to be that there is some $\rho>0$ such that for all balls $B_{\rho}\subset M$ if $g\in L^{2}_{1}(M,N)$ is equal to $f$ outside $B_{\rho}$ then $E(g)\geq E(f)$. (Any smooth harmonic map is energy-minimising in this sense.)  In fact,  in \cite{kn:SchU1} Schoen and Uhlenbeck consider a more general class of equations, adding a perturbation term to the energy, and in \cite{kn:SchU2} they extend the theory to the Dirichlet problem on a manifold
 $M$ with boundary.

The foundation of the Schoen-Uhlenbeck work is to establish versions of monotonicity and the small energy property for energy-minimising harmonic maps. In the end the statements are essentially the same as for the smooth case but the proofs are  different  because many of the constructions discussed above do not make sense in this wider class. 

For our discussion 
 we assume that $M$ is locally Euclidean, so by scaling we may regard  $f$ as being defined on the unit ball $B$ in $\bR^{n}$ and  we can assume that any variation supported in $B$ increases energy.
The Schoen-Uhlenbeck proof of monotonicity proceeds as follows. Suppose that the restriction of $f$ to the unit sphere $S^{n-1}=\partial B$ is also in $L^{2}_{1}$ and  define a map
$g$ to be equal to $f$ outside $B$ and by
$g(x)= f(\frac{x}{\vert x\vert})$ 
for $\vert x\vert\leq 1$. Of course $g$ is not defined at the origin but   when $n>2$ it is an $L^{2}_{1}$ map. Simple calculus gives
  $$    \int_{B}\vert dg\vert^{2} = \int_{S^{n-1}} \vert d_{S^{n-1}} g\vert^{2} \int_{0}^{1} r^{n-3} dr = (n-2)^{-1} \int_{S^{n-1}} \vert  d_{S^{n-1}}g\vert^{2} . $$
 Here the notation $d_{S^{n-1}}$ refers to the derivative of the restriction of the map to the sphere. The energy-minimising property gives
$$ (n-2) \int_{B}\vert df\vert^{2} \leq \int_{S^{n-1}}\vert df\vert^{2}- \int_{S^{n-1}} \vert \nabla_{r} f \vert^{2}. $$
(Notice that the term involving the radial derivative enters with a different factor compared with (24), but this will not matter.)
Let $E(r)$ be the energy of $f$ on the ball $B_{r}$.
Then $E(r)$ is an increasing function of $r$ and so differentiable at almost all $r$.  Similarly for almost all $r$ the restriction of $\vert df\vert^{2} $ to the boundary of the ball $B_{r}$ is in $L^{2}$ and at such values of $r$ 
  $$  E'(r)= \int_{\partial B_{r}} \vert df\vert^{2}. $$
 For such $r$ we can apply the preceding discussion for the unit ball, after rescaling, and we obtain the inequality
$$ (n-2) E(r) \leq  r E'(r). $$
  As before, this is the monotonicity statement that the normalised energy $\hatE(r)=r^{2-n}E(r)$ is increasing.  Moreover we have, for $r_{1}<r_{2}$,
 \begin{equation} \hatE(r_{1})\leq \hatE(r_{2})- \int_{r_{1}<r<r_{2}} r^{2} \vert \nabla_{r} f\vert^{2} \end{equation}

 \

The work of Schoen and Uhlenbeck develops an important analogy between the theories of harmonic maps and of minimal submanifolds, and more general volume-minimisig sets.  The analogue of the argument above in the latter case, for a $d$-dimensional volume-minimising set $X\subset \bR^{m}$ and a point $x$ in $X$, is to consider, for small $r$, the intersection
$Y_{r}$ of  $X$ with the sphere of radius $r$ centred at $x$. Let $CY_{r}$ be the cone over $Y_{r}$ with vertex at $x$ and let$\tilde{X}_{r}$ be the set obtained by removing the intersection $X\cap B_{x,r}$ from $X$ and replacing it with the cone $CY_{r}$. Comparing $\tilde{X}_{r}$ with $X$, the volume-minimising property shows that ${\rm Vol}( X\cap B_{x,r}) \leq {\rm Vol}\ CY_{r}$ and this leads to the monotonicity of the normalised volumes $r^{-d} {\rm Vol}(B_{x,r}\cap X)$ with respect to $r$, for fixed $x$.

The Schoen-Uhlenbeck proof of the small energy result for energy-minimising maps is more involved and we will outline  it in subsection 4.3 below. But before that we discuss the general structural results and  overall picture which Schoen and Uhlenbeck obtained. The first consequence is that a minimising map $f$ is smooth outside a closed singular set $\Sigma$ with ${\rm dim}\ \Sigma <n-2$. For example when $n=2$ this says that $\Sigma$ is empty, which is immediate from  the small energy result. For $n>2$ it is proved by a covering argument similar to the one  we described above for Proposition 6.  A refined result is
\begin{thm}
\begin{itemize}
\item ${\rm dim}\ \Sigma\leq n-3$.
\item Suppose that  for some $k\geq 2$ and for all $\nu$ with $2\leq\nu\leq k$ there is no smooth minimising harmonic map from the sphere $S^{\nu}$ to $N$. Then
${\rm dim}\ \Sigma \leq n-k-2$.
\end{itemize}
\end{thm}
In the second item here, if $k=n-1$ then the statement is that $\Sigma$ is empty, so the map is smooth.

The proofs of these refined results depend on the important notion of a \lq\lq tangent map'', introduced by Schoen and Uhlenbeck.  This is  analogous to  the notion of a tangent cone in submanifold geometry. Suppose for simplicity that the domain $M$ of the  minimising map $f$ is the unit ball in $\bR^{n}$ and for $\lambda>1$
let $f^{\lambda}$ be the composite of $f$ with the scaling map; so 
$f^{\lambda} :\lambda B \rightarrow N$. Take any sequence $\lambda_{i}\rightarrow \infty$, so for any compact set $K\subset \bR^{n}$ the map $f^{\lambda_{i}}$ is defined over $K$ for large enough $i$. Monotonicity implies that the energies of the $f^{\lambda_{i}}$ are bounded on compact sets so, possibly passing to  a subsequence, there is a weak $L^{2}_{1, {\rm loc}}$ limit $f_{\infty}$.  The inequality (25)  implies that $f_{\infty}$ is  {\it radially homogeneous} in that $\nabla_{r} f_{\infty}=0$ almost everywhere. Such a map is called a {\it tangent map} to $f$ at $0$.
 
 Uniqueness of the tangent map, {\it i.e.} that one gets the same limit for any sequence of scalings, is a major question in general---there are examples  \cite{kn:White} where it is not unique--- but is not
 directly relevant to the discussion here. The main difficulties are  that the convergence  obtained is only in the weak topology and that it is not clear that the limit will again be minimising. Much of the work in the paper of Schoen and Uhlenbeck goes into overcoming these difficulties. They show that the convergence can be improved to $L^{2}_{1,{\rm loc}}$  and that, at least in certain restricted situations, the limit is minimising. Glossing over many details, we illustrate the argument for the case when $n=3$. Then the statement (1) of Theorem 5 can be improved to the statement that $\Sigma$ is a discrete set. To see this, let $p$ be a singular point and consider a tangent map at $p$. By radial homogeneity, this is equivalent to a map
$g:S^{2}\rightarrow N$. Schoen and Uhlenbeck show that $g$ is minimising, so the singular set is empty by the previous discussion in dimension  $n=2$. Thus the tangent map $f_{\infty}$ has an isolated singularity at the origin and this implies that the singularity $p$ of $f$ is isolated.  For another illustration of the argument, consider the case when there are no minimising harmonic maps
$S^{\nu}\rightarrow N$ for $2\leq\nu\leq n-1$. Suppose there were a singular point $p$ of $f$. From the tangent map we get a map $g:S^{n-1} \rightarrow N$. Assume for  simplicity that  that this is minimising. The hypothesis implies that $g$ cannot be  smooth,  so we can go a singular point of $g$ in $S^{n-1}$ and take a second tangent map there. This gives a map from $S^{n-2}$ to $N$ and the hypothesis implies that this cannot be smooth, so we can a tangent map at a singular point. Continuing in this way, with these iterated tangent maps, we get a contradiction to the existence of $p$, so the original map is smooth.

\subsection{Small energy}

In this subsection we outline the Schoen and Uhlenbeck proof of the small energy result.

  First, it was established before that a H\"older continuous weak harmonic map is smooth \cite{kn:Hild}. (This is similar to Proposition 5 in that {\it a posteriori} the estimate in that Proposition holds, but of course the proof is much harder.) Next, it suffices to get a bound for the growth of the normalised energy function. Morrey's Lemma states that if a function $g$ on $\bR^{n}$
with weak derivative in $L^{1}$ satisfies an estimate:
$$  r^{-n} \int_{B_{x,r}} \vert dg\vert \leq C r^{-\beta}, $$
for all balls $B_{x,r}$ then $g$ is in $C^{,1-\beta}$. For our map $f$
$$  \left( \int_{B_{x,r}}\vert df\vert\right)^{2} \leq \int_{B_{r,x}} \vert df\vert^{2}\  {\rm Vol}\  B_{x,r}, $$
so if the normalised energy on all $B_{x,r}$ is less than $C r^{\alpha}$ then  $f$ is in $C^{,\alpha/2}$. 

The essential statement in Schoen and Uhlenbeck's proof of the small energy result  is then: 
\begin{thm}
There is an $\epsilon_{0}>0$ and $\theta_{0}\in (0,1)$ such that  if $f:B\rightarrow
N$ is an energy minimising map with normalised energy less than $\epsilon_{0}$ on all interior balls then
$$   \hatE(\theta_{0})\leq \thalf \hatE(1)= \thalf E $$
\end{thm}
(The factor $\thalf$ here could be replaced by any fixed number in $(0,1)$.)
Given this, it follows from an elementary argument (similar to that in subsection 2.3 above)
that $\hatE(r)\leq C r^{\alpha}$ for a suitable $\alpha$, depending on $\theta$.

The proof by Schoen and Uhlenbeck of Theorem 6 involves the choice of four  parameters
$\theta_{0}, \theta, \tau, h$. Here $\theta_{0}$ will be as in the statement of the Theorem, and we will choose $\theta_{0}<\teighth$ say. The parameter $\theta$ will be chosen in the interval $[\theta_{0}, 2\theta_{0}]$.
The parameter $\tau$ will be  much smaller than $\theta_{0}$. Given $\theta,\tau$ we write $A$ for the annulus 
$A= \{ x:\theta\leq \vert x\vert \leq \theta+\tau\}$. The idea is to construct a comparison map $\tilde{f}:B\rightarrow N$ such that $\tilde{f}(x)=f(x)$ for $\vert x\vert\geq \theta+\tau$. Then the minimising property of $f$ gives
\begin{equation} \int_{B_{\theta}} \vert df\vert^{2}\leq I+II\end{equation}
where
$$  I=\int_{B_{\theta}} \vert d\tilde{f}\vert^{2}\ \ \ ,\ \ \ 
  II= \int_{A}\vert d\tilde{f}\vert^{2}. $$
and the task will be to bound these terms $I$ and $II$.

We write $E$ for the energy of $f$ on the unit ball $B$, so $E\leq \epsilon_{0}$.
and we use the convention that $c$ is a constant which changes from line to line. We also fix a tubular neighbourhood $\Omega\subset V$ of $N$ in $V$ and let $\pi:\Omega\rightarrow N$ be the standard projection.

 The construction of $\tilde{f}$ goes through three other maps $f_{1}, f_{2}, f_{3}$ and depends on the parameter $h>0$, which is a
 small smoothing parameter.

{\bf The map $f_{1}$}

 The map $f_{1}:B\rightarrow V$ is the smooth map obtained as a standard mollification of $f$ by convolution with a function supported in the ball of radius $h$. (More precisely, $f_{1}$ will be defined on a slightly smaller ball than $B_{1}$ but this does not matter.) Thus if $\phi$ is  positive function on $V$ with integral $1$ and supported in the unit ball:
$$   f_{1}(x)= h^{-n} \int \phi( h^{-1} (x-y)) f(y) dy= \int \phi(z) f(x- hz) dz, $$
where the two formulae are related by the change of variable $y=x- hz$.

{\bf The map $f_{2}$}

The map $f_{1}$ does not map into $N$ and the distance between $f_{1}(x)$
and $f(x)$ need not be small {\it a priori}.  But we will see that we can arrange, by choosing $\epsilon_{0}$ small, that $f_{1}$ maps into the tubular neighbourhood $\Omega$. Then we define $f_{2}= \pi \circ f_{1}: B\rightarrow N$.

{\bf The map $f_{3}$}

We modify the smoothing construction of $f_{1}$ to make a new map
$f_{3}:B\rightarrow V$ equal to $f$ outside $B_{\theta+\tau}$ and to $f_{2}$ in $B_{\theta}$. Let $\eta$ be a function on $B$ equal to the constant $h$ on $B_{\theta}$ and to $0$ outside $B_{\theta+\tau}$ and set
$$  f_{3}(x)= \int \phi(z) f(x-\eta(x) z) dz . $$
We will arrange that $f_{3}$  maps into the neighbourhood $\Omega$ and  we define
$\tilde{f}= \pi\circ f_{3}$.

\

{\bf Bounds on the convolution}

The first business is to arrange that $f_{1}$ and $f_{3}$ map into the tubular neighbourhood $\Omega$ of $N$. This is a crucial insight in the Schoen and Uhlenbeck proof and depends on  the following lemma.
\begin{lem}
For $\phi$ as above there is a constant $C$ such that  all functions
$g$ on the unit ball $B$ with
$$  \int_{B} \phi(y) g(y)  =0$$
satisfy $\Vert g\Vert_{L^{2}}\leq C \Vert dg \Vert_{L^{2}}$.
\end{lem}
If $\phi$ were replaced by a constant this becomes the Poincar\'e Lemma. The proof of the generalisation is straightforward. An immediate consequence is that for any function $g$ on $B$
$$  \Vert g - g_{*}\Vert_{L^{2}} \leq C \Vert dg\Vert_{L^{2}}$$
where $g_{*}$ is the constant function equal to
$$  g_{*}= \int_{B} \phi(y) g(y) dy . $$
In particular there is some $y$ in $B$ such that 
$$ \vert g(y)-g_{*}\vert \leq \frac{C}{\sqrt{{\rm Vol} B}} \Vert dg\Vert_{L^{2}}. $$
Now consider $x_{0}\in B$ and apply this to $g(y)= f(x_{0}- h y)$ 
so $g_{*}=f_{1}(x_{0})$.
We have
$$  \Vert dg\Vert^{2}_{L^{2}}= h^{2-n}\int_{\vert x-x_{0}\vert \leq h}\vert df\vert^{2} \leq \epsilon_{0}. $$
 We deduce that there is some $x$ with $\vert x-x_{0}\vert\leq h$ such that
 $\vert f(x)- f_{1}(x_{0})\vert \leq c\ \epsilon_{0}^{1/2}$. So $f_{1}$ maps into the $c\ \epsilon_{0}^{1/2}$ neighbourhood of $N$. This argument also applies to the map $f_{3}$, because $x_{0}$ is fixed.

 The standard convolution formula shows that $\Vert df_{1}\Vert^{2}_{L^{2}}\leq \Vert df\Vert^{2}_{L^{2}}=  E$. We also have a pointwise bound
\begin{equation}  \vert df_{1}\vert \leq c \epsilon_{0}^{1/2} h^{-1}, \end{equation}
which follows easily from the bound on the normalised energy. 

\

{\bf The term I}

By construction 
$$  I= \int_{B_{\theta}} \vert df_{2}\vert^{2}. $$

Composition with the projection $\pi$ can increase the norm of the derivative by at most a small factor so it suffices to bound
 the integral of $\vert df_{1}\vert^{2}$ over $B_{\theta}$. We need a better bound than that given by (27). To achieve this, Schoen and Uhlenbeck consider the harmonic function $v$ on the ball $B_{\thalf}$ with the same boundary values as $f_{1}$. Then $v$ minimises the Dirichlet energy over all functions with these boundary values so
$$  \int_{B_{\thalf}} \vert dv\vert^{2}\leq \int_{B_{\thalf}} \vert df_{1}\vert^{2}\leq  E. $$
Standard theory promotes this to a pointwise bound on the interior ball $B_{\theta}\subset B_{\tquart}\subset B_{\thalf}$ so
\begin{equation}  \int_{B_{\theta}} \vert dv\vert^{2} \leq c E \ {\rm Vol} \ B_{\theta}= c E \ \theta^{n}. \end{equation}
Now write $w=f_{1}- v$, so $\Delta w =\Delta f_{1}$. Recall that $f_{1}$ is the convolution $\phi_{h}* f$ of f with a function $\phi_{h}$ of $L^{1}$ norm $1$.
 The Laplace operator $\Delta$ on $\bR^{n}$ commutes with convolution so $$\Delta f_{1}= \phi_{h}*(\Delta f)= \phi_{h}*(A(df,df)),$$ 
and hence
$$  \Vert \Delta f_{1} \Vert_{L^{1}}\leq \Vert \phi_{h}\Vert_{L^{1}} \Vert A(df,df)\Vert_{L^{1}} =  \Vert
A(df,df)\Vert_{L^{1}}. $$ 
 Clearly $\Vert A(df,df)\Vert_{L^{1}}\leq c E$ so we get
 $$  \Vert \Delta w\Vert_{L^{1}}\leq c E. $$
 Hence
 $$  \int_{B_{\thalf}}\vert dw\vert^{2}= \int_{B_{\thalf}} (w,\Delta w)\leq c E \ {\rm sup}_{B_{1/2}}\  \vert w\vert. $$
The bound (27) and maximum principle considerations imply that $\vert w\vert \leq c\epsilon_{0}^{1/2} h^{-1}$ on $B_{\thalf}$. So we conclude that
$$\int_{B_{\theta}} \vert dw\vert^{2}\leq\int_{B_{\thalf}} \vert dw\vert^{2}\leq c \ \epsilon_{0} h^{-1} E. $$ Combined with (28) this gives a bound on the $L^{2}$ norm of $df_{1}$ and hence of $\tilde{f}$ over $B_{\theta}$:
\begin{equation} I= \int_{B_{\theta}} \vert d\tilde{f}\vert^{2}\leq c( \epsilon_{0} h^{-1}  + \theta^{n}) E. \end{equation}
 
 {\bf The term II}
 
 As before, it suffices to work with $f_{3}$. Recall that $A$ is the annulus $\theta\leq \vert x\vert\leq \theta+\tau$.
 Let $A_{+}$ be the $h$ neighbourhood of $A$, so the maps $f_{3}$ and $\tilde{f}$ on $A$ are determined by the restriction of $f$ to $A_{+}$. 
\begin{lem}
Suppose that the derivative of $\eta$ is bounded by $\vert d\eta\vert\leq k$. Then
$$  \int_{A}\vert d f_{3}\vert^{2}\leq c(1+k)^{2} \int_{A_{+}}\vert df\vert^{2}. $$
\end{lem}
We have
$$  f_{3}(x)=\int \phi(y) f(x- \eta(x) y) dy. $$
When we differentiate with respect to $x$ we get a term from the derivative of $\eta$, which is bounded by $k$. This gives 
$$  \vert d f_{3} (x)\vert \leq  (1+k)\int\ \phi(y)\   \vert df\vert(x-\eta(x)y) \ dy. $$
A quick route from here is to use the theory of the maximal function. For points $x$ where $\eta(x)>0$
$$  \vert df_{3}(x)\vert \leq (1+k) \int_{\vert y\vert \leq 1} \vert df\vert(x-\eta(x)y)\ 
dy\leq (1+k)  \eta(x)^{-n} \int_{B_{x,\eta(x)}} \vert df\vert\leq (1+k) M(\vert df\vert),
$$ where $M(\vert df\vert)$ is the maximal function of $\vert df\vert$.
Then 
$$  \Vert df_{3}\Vert_{L^{2}}\leq c(1+k) \Vert M(df)\Vert_{L^{2}}\leq c (1+k)
\Vert df\Vert_{L^{2}}, $$
and the Lemma follows. (Schoen and Uhlenbeck give a 
direct calculus proof of this Lemma.)

\

As before, the projection $\pi$ only changes the energy by  small amount, so Lemma 2 implies that, if we choose $\eta$ so that $\vert d\eta\vert\leq 1$, we have
$$  II= \int_{A} \vert d \tilde{f}\vert^{2}\leq c \int_{A^{+}} \vert df\vert^{2}. $$
 The right hand side here is bounded by $c E$ but this does not suffice since the constant $c$ could be large. To overcome this Schoen and Uhlenbeck bring in another idea: the choice of $\theta$. The condition that $\vert d\eta\vert
\leq 1$ implies that $h$ cannot be more than the annulus thickness $\tau$. Let us now fix $h=\tau/2$ say. So $A_{+}$ is an annulus of thickness $5\tau$.
Recall that
 $\tau$ is to be much smaller than $\theta_{0}$. There are approximately
 $Q= \theta_{0}/ 5\tau$ disjoint annuli of the form $A_{+}$ for different values of $\theta$ in $[\theta_{0}, 2\theta_{0}]$. So we can make a choice of one of these such that
$$  \int_{A_{+}}\vert df\vert^{2}\leq E/Q\leq c E \tau/\theta_{0} . $$

Making this choice, combining with the bound (29) for the term I and setting $h=\tau/2$ we get 
$$  I+II\leq c( \epsilon_{0}/ \tau + \theta_{0}^{n}+ \tau/\theta_{0}) E, $$
(By adjusting constants it does not matter if we write $\theta$ or $\theta_{0}$ here since $\theta_{0}<\theta<2\theta_{0}$.) Then (26) gives a bound on the normalised energy
$$  \hatE(\theta)\leq c \left(\epsilon_{0}/(\tau \theta_{0}^{n-2}) + \theta_{0}^{2} + \tau/\theta_{0}^{n-1}\right)\  E. $$

By making $\theta_{0}$ small, then $\tau$, then $\epsilon_{0}$, we get a  $\theta\in [\theta_{0}, 2\theta_{0}
 ]$ with $\hatE(\theta)\leq \thalf E$ and by monotonicity $\hatE(\theta_{0})\leq \thalf E$.

\subsection{Some further developments}

The influence of the work  of Schoen and Uhlenbeck has been immense and extends in many directions (the paper \cite{kn:SchU1} has 303 citations on MathSciNet at the time of writing).  Singularities of the kind which which came to the fore in their paper arise in various models in Mathematical Physics, for point singularities in $\bR^{3}$. They  also appear in complex algebraic geometry as meromorphic maps (related to the work of Uhlenbeck and Yau that we discuss in subsection 6.1 below). In  Hardt's  survey \cite{kn:Hardt} of
 developments on singularities  of harmonic maps in the decade following the Schoen and Uhlenbeck paper he writes {\it \lq\lq the paper [of Schoen and Uhlenbeck] has many ideas and techniques that have proved to have wide influence in geometric analysis.''}

 {\it Stationary maps} form another important subclass of weak harmonic maps. Such a map is called stationary if the first variation of the energy vanishes for variations induced by $1$-parameter families of compactly supported diffeomorphisms of the domain. This includes the minimising maps considered by Schoen and Uhlenbeck but forms a larger class. Many of Schoen and Uhlenbeck's results were later extended to stationary maps. An important paper
 \cite{kn:Brezisetal} on maps with point singularities in $\bR^{3}$, includes  examples of  tangent maps which are not stationary or mininimising. If  $f:\bR^{3}\setminus \{0\}$ is the radial extension of a degree $1$ holomorphic (hence harmonic) map $\phi$ from $S^{2}$ to $S^{2}$ then $f$ is stationary if and only if $\phi$ is a rotation; otherwise the  first variation of the energy under motion of the singularity (with fixed boundary values) does not vanish.

The paper \cite{kn:NV} of Naber and Valtorta is one  notable more recent development in the study of singularities of harmonic maps.  A tangent map  from $\bR^{n}$ to $ N$ is called $k$-symmetric if it factors through an orthogonal  projection
$\bR^{n}\rightarrow \bR^{n-k}$. Given a minimising weakly harmonic map $f$
let $  \Sigma_{k}\subset M$  be the set of points  in $x$ such that no tangent map is $(k+1)$-symmetric. Thus $\Sigma_{0}\subset \Sigma_{1}\dots \subset \Sigma$ where $\Sigma$ is the singular set of $f$: the $\Sigma_{k}$ give a stratification of the singular set. Schoen and Uhlenbeck's work, in the proof of Theorem 5 above, shows that $\Sigma_{k}$ has Hausdorff dimension at most $k$. Naber and Valtorta prove the stronger statement that $\Sigma_{k}$ is $k$-rectifiable (in fact they show this for  stationary maps $f$). They also prove a  weak $L^{3}$ result for minimising maps
$f$:
$$  {\rm Vol}\  \{ x\in M: \vert df\vert\geq \epsilon^{-1}\}\leq C \epsilon^{3}. $$

\

\

\

The ideas and techniques in  the analysis of harmonic maps which we have discussed in  Section 3 and this Section 4 have been important in other branches of differential geometry and PDE theory. For the latter we just mention the large body of work, for example \cite{kn:BN},  on bubbling phenomena in critical exponent problems.  As we have mentioned, there are close analogues between harmonic maps and  minimal submanifold theory. The next two sections of this article will describe analogous in 
gauge theory. Another area is Riemannian geometry. In the case of 4-dimensional manifolds the $L^{2}$ norm of the Riemann curvature serves an energy functional which has analogous properties to the harmonic maps energy in dimension $2$. Partial compactness results for solutions of  the 
Einstein equations on $4$-manifolds satisfying suitable bounds on the volume, diameter and this $L^{2}$ norm were obtained by Anderson \cite{kn:Ander} and Nakajima \cite{kn:Nak2}. These results follow a similar pattern to the Sacks-Uhlenbeck theory, with \lq\lq bubbling'' at a finite set of points.
The results were extended to metrics on 4-manifolds satisfying various other equations such as  extremal K\"ahler metrics, assuming a bound on the Sobolev constant,  by Tian and Viaclovsky\cite{kn:TV}. 

The $L^{2}$-norm of the Riemann curvature is  much less effective in higher dimensions. But for limits of manifolds satisfying   Ricci curvature bounds a theory analogous to that of Schoen and Uhlenbeck was developed  by Cheeger and Colding \cite{kn:CC}. Here the volume ratio of metric balls plays a role analogous to the normalised energy.

\section{Gauge Theory}
\subsection{Background}
In the mid-1970's \lq\lq gauge theory' or \lq\lq Yang-Mills theory'' entered mathematics as a new subject, propelled by interactions with physics, and this subject became the scene for many of Uhlenbeck's most prominent achievements.

\

 We begin by reviewing the basic differential geometry. To simplify notation slightly, we will nearly always consider connections on vector bundles, say complex vector bundles, usually with Hermitian metrics on the fibres. For such a bundle $E\rightarrow M$ a connection $A$ can be identified with a covariant derivative, a differential operator
$$  \nabla_{A}: \Omega^{0}(E)\rightarrow \Omega^{1}(E), $$
and we will often not distinguish between $A$ and $\nabla_{A}$. In a local trivialisation of $E$ and local coordinates $x_{i}$ on $M$ the connection is represented by a matrix-valued $1$-form $\uA=\sum \uA_{i} dx_{i}$  and the covariant derivative, thought of as acting on vector-valued functions via the trivialisation, has components
$$   \nabla^{\uA}_{i}= \frac{\partial}{\partial x_{i}} + \uA_{i}. $$
If the trivialisation is unitary then the $\uA_{i}$ take values in the  skew-adjoint matrices ({\it i.e.} in the Lie algebra of the unitary group). The curvature of the connection is a bundle valued $2$-form $F_{A}\in \Omega^{2}({\rm End} E)$ and if the connection is unitary it lies in $\Omega^{2}({\rm ad}_{E})$, where ${\rm ad}_{E}$ is the bundle of skew-adjoint endomorphisms. In a local trivialisation, as above, the curvature is the operator given by the commutator   
$$[\nabla^{\uA}_{i},\nabla^{\uA}_{j}]= \frac{\partial \uA_{j}}{\partial x_{i}}-  \frac{\partial \uA_{i}}{\partial x_{j}}+ [\uA_{i}, \uA_{j}]. $$
 Written as a matrix-valued $2$-form
\begin{equation}  F = d\uA+ \uA\wedge \uA. \end{equation}
 A change in local trivialisation is given by a map $g$ to the structure group $U(r)$. This acts on the covariant derivative by conjugation and changes $\nabla^{\uA} $ to 
\begin{equation} g \nabla^{\uA} g^{-1}= \frac{\partial}{\partial x_{i}} + g\uA_{i}g^{-1}  - (dg)g^{-1}. \end{equation}

For any bundle-valued $1$-form $a\in \Omega^{1}({\rm End}(E))$ the operator
$  \nabla_{A}+a$ is again a covariant derivative.  We regard the space of connections $A$ as an affine space and we just write the new connection as $A+a$.
A slight variant of the preceding discussion  is to consider an automorphism $g$ of the bundle $E\rightarrow M$ covering the identity on $M$. This acts on covariant derivatives by conjugation and we can write
\begin{equation}  g(A)= A - (d_{A}g)g^{-1}. \end{equation}
In the same vein the global version of the formula (30) is
\begin{equation}  F(A+a)= F(A) + d_{A} a + a\wedge a. \end{equation}

In (32) and (33) $d_{A}$ denotes the coupled exterior derivative defined by the connection (so on bundle valued $0$-forms we could write this also as $\nabla_{A}$).

The Yang-Mills equations (for the structure group $U(r)$) arise from the functional on the space of unitary connections over a Riemannian or pseudo-Riemannian manifold $M$
$$  {\cE}(A)= \int_{M}\vert F(A)\vert^{2}. $$
Here $\vert F\vert^{2}$ is computed using the standard norm on skew-adjoint matrices and the quadratic form on $2$-forms induced by the Riemannian or pseudo-Riemannian structure. The Yang-Mills equations are the Euler-Lagrange equations associated to this functional which have the form
\begin{equation}  d^{*}_{A} F_{A}=0 \end{equation} 
where $d^{*}_{A}$ is the formal adjoint of $d_{A}$. This follows from (33) since
$$  \langle F(A+a), F(A+a)\rangle = \langle F(A), F(A)\rangle + 2\langle d_{A}a, F(A)\rangle + O(a^{2}), $$
which is $$ \langle F(A), F(A)\rangle + 2\langle
a, d^{*}_{A} F(A)\rangle + O(a^{2}). $$   
 In the case of a rank $1$ bundle, with structure group the circle $U(1)$, the Yang-Mills equations are linear. When $M$ is space-time, with the
 Lorentzian metric, we get Maxwell's equations for the electromagnetic field, but from henceforth in this article we will consider only Riemannian base manifolds. 

In the Riemannian case the functional ${\cE}$ is a positive \lq\lq energy'' functional. Just as the harmonic maps functional can be thought of loosely  as measuring the deviation of a map from a constant so, at least over a simply connected manifold $M$, the Yang-Mills functional can be thought of
measuring the deviation from a product connection. The analogy with harmonic maps has been an important guiding theme in the development of Yang-Mills theory, and one emphasised in Uhlenbeck's work, as we will see below. The critical dimension for the base manifold $M$
 in the Yang-Mills case is $4$, analogous to the critical domain dimension $2$ for harmonic map theory. The Yang-Mills functional is conformally invariant in dimension $4$. Said in another way, working over a  $\rho$-ball $B_{\rho}$ in $M=\bR^{n}$, if we rescale the ball to unit size the energy changes by a factor $\rho^{4-n}$.

Another connection between Yang-Mills theory and harmonic maps comes in the question of  \lq\lq  gauge fixing''. The local representation $\uA$ of a connection over an open set $\Omega\subset M$ depends on a choice of bundle trivialisation. By changing the trivialisation we can make $\uA$ as \lq\lq bad'' as we like; conversely we would like to choose a good representation for a given connection. A natural way to do this is to seek a \lq\lq Coulomb gauge'' in which $d^{*}\uA=0$.
This is traditional in electromagnetic theory. Staying in positive
 signature, we view a magnetic field ${\bf B}$ on $\bR^{3}$ as the curvature of  a connection on a Hermitian complex line bundle. Transferring to vector-field notation, $\uA$ becomes the magnetic potential ${\bf A}$ with ${\rm curl}\ {\bf A}= {\bf B}$ and the Coulomb gauge condition is ${\rm div}\ {\bf A}=0$. Then $\Delta {\bf A}$ is the current
${\bf J}= {\rm curl}\  {\bf B}$ and ${\bf A}= G* {\bf J}$ for the Newton potential $G$. Differentiating, this gives the Biot-Savart formula for the magnetic field generated by a current.

Going back to the general situation,
if we start with some arbitrary representation $\uA_{0}$ for the connection over $\Omega \subset M$ a representation in Coulomb gauge corresponds to
 a solution of the equation
\begin{equation}  d^{*}( g \uA_{0} g^{-1} - dg g^{-1}) =0, \end{equation} for a map $g: \Omega\rightarrow U(r)$. 
This is the Euler-Lagrange equation associated to the functional
$$  \Vert g  \uA_{0} g^{-1} - dg g^{-1}\Vert^{2}_{L^{2}}. $$
When $\uA_{0}=0$ this is the harmonic maps energy for the map $g:M_{0}\rightarrow U(r)$ and the equation (35) is the harmonic map equation.
 For general $\uA_{0}$ we have a deformation of that equation.  

\subsection{The 1982 papers in Commun. Math. Phys.}

The title of this subsection refers to the two papers of Uhlenbeck 
\cite{kn:UCMP1}, \cite{kn:UCMP2}. Along with work of Taubes from around the same time, such as \cite{kn:Taubes1}, these papers initiated the study of the analytic and PDE aspects of Yang-Mills theory, to set alongside the developments of that period of a more differential-geometric and algebro-geometric nature. 
\

The paper \cite{kn:UCMP2} bears on the gauge-fixing problem indicated at the end of the previous subsection. While the overall aim is to obtain global results, the main work takes place locally, for connections over a ball. There is a simple way to fix a gauge (the \lq\lq exponential gauge'') for a connection over  a ball using parallel transport along rays through the origin. In other words, in polar coordinates the connection form $\uA$ is determined (up to an overall conjugation) by the condition that it contains no $dr$ component. This is convenient for many purposes but is not well-suited to elliptic analysis. The curvature
 depends on one derivative of the connection so we would hope, roughly speaking, that we can choose a gauge in which $\uA$ gains one derivative compared with the curvature. But, for example, an $L^{\infty}$ bound on the curvature gives only an $L^{\infty}$ bound on the connection form in an exponential gauge, it does not control the derivatives. Similarly, the Yang-Mills equations for
$\uA$ are not elliptic in exponential gauge. On the other hand,
 if the Coulomb gauge condition $d^{*}\uA=0$ is satisfied then, in Sobolev spaces, 
 the leading term $d\uA$ in the curvature does control roughly speaking one more derivative of $\uA$ and the Yang-Mills equations are elliptic. There is a parallel discussion in Riemannian geometry,  with the traditional geodesic coordinates compared with harmonic coordinates. In the latter the Einstein equations for the metric tensor are elliptic.

The central result of \cite{kn:UCMP2} is a \lq\lq small energy'' theorem.
\begin{thm}
For $p>1$  there are $\epsilon,C>0$  such that if $A$ is a connection over $B^{n}$ with $\Vert F\Vert_{L^{n/2}}\leq \epsilon$ then there is local trivialisation in which the connection form $\uA$ has the following properties:
\begin{enumerate}\item  $d^{*}\uA=0$,
   \item On the boundary,  the contraction of $\uA$ by the normal vector vanishes,
\item   $$  \Vert \uA \Vert_{{p},{1}} \leq C \Vert F\Vert_{L^{p}}$$
\end{enumerate}
\end{thm}

Thus the curvature does control one more derivative of the connection form in this $L^{p}$ sense.  Uhlenbeck proves a stronger statement, for connections which are only in $L^{p}_{1}$ (for $p\geq n/2$) but in 
 our discussion we work with smooth connections, which makes things
 a bit simpler. Also, to simplify notation we write the proof for the case $n=4$.

Uhlenbeck uses a continuity argument to establish this result. The estimates are another instance of what we are calling
\lq\lq critical quadratic rearrangement'', as in the proof of  Theorem 3 for harmonic maps. For $\rho\in [0,1]$ let $A_{\rho}$ be the connection over the unit ball obtained by pulling back the restriction of $A$ to the $\rho$ ball by the dilation map. Then the scaling behaviour of the $L^{2}$ norm on $2$-forms in dimension $4$ shows that $\Vert F(A_{\rho})\Vert_{L^{2}}\leq \Vert F(A)\Vert_{L^{2}}$, so we have a path of connections joining $A$ to the trivial connection, all with 
$\Vert F\Vert_{L^{2}}\leq \epsilon $. The strategy is to construct a corresponding path $\uA_{\rho}$ satisfying the conditions in the statement. To set this up we consider a variant of the third condition in Theorem 7 (for $p=2$),  depending on a small number $\eta$ to be chosen below. The variant is
$$   \Vert \uA\Vert_{2,{1}}< \eta,  \ \ \ \ \ \ \ \ \ \ \ \ \ \ \ \ \ (3') $$
which is an open condition. 

 Suppose that for some $\rho$ we have  a  $\uA_{\rho}$ satisfying (1),(2) of Theorem 7 and ($3'$), and just write 
$\uA_{\rho}=\uA$. The boundary condition is elliptic for the operator $d^{*}\oplus d$ and  if $a$ is a $1$-form satisfying the boundary condition in (2) of Theorem 7 and $d^{*}a=da=0$ then $a=0$. (For we can write $a=df$ and then $f$ satisfies the Laplace equation with Neumann boundary conditions and hence is constant, so $a=0$.) Then  elliptic theory gives estimates
   \begin{equation}  \Vert \uA\Vert_{p,1}\leq K_{p} \Vert d\uA\Vert_{L^{p}}. \end{equation}
   For $p<4$ the  Sobolev embedding $L^{p}_{1}\rightarrow L^{q}$ with $q=4p/(4-p)$ combined with (36) tells us that 
   $$  \Vert \uA \Vert_{L^{q}}\leq K'_{p} \Vert d\uA\Vert_{L^{p}}. $$
 Now the formula $d\uA= F- \uA\wedge\uA$ leads to 
 \begin{equation}  \Vert \uA\Vert_{p,1}\leq  c_{1} \Vert F\Vert_{L^{p}}  + c_{2} \Vert \uA\Vert_{p,1}\Vert \uA\Vert_{2,1}, \end{equation}
where we have used H\"older's inequality,  exploiting the fact that $1/p= 1/q+ 2$. 

Now for the crucial step take $p=2$. Then  (37) gives
 
$$   \Vert \uA\Vert_{L^{2}_{1}}\leq  c_{3} \Vert F\Vert_{L^{2}} + c_{4} \Vert \uA\Vert^{2}_{L^{2}_{1}}.
$$

Thus if $\eta$ is chosen so that $c_{4} \eta<\thalf$, say, we get $\Vert \uA\Vert_{L^{2}_{1}}\leq C\Vert F\Vert_{L^{2}}$ with $C= 2 c_{3}$. In other words (1),(2) and ($3'$) imply (3).
Now choose $\epsilon< \eta/2C$ so that condition (3) implies ($3'$).
Going back to (37) and applying the same rearrangment argument 
for the quadratic term we get $L^{p}_{1}$ bounds on $\uA$ for all $p<4$ (provided $\epsilon$ is chosen suitably small).  A similar argument works for $p\geq 4$. Note that we are {\it not} supposing that the curvature is small in $L^{p}$ for $p>2$, so $\uA$ could be large in $L^{p}_{1}$, but we have some bound. Similarly, if we want to stay in the smooth category, we can estimate higher derivatives. Then it is straightforward to show that the set of $\rho\in [0,1]$ for which a solution exists is closed. The point is that the open condition ($3'$) cannot be violated in taking a limit because it is implied by the closed condition (3).
(In this case of dimension $n=4$ it is not essential to invoke $L^{p}$ theory for $p\neq 2$ to prove the main result; the proof can be done in Sobolev spaces $L^{2}_{k}$, as in \cite{kn:DK}.)

The openness part of the continuity proof uses, as usual, the implicit function theorem. Condition (3') is open by its nature so we just have to deform the solution to the equations (1),(2). This requires some technical work to set up due to the boundary condition. At a solution $\uA=\uA_{\rho}$ the linearised equation  for a matrix valued function $\psi$ is
$$ d^{*}( d\psi + [\uA,\psi]) = \sigma, $$
with Neumann boundary condition and where the integral of the given $\sigma$ is zero.
The operator on the left hand side can be written as $\Delta \psi + \{\uA,d\psi\}$
where $\{\ ,\ \}$ combines the inner product on $1$-forms with the matrix bracket. This is treated as a deformation of the ordinary Poisson equation  with Neumann boundary conditions.  Sobolev estimates similar to those used above show that the linearised equation is soluble and the implicit function theorem can be applied. 

\

The  global result in Uhlenbeck's paper \cite{kn:UCMP2} concerns the \lq\lq subcritical'' case,  with an $L^{p}$ bound on the curvature for $p>n/2$.

\begin{thm}
Let $A_{i}$ be  a sequence of unitary connections on a bundle $E$ over a compact Riemannian $n$-manifold $M$ satisfying a bound $\Vert F(A_{i})\Vert_{L^{p}}\leq C$ for some $p>n/2$. There is a subsequence $\{i'\}$ and bundle automorphisms
 $g_{i'}$ such that the  transformed connections $g_{i'}(A_{i'})$ converge weakly in $L^{p}_{1}$ to an $L^{p}_{1}$ limit $A_{\infty}$.
\end{thm}

 This is a relatively elementary consequence of Theorem 7.  First, the scaling behaviour of the $L^{p}$ norm means that there is an $r_{0}$
 such that  if the restriction of $A_{i}$
to any $r_{0}$-ball  in $M$ is pulled back to the unit ball $B\subset \bR^{n}$ via geodesic
coordinates then the pulled back connection satisfies the small-curvature
hypothesis of Theorem 7. Cover the manifold $M$ by a finite collection of such small balls. Then after passing to a subsequence and applying a sequence of gauge transformations  {\it over each ball} we can suppose that the connections converge weakly in $L^{p}_{1}$ over the ball. The problem is to convert this local convergence to the global result.

\

It is convenient  to take a different point of view here and consider
 a bundle-with-connection  over a manifold $M$ presented by the data
 \begin{itemize}
 \item an open cover $M=\bigcup U_{\alpha}$;
 \item on each overlap $U_{\alpha}\cap U_{\beta}$ a transition function
 $g^{\alpha \beta}$, 
  taking values in the unitary group, such that  $g^{\alpha \gamma}= g^{\alpha \beta} g^{\beta \gamma}$ on $U_{\alpha}\cap U_{\beta}\cap U_{\gamma}$; 
 \item on each $U_{\alpha}$ a connection form $\uA^{\alpha}$, such that
 \begin{equation} 
   dg_{\alpha \beta}= g^{\alpha \beta} \uA^{\alpha} - \uA^{\beta} g^{\alpha \beta} \end{equation} on $U_{\alpha}\cap U_{\beta}$.
\end{itemize}
(Note that (38) is equivalent to (32).)

Now suppose that we have a sequence of such data (for a fixed cover)
$g^{\alpha \beta}_{i}, \uA^{\alpha}_{i}$ and that the $\uA^{\alpha}_{i}$ converge over $U_{\alpha}$ to $\uA^{\alpha}_{\infty}$.
For the moment let us suppose that this is $C^{\infty}$ convergence. Since the unitary group is compact the formula
(38) implies that the derivatives $dg^{\alpha \beta}_{i}$ are bounded, so there is a subsequence $\{i'\}$ such that the $g^{\alpha \beta}_{i'}$ converge in $C^{0}$ 
on compact subsets of $U_{\alpha}\cap U_{\beta}$. Differentiating (38) shows that this convergence is in $C^{\infty}$ on compact subsets.
 We can slightly shrink the $U_{\alpha}$ so, without loss of generality, we can suppose that the $g^{\alpha \beta}_{i}$ converge in $C^{\infty}$ on $U_{\alpha}\cap U_{\beta}$ to a limit $g^{\alpha \beta}_{\infty}$.  This system of data $( g^{\alpha \beta}_{\infty}, \uA^{\alpha}_{\infty})$ satisfies all the conditions to define a bundle-with-connection. Let $E_{i}$ be the bundle defined by the transition functions $g^{\alpha \beta}_{i}$, for $i$ finite or infinite, and $A_{i}$ the connection on $E_{i}$.
What we want to show is that for large enough $i$ there is a bundle isomorphism
$h_{i}:E_{\infty}\rightarrow E_{i}$ such that the pull-backs $h^{*}_{i}(A_{i})$ converge to $A_{\infty}$. This boils down to the problem of choosing $U(r)$-valued functions
$h_{i}^{\alpha}$  on possibly slightly smaller sets  $U'_{\alpha}$ (which still cover)  such that
\begin{equation}    h_{i}^{\alpha} g^{\alpha \beta}_{i}= g^{\alpha \beta}_{\infty} h_{i}^{\beta}\end{equation}
on the intersections,  and with $h_{\alpha_{i}}\rightarrow 1$ as $i\rightarrow \infty$. 

By an induction argument it suffices to treat the case when the cover is by two open sets $M=U_{\alpha}\cup U_{\beta}$. We choose $h_{\beta}=1$ so the condition to solve is
$$    h_{i} = g_{\infty} g_{i}^{-1}, $$
where we write $h=h_{\alpha}, g_{i}= g_{i}^{\alpha \beta}, g_{\infty}= g_{\infty}^{\alpha \beta}$. Now $g_{i}\rightarrow g_{\infty}$ in $C^{0}$ so for $i$ large $g_{\infty}g_{i}^{-1}$ takes values in a small neighbourhood of the identity in $U(r)$ and we can write
$$   g_{\infty}g_{i}^{-1}= \exp( L_{i}), $$
for a matrix-values function $L_{i}$ on $U_{\alpha}\cap U_{\beta}$. We take a suitable  shrunken cover $U'_{\alpha}\cup U'_{\beta}$ and a cut-off function $\chi$ so that $\chi L_{i}$ is equal to $L_{i}$
 on $U'_{\alpha}\cap U'_{\beta}$ and $\chi L_{i}$ can be extended smoothly by $0$ over $U'_{\alpha}$. Then we can take $h_{i}=\exp(\chi L_{i})$ and clearly $h_{i}\rightarrow 1$ as $i\rightarrow \infty$.

In the setting of Theorem 8 we do not have $C^{\infty}$ convergence of the connection forms $\uA^{\alpha}_{i}$ but only weak $L^{p}_{1}$ convergence.
 However this implies weak $L^{p}_{2}$ convergence of the $g^{\alpha \beta}_{i}$ (after taking a subsequence) which implies $C^{0}$ convergence, since evaluation at a point is bounded on $L^{p}_{2}$ when $p>n/2$. Then the whole argument goes through unchanged. This is the crucial point 
where the condition $p>n/2$ is required.

\

\

We now turn to the other {\it Commun. Math. Phys.}  paper \cite{kn:UCMP1}.
The main result is the removability of point singularities:
\begin{thm} A Yang-Mills connection  over the punctured ball
$B^{4}\setminus \{0\}$ with curvature in $L^{2}$  extends to a smooth Yang-Mills connection over $B^{4}$.
\end{thm}

More precisely, if $A$ is a finite-energy  Yang-Mills connection on a bundle
$E\rightarrow B^{4}\setminus\{0\}$ then there is a bundle $\tilde{E}\rightarrow B^{4}$ and an isomorphism ${\iota}:\tilde{E}\vert_{B^{4}\setminus\{0\}}\rightarrow E$ such that ${\iota}^{*}(A)$ extends smoothly over the origin.

 This is the analogue of Sack's and Uhlenbeck's Theorem 4 for harmonic maps of the punctured disc and Uhlenbeck's strategy of proof is similar. For $r\leq 1$ let $$
  \cE(r)=\int_{\vert x\vert \leq r} \vert F\vert^{2}. $$
  The strategy is to derive  differential inequalities relating $\cE$ and $\frac{d\cE}{dr}$.

This paper \cite{kn:UCMP1} introduced a number of important techniques and results, on the way to the proof of Theorem 9. One was the use of exponential gauges discussed above, Another was a  small energy result:
\begin{thm} There are $\epsilon,C>0$ such if $A$ is a Yang-Mills connection over the unit ball $B\subset \bR^{4}$ with energy $\cE$ less than $\epsilon$ then $\vert F\vert^{2}\leq  C \ \cE$ on  $B_{\thalf}$.
\end{thm}
Again, this is the analogue of what we discussed for harmonic maps. It can be proved using Theorem 7 above to find a Coulomb gauge and then applying elliptic estimates, very much like the argument for Theorem 3 . But the proof in \cite{kn:UCMP1} is different. It does not require gauge fixing and introduces another important technique.   Recall that the Yang-Mills equations are $d_{A}^{*}F=0$. The curvature of any connection satisfies the Bianchi identity $d_{A}F=0$, so $\Delta_{A}F=0$ where $\Delta_{A}$ is the coupled \lq\lq Hodge'' Laplace operator
$$  \Delta_{A}= - (d_{A}d^{*}_{A}+ d^{*}_{A} d_{A}). $$
There is another  Laplace-type operator $  -\nabla^{*}_{A}\nabla_{A}$ acting on the bundle-valued forms. In a local trivialisation this is
$$   \sum \left(\nabla^{\uA}_{i}\right)^{2}. $$
The two are related by a Weitzenbock formula which, over a flat base manifold, is
   $$  \Delta_{A} \phi=\nabla_{A}^{*}\nabla_{A}\phi + \{ F, \phi\},$$
   
   where the pointwise bilinear operation $\{\ ,\ \}$ combines the  bracket on bundle endomorphisms with the map $\Lambda^{2}\otimes \Lambda^{k}\rightarrow \Lambda^{k}$ furnished by the derivative of the action of the orthogonal group on $k$-forms.

The upshot of this differential geometry is that the curvature of a Yang-Mills connection (over a flat base manifold) satisfies the equation
\begin{equation}  \nabla^{*}_{A}\nabla_{A} F= \{F,F\}. \end{equation}
One computes that $\vert \{F,F\}\vert \leq 4 \vert F\vert^{2}$ and then (40) leads to a differential inequality for $\vert F\vert$:
\begin{equation}   \Delta\vert F\vert \geq 4 \vert F\vert^{2}. \end{equation}
(The function $\vert F\vert$ may not be smooth at zeros of $F$ but this difficulty can be got around in standard ways. The exact constant $4$ in (41) will not be important and depends on conventions for defining the norm $\vert F\vert$.)

One can then apply the Nash-Moser iteration technique to derive interior estimates on $f=\vert F\vert$ from (41). At the first step, let $\chi$ be a fixed cut-off function, equal to $1$ on the ball of radius $\tthreequart$ say. Then
multiplying by $\chi^{2} f$ and integrating by parts:
$$  \int_{B^{4}}  (\nabla(\chi^{2} f), \nabla f) \leq 4 \int \chi^{2} f^{3}.$$
We have
$$   (\nabla(\chi^{2} f),\nabla f)= \vert \nabla(\chi f)\vert^{2}- \vert \nabla \chi\vert^{2} f^{2}, $$
so we get 
$$     \int\vert \nabla(\chi f)\vert^{2}\leq 4 \int (\chi f)^{2} f + \int \vert \nabla \chi\vert^{2} f^{2}. $$
Invoking the Sobolev embedding $L^{2}_{1}\rightarrow L^{4}$ in dimension four and H\"older's inequality, we see that if the $L^{2}$ norm of $f$ is sufficiently small we can apply quadratic rearrangement to get an $L^{4}$ bound on $\chi f$, hence an $L^{4}$ bound on $f$ in the $\tthreequart$-ball. The relevant manipulation is again, similar to that in the proof of Theorem 3.  Repeating the process, with a suitable sequence of concentric balls and cut-off functions and keeping track of the constants, leads to an $L^{\infty}$ bound on $f$ in the $\thalf$-ball. If one only needs to work with Yang-Mills solutions then it is  usually possible to avoid the use of the sharp Coulomb gauge-fixing result of \cite{kn:UCMP2}, using this alternative approach from \cite{kn:UCMP1}.

\

Returning to the removal of singularities problem, Uhlenbeck explains in \cite{kn:UCMP1}, for a connection over $B^{n}\setminus \{0\}$, the critical nature of the curvature decay condition $\vert F\vert \leq C \vert x\vert^{-2}$.
If we take a non-trivial  Yang-Mills connection over $S^{n-1}$ (for example the Levi-Civita connection on the tangent bundle) and pull it back  by radial projection we get a Yang-Mills connection over the punctured ball whose curvature is exactly $O(\vert x\vert^{-2})$. In dimension $n>4$ this curvature is in $L^{2}$ but when $n=4$ it is in $L^{p}$ for any $p<2$ but not in $L^{2}$.  For $n=4$ the small energy result, applied to a ball of radius $\vert x\vert/2$ centred at a point $x$ leads to a bound 
\begin{equation} \vert F(x)\vert = o(\vert x\vert^{-2}) \end{equation} so we get a little above the critical $O(\vert x\vert^{-2})$ threshold but, as in the case of  Theorem 4 for harmonic maps, more is needed. A differential inequality
$$   (1-\delta) \  \cE \leq \frac{1}{4r} \frac{d\cE}{dr} $$
implies that $\cE =O(r^{4-4\delta})$, which gives $\vert F\vert=  O(r^{-2\delta})$
and then $F$ is in $L^{p}$ for $p<2/\delta$. Uhlenbeck  establishes a more complicated inequality
\begin{equation}    \left( 1- \omega\  \cE(2r)^{\thalf}\right)\cE(r)\leq \frac{1}{4r} \frac{d\cE}{dr}, \end{equation} for a constant $\omega$, from which she is able to deduce that $\vert F\vert$ is bounded.
 Then an exponential gauge produces a bounded connection form, which can be adjusted to satisfy the Coulomb condition, and elliptic regularity shows that the connection extends smoothly over the origin. (See the further discussion in 6.3 below.)

\

 We will not go much further into the details of Uhlenbeck's proof in 
 \cite{kn:UCMP1},
partly because we will discuss another proof, of Uhlenbeck and Smith, in 6.3 below. One important idea in the proof of \cite{kn:UCMP1} is  that on any small annulus the connection is close to flat, in the sense that when the annulus is scaled to standard size the curvature is small, by (42) above. This means that the nonlinear equation can be approximated by its linearisation, provided a suitable gauge is used over the annulus. The construction of these gauges takes up much of the work in the paper. In the case of  the Abelian gauge group $U(1)$ we have $F= d\uA$ and the Yang-Mills equation is $d^{*}\uA=0$. Then over a domain $\Omega$:
$$  \int_{\Omega}\vert F\vert^{2}= \int_{\partial \Omega} \uA\wedge *F. $$
In the general non-Abelian case there is a similar formula with lower order terms, and by estimating these Uhlenbeck  obtains the  differential inequality (43). The crucial constant $4$ in (43) appears as the first eigenvalue of the Laplacian on co-closed $1$-forms on $S^{3}$. 

\subsection{Applications}
 
 Much of the original motivation for the removal of singularities theorem
 was to answer a question raised by physicists. A finite-energy Yang-Mills $U(r)$ connection over $\bR^{4}$ extends smoothly to $S^{4}$, in particular it has a  topological invariant Chern number. In a later paper   \cite{kn:UCMP3} Uhlenbeck showed that for any finite energy connection the Chern number, defined by integrating the Chern-Weil form, is an integer. This used the full force of Theorem 7, in fact extended to $L^{2}_{1}$ connections.
 
The results of these two papers of Uhlenbeck were foundational in the development of Yang-Mills theory over  Riemannian manifolds of dimension at most 4. In one direction they opened the way to the use of variational methods. In the subcritical case, over  manifolds of dimension  $2$ and $3$, the analysis is relatively straightforward. In particular any bundle admits a connection which minimises the Yang-Mills functional. This is an easy consequence of Theorem 8 which implies that there is a minimising sequence which converges in $L^{2}_{1}$ and in these dimensions the bundle theory works in a straightforward way with $L^{2}_{1}$ connections and $L^{2}_{2}$ gauge transformations because  $L^{2}_{2}$ maps are continuous.  The existence of minimisers is not the end of the story: one would like to go on to the relate the Yang-Mills connections to the topology of the space of connections modulo gauge equivalence. There has been a lot of work on this in the case of bundles over surfaces, following Atiyah and Bott, including a paper of Daskalopoulus  and Uhlenbeck \cite{kn:DU1}.  Rade (who, like Daskalopoulos, was a PhD student of Uhlenbeck) obtained complete results about the Yang-Mills gradient flow  in these dimensions \cite{kn:Rade}.

\

The variational theory in dimension $4$ is much harder. One early and clear cut result was obtained by Sedlacek (another PhD student of Uhlenbeck)
\cite{kn:Sed}.
Let $G$ be a compact Lie group. A principal $G$-bundle $P$ over a compact oriented 4-manifold is determined by two characteristic classes $\kappa(P), w(P)$ where
$\kappa(P)\in H^{4}(X, \pi_{3}(G))$ and $w(P)\in H^{2}(X,\pi_{1}(G))$. For example if $G=SU(m)$ then $\kappa$ is the second Chern class and $w$ is trivial, since $SU(m)$ is simply connected, while if $G=SO(m)$ for $m=3$ or $m\geq 5$ then $\kappa$ is a multiple of the first Pontrayagin class and $w$ is the second Stiefel-Whitney class in $H^{2}(X;\bZ/2)$. Now let $X$ have a Riemannian metric, so the Yang-Mills equations are defined. Sedlacek's result is that for any $P$ there is a Yang-Mills connection on a $G$-bundle $P'\rightarrow X$ with $w(P)=w(P')$. It might happen that $P'$ is isomorphic to $P$ but the result allows the possibility that they are different. (This is the Yang-Mills 
analogue of the fact that for any homotopy class of maps from a surface 
 there is a harmonic map inducing the same homomorphism on fundamental groups, but possibly in a different homotopy class.)

To prove this result, Sedlacek considers a minimising sequence $A_{i}$ for the Yang-Mills functional on $P$. A covering argument just like that for Theorem 2 shows that after passing to a subsequence $i'$, there is a finite subset $S\subset X$ (possibly empty) such that each point in $X\setminus S$ is the centre of a ball on which the connections $A_{i'}$ satisfy the small energy condition for the Coulomb gauge fixing result, Theorem 7. Then we arrive at a situation like that we considered in the proof of Theorem 8, with  a subsequence $i'$, a cover $B_{\alpha}$ of $X\setminus S$ and connection $1$-forms
$\uA^{\alpha}_{i'}$ in Coulomb gauge and with $L^{2}_{1}$ limits
$\uA^{\alpha}_{\infty}$ and $L^{2}_{2}$ limits $g^{\alpha \beta}_{\infty}$ of the transition functions $g^{\alpha \beta}_{i}$. Sedlacek shows that these limits $\uA^{\alpha}_{\infty}$ are  weak solutions of the Yang-Mills equations and then, using ellipticity
in Coulomb gauge, that they are in fact smooth.

The equation $(38)$, and the cocycle conditions, are preserved in the limit so the smoothness of the $\uA^{\alpha}_{\infty}$ implies that  of the $g^{\alpha \beta}_{\infty}$. Thus the limiting data defines a smooth connection
 $A_{\infty}$ on a bundle $P_{\infty}$ over $X\setminus S$. This is a finite-energy
 Yang-Mills connection  so the removal of singularities theorem shows that the bundle and connection extend smoothly over the finite set $S$ to a connection on a bundle $P'$.  The remaining step is to show that the characteristic classes
$w(P), w(P')$ are equal. The point here is that $P'$ may not be isomorphic to $P$ over $X$. Even over $X\setminus S$ the $L^{2}_{2}$ convergence of the transition functions does not give the $C^{0}$ convergence, so the last part of the proof of Theorem 8 does not extend to this situation (although, by algebraic topology, it does turn out in the end that the bundles are isomorphic over $X\setminus S$). 

\

Consider for example  the case when $X=\bC\bP^{2}$, with its standard metric and orientation, and the gauge group $G=SU(2)$. If $P$ is the bundle with $c_{2}(P)=1$  there is no
minimiser of the Yang-Mills functional: a  minimising sequence will converge to the flat connection away
from a point in $\bC\bP^{2}$, displaying the same kind of bubbling behaviour
that we described for harmonic maps from $T^{2}$ to $S^{2}$ in Section 3, If $c_{2}(P)=-1$ on the other hand then  minimisers can be constructed (they are instantons, see below). But there is a  non-compact moduli space of these minimisers and we could choose a \lq\lq bad''
 minimising sequence exhibiting bubbling over a point.

The variational theory was developed much further in a sequence of papers by Taubes such as \cite{kn:Taubes2},\cite{kn:Taubes3},  relating the solutions to the topology of the space of connections modulo equivalence. We refer to the article \cite{kn:SKDNotices}  for a
 discussions of those   developments.  
   
\

Perhaps the largest impact of Uhlenbeck's papers \cite{kn:UCMP1},\cite{kn:UCMP2} came 
in the study of  the \lq\lq instanton'' solutions to the
 Yang-Mills equations in 4-dimensions. Instantons, over an oriented, Riemannian 4-manifold $X$ are connections whose curvature is self-dual or anti-self-dual (the two being interchanged by switching orientation). In the analogy between 4-dimensional Yang-Mills theory and harmonic maps of surfaces, these correspond to holomorphic maps to a K\"ahler (or just almost-K\"ahler) manifold. The instantons have consequences in 4-manifold topology, analogous to those of the mapping theory in symplectic topology. The analogue of Theorem 2 for sequences of instantons leads to the \lq\lq Uhlenbeck compactifications'' of instanton moduli spaces. These are made up of pairs $([A], D)$ where
$[A]$ is the gauge equivalence class of an instanton and $D$ is a formal sum of points $q_{j}$ of $X$ with multiplicities $\kappa_{j}$. A sequence of
instantons converges to $([A],D)$ if the connections converge on the complement of the points $q_{j}$ and exhibit \lq\lq bubbling'' over the $q_{j}$ with $\kappa_{j}$ units of energy (suitably normalised) concentrating at $q_{j}$.
We refer to the books \cite{kn:DK}, \cite{kn:FU} for  detailed accounts of these developments.

\section{The Yang-Mills equations in higher dimensions}

 \subsection{Hermitian Yang-Mills connections on stable bundles}

 Much of the work involving Yang-Mills theory over manifolds of dimension greater than four focuses on manifolds with some extra structure, as opposed to general Riemannian manifolds. One of the most important such developments came in work of Uhlenbeck and Yau \cite{kn:UY}, establishing the \lq\lq Kobayashi-Hitchin conjecture'' in complex differential geometry. 

Let $X$ be a compact K\"ahler manifold of complex dimension $m$ with K\"ahler form $\omega$. The forms on $X$ decompose into bi-type and the metric defines a contraction operator $\Lambda: \Omega^{1,1}\rightarrow \Omega^{0}$. This is just the trace with respect to the metric $\omega$. We recall three  identities:
\begin{enumerate}
\item For a function $f$:\ \ \  $\Delta f= -2i\Lambda \db \partial f$.
\item For a $(0,1)$ form $\alpha$:\ \ \  $ i \Lambda (\alpha\wedge \overline{\alpha})= -\vert \alpha\vert^{2}$.
\item For a $(1,1)$ form $\theta$:\ \ \ $ \theta\wedge \omega^{n-1}= (\Lambda \theta)\  \omega^{m}/m= (m-1)! (\Lambda \theta) \ {\rm vol}$.
\end{enumerate}

 Now consider a unitary connection on a complex vector bundle $E\rightarrow X$. The curvature decomposes into $F=F^{0,2}+F^{1,1}+ F^{2,0}$ and we define $\hatF= i \Lambda F^{1,1}$. So $\hatF$ is a section of the bundle of self-adjoint endomorphisms of $E$. The connection is
 called a {\it Hermitian-Yang-Mills connection} if $F^{0,2}$ and $ F^{2,0}$ vanish and  
\begin{equation}   \hatF =  \mu 1_{E}, \end{equation}
for a constant $\mu$. The constant $\mu$ is determined by topology. By the  third item above
$$ (m-1)! \left( {\rm Tr}\  \hatF\right) {\rm vol} = i {\rm Tr} F \wedge \omega^{m-1}, $$
and by Chern-Weil theory the $2$-form $ (i/2\pi) {\rm Tr}\ F$ represents the first Chern class $c_{1}(E)$. So if a solution to exists we have
$$  \mu = \frac{2\pi}{ (m-1)! {\rm Vol}\ (X)} \ \ \frac {{\rm deg} (E)}{{\rm rank} E}, $$
where  the degree ${\rm deg}(E)$ is defined to be the pairing 
$(c_{1}(E)\cup \omega^{m-1})[X]$, which is a topological invariant of the bundle $E\rightarrow X$  and the  K\"ahler class $[\omega]$.
The ratio ${\rm deg}\ (E)/{{\rm rank}\  E}$ is called the {\it slope} of the
 bundle $E$.  Hermitian-Yang-Mills connections are Yang-Mills connections: in fact they are absolute minimisers of the  Yang-Mills functional on the given bundle. If $E$ is the tangent bundle of $X$ and the connection is the Levi-Civita connection then $\hatF$ is the Ricci tensor, so (44) is related to the Einstein equations and solutions are often called Hermitian-Einstein connections in the literature.

The significance of the condition that the curvature $F$ has type $(1,1)$ is that this implies that the connection is compatible with a holomorphic structure on the bundle $E$. For any connection $\nabla$ we can write
$$  \nabla= \partial_{\nabla} \oplus \db_{\nabla}:\Omega^{0}(E)\rightarrow \Omega^{1,0}(E)\oplus \Omega^{0,1}(E). $$
The sheaf of local solutions of the equation $\db_{\nabla}s=0$ is a sheaf of modules over the structure sheaf of the complex manifold $X$:  the local holomorphic functions. But when ${\rm dim}_{\bC} X>1$ the equation $\db_{\nabla} s=0$ is overdetermined and for a general connection the only solution will be $s=0$. The condition $F^{0,2}=0$ is the integrability condition for this equation, which implies the existence of solutions generating the bundle $E$, thus defining a holomorphic structure on the bundle. For a unitary connection the component $F^{2,0}$ is $- (F^{0,2})^{*}$ so the vanishing of one implies the same for the other.

 We have then the existence question: given a holomorphic bundle $E$ does it admit a compatible Hermitian-Yang-Mills connection? The Kobayashi-Hitchin conjecture, formulated independently by Kobayashi and Hitchin around 1980, is that  such a connection exists if and only if $E$ is a direct sum of {\it stable} holomorphic bundles of equal slope. The notion of stability here was introduced before by algebraic geometers in the context of  moduli problems. A holomorphic bundle of $V$ is defined to be stable if  every non-trivial coherent subsheaf ${\cal S}$ of rank less than ${\rm rank} V$ satisfies the condition
\begin{equation}  {\rm slope} \ ({\cal S})< {\rm slope} \ (V). \end{equation}
(By general theory, such a subsheaf is given by a proper subbundle of $V$ outside a singular set of complex codimension $2$ or more, so the first Chern class of ${\cal S}$ is defined in $H^{2}(X)$.) 

Part of the evidence for this conjecture came from results of Narasimhan and Seshadri from the 1960's which covers the case when $X$ is a complex curve. The fact that the existence of a Hermitian-Yang-Mills connection implies that the bundle is a sum of stable bundles is relatively straightforward and was proved by Kobayashi \cite{kn:Kob}. In this subsection we discuss Uhlenbeck and Yau's proof of the existence result, for general Kahler manifolds $(X,\omega)$. The essential statement can be put in the form:
\begin{thm}
If a holomorphic bundle $E$ does not admit a Hermitian-Yang-Mills connection then there is a subsheaf ${\cal S}$, as above, with ${\rm slope}\ {\cal S}\geq {\rm slope}\ (V)$.
\end{thm}

 There are two ways of setting up differential geometry on holomorphic vector bundles. In one---which is the traditional point of view in complex differential geometry---one has a fixed holomorphic bundle and varies the Hermitian metric. A metric defines a unique compatible connection
(often called the \lq\lq Chern connection''). In the other, we fix the metric $(\ ,\ )$ on a $C^{\infty}$ bundle $E$ and vary the connection. Let $\nabla_{0}$ be some unitary reference connection and write
$  \nabla_{0} = \partial_{0} + \db_{0}$.
If $g$ is any automorphism of $E$, not necessarily unitary, we define a new covariant derivative by
$$   \nabla^{g}=  (g^{*})^{-1}\circ \partial_{0}\circ g^{*}+ g \circ \db_{0} \circ g^{-1} $$
That is, we conjugate $\db_{0}$ by $g$ and define the $(1,0)$ part in the unique way to make a unitary connection $\nabla^{g}$. If $g$ is unitary then $g=(g*)^{-1}$ and we have the ordinary gauge transformation, so $\nabla^{g}$ geometrically equivalent to $\nabla_{0}$. In general, the $\db$-operators of $\nabla^{g}, \nabla_{0}$ are equivalent, in the sense that they define isomorphic holomorphic structures on the bundle, but the connections are essentially different. Working modulo the unitary gauge transformations, we can restrict attention to self-adjoint automorphisms which we write as $h$. Then 
$$  \nabla^{h}= h  \circ \db_{0} \circ h^{-1} + h^{-1}\circ \partial_{0}\circ
h. $$
 To match up with the first point of view, it is equivalent to fix the holomorphic structure with $\db$-operator $\db_{0}$ and vary the metric to $(hs,hs)= (s,h^{2}s)$.  We will use this second point of view, which means that some formulae will look different from those in \cite{kn:UY}.

The curvature of the connection $\nabla^{h}$ is
\begin{equation}F(\nabla^{h})= F_{0}+ \db_{0}(h^{-1}\  \partial_{0} h)- \partial_{0}(\db_{0} h \ h^{-1}) - \left( h^{-1} \partial_{0} h\  \db_{0}h h^{-1} + \db_{0} h \  h^{_2} \  \partial_{0} h\right) \end{equation}
where $F_{0}$ is the curvature of $\nabla_{0}$. 

Write $h=e^{u}$, so $u$ is a section of the bundle of self-adjoint endomorphisms of $E$. To gain understanding of  the nature of the Hermitian-Yang-Mills equation we can consider the linearisation about $u=0$. Using the connection between the $\partial$ and $\db$ operators and the Laplacian one sees that this is
\begin{equation}  \hatF(e^{u})= \hatF_{0} + i\Lambda(\db_{0}\partial_{0}- \partial_{0} \db_{0}) u + O(u^{2}) =\hatF_{0} + \nabla^{*}_{0}\nabla_{0} u + O(u^{2}).  \end{equation}
So the linearisation is the coupled Laplacian. The Hermitian-Yang-Mills equation has many similarities with the harmonic equation for a map into the space of hermitian matrices: the nonlinear term is quadratic in the first derivatives of $h$.

 The proof by Uhlenbeck and Yau of Theorem 11 uses a continuity method,  with the   family of equations, for
$t\geq 0$:
     \begin{equation}  \hatF(e^{u})= -t u . \end{equation}
 To set things up they prove:
 \begin{prop}
 \begin{enumerate}  
 \item The equation (48) has a solution for large $t$.  
    \item The set of $t\in [0,\infty)$ for which a solution to (48) exists is open.
    \item  If there is a smooth family of solutions $u_{t}$ to (48) for $t$ in an interval $(t_{0}, t_{1})$ satisfying a bound $\Vert u_{t}\Vert_{L^{\infty}}\leq C$ then the solution extends to the closed interval $[t_{0}, t_{1}]$.
\end{enumerate}
\end{prop}

Item (1) is proved by Uhlenbeck and Yau with an auxiliary
 continuity argument.
  It can also be established using the implicit function theorem, writing
  $u= \epsilon \hatF_{0} + w$ with $\epsilon=t^{-1}$. When $\epsilon=0$ there is a trivial solution $w=0$ and this can be deformed to a solution for small
$\epsilon$.

The proof of item (2) also uses the implicit function theorem, in a standard way once one knows the invertibility of the linearised operator. This involves slightly complicated calculations, extending the formula (47), which we pass over here.

Uhlenbeck and Yau use  an interesting technique to prove item (3), based on an interpolation inequality. For large $p$:
\begin{equation}  \Vert v \Vert^{2}_{L^{2p}_{1}}\leq c \Vert v \Vert_{L^{p}_{2}}\ \Vert v \Vert_{L^{\infty}}. \end{equation}
  Let $\dot{u},\dot{h}$ be the $t$-derivatives of $u$ and $h=e^{u}$ on the interval $(t_{0}, t_{1})$.
  Thus $\dot{u}$ satisfies the linear equation obtained by differentiating
  (48). Applying the maximum principle to this equation they show that
  $\Vert \dot{u}\Vert_{L^{\infty}}$ satisfies a fixed  bound  (this step is related to the invertibility of the linearised operator for item (2)).
The formula (46) for the curvature leads to an expression for   $\nabla^{*}_{0} \nabla_{0} \dot{h}$ in terms of $h, \hatF(h)$ and the derivatives of $h$. The hypothesis means that $\hatF(h)$ is bounded and one gets 
$$   \vert \nabla_{0}^{*} \nabla \dot{h}\vert \leq c\left( 1+ \vert \nabla_{0} h\vert \vert \nabla_{0} \dot{h}\vert + \vert \nabla_{0} h\vert^{2} \vert \dot{h}\vert \right) . $$
The $L^{\infty}$ bound on $\dot{u}$ gives one on $\dot{h}$ and we get
 $$  \Vert \dot{h}\Vert_{p,2}\leq c\left( 1 + \Vert  h\Vert_{2p,1} \Vert  \dot{h}\Vert_{2p,1}+ \Vert h \Vert^{2}_{2p,1}\right). $$
Let $M(t)= \Vert  h_{t}\Vert_{p,2}$, so $\vert \frac{dM}{dt}\vert \leq \Vert \dot{h}\Vert_{p,2}$.
The interpolation inequality (49), applied to $h$ and to $\dot{h}$, gives
$$  \vert \frac{dM}{dt} \vert \leq c\left( 1+ M+ \sqrt{ M\  \vert\frac{dM}{dt}\vert} \right), $$
which implies that
$$   \vert\frac{dM}{dt}\vert \leq c (1+ M). $$
 This gives a bound on $\Vert  h_{t}\Vert_{L^{p}_{2}} $ over the finite interval $(t_{1}, t_{2})$. Then it is straightforward to obtain bounds on all higher derivatives,  which implies that the solution extends to the end points.

Uhlenbeck and Yau show that a solution to  (48) for $t>0$ satisfies an {\it a priori}
 $L^{\infty}$ bound 
 \begin{equation} \Vert u_{t}\Vert_{L^{\infty}}\leq c t^{-1},
  \end{equation}
  (see Item 1 of Proposition 8 below). 
   Then Proposition 7 implies that solutions exist for all $t>0$.  If there is a $C$ such that $\Vert u_{t}\Vert_{L^{\infty}}\leq C$  for all small $t>0$ then Proposition 7 implies that the solution exists for $t=0$; so we have a Hermitian-Yang-Mills connection. The plan of the proof of Theorem 11 is to show that
 if there is no such $C$---so there is a sequence $t_{i}\rightarrow 0$ such that $\Vert u_{t_{i}}\Vert_{L^{\infty}}\rightarrow \infty$---then there is a subsheaf ${\cal S}$ with   ${\rm slope}\ {\cal S}\geq
{\rm slope}\ (V)$.

 Uhlenbeck and Yau establish the following {\it a priori} estimates for  a solution $u=u_{t}$ of (48), writing (for convenience below) $f=2^{-\thalf} \vert u\vert$.
\begin{prop}
\begin{enumerate}
\item $\Vert f \Vert_{L^{\infty}}\leq c t^{-1}$;
\item $\Vert f \Vert_{L^{\infty}}\leq c \Vert f \Vert_{L^{1}}$;
\item $\Vert \nabla f \Vert^{2}_{L^{2}}\leq c \Vert f \Vert_{L^{\infty}}$;
\item $\Vert \nabla_{0} u \Vert^{2}_{L^{2}}\leq c (1+ \Vert f \Vert^{2}_{L^{\infty}})$. 
\end{enumerate}
\end{prop}

     For simplicity, we will discuss the proofs in the case of a rank $2$ bundle $E$ with
trivial determinant. Then we can suppose that the trace of $\hatF_{0}$ is
zero and we restrict to trace-free $u$.
We recall the differential geometric theory for a bundle $E$ decomposed as an orthogonal direct sum
$E_{I}\oplus E_{II}$. A unitary connection on $E$ is defined by connections on $E_{i}$ and a second fundamental form $B$, which is a $1$-form with values in ${\rm Hom}(E_{II}, E_{I})$. In matrix notation we can write our connection as
\begin{equation} \left(\begin{array}{cc} \nabla_{I} & B\\-B^{*} & \nabla_{II}\end{array}\right) .\end{equation}
The curvature is 
\begin{equation} \left(\begin{array}{cc} F_{I}- BB^{*} & d_{I, II} B\\-d_{I, II}B^{*} & F_{II}-B^{*}B\end{array}\right)
\end{equation}
where $d_{I, II}$ is the coupled exterior derivative defined by $\nabla_{I}$ and $ \nabla_{II}$. Over a complex manifold we can write the $\db$-operator on the direct sum as
$$ \left(\begin{array}{cc} \db_{I} & \beta\\ \gamma & \db_{II}\end{array}\right)
$$
where $\beta,\gamma$ are bundle-valued $(0,1)$-forms. Then $B=\beta -\gamma^{*}$. For a connection defining a holomorphic structure on $E$ the component $\beta$ vanishes if and only if $E_{II}$ is a holomorphic subbundle and similarly for $\gamma$ and $E_{I}$. The (1,1) parts of the quadratic terms in (52) are 
$$  \left( BB^{*}\right)_{1,1}= \beta\beta^{*}+\gamma^{*}\gamma\ \ ,\ \left( B^{*}B\right)_{1,1}= \gamma\gamma^{*}+\beta^{*}\beta. $$
The crucial point for us is that the constituents have a definite sign. Using  the second of the three formulae stated 
 at the beginning of this subsection we get
$$ \vert \beta\vert^{2}= -{\rm Tr} (i\Lambda (\beta\beta^{*}))=  {\rm Tr} (i\Lambda \beta^{*}\beta)\ \ \ ,\ \ \ \vert \gamma\vert^{2} =
 {\rm Tr}\ (i\Lambda (\gamma\gamma^{*}))\vert = -{\rm Tr}\ (i\Lambda (\gamma^{*}\gamma))
 .$$
 When $E_{I}$ is a holomorphic subbundle, so $\gamma=0$, this is an aspect of the  principle that \lq\lq curvature decreases in holomorphic subbundles and increases in holomorphic quotients''.  
 
 To apply this in our situation, work initially over the open set $\Omega$ in $X$ where $u\neq 0$. Then $u$ has eigenvalues $-f,f$ 
 and the bundle $E$ is decomposed into a sum of eigenspace line bundles $E_{I},E_{II}$. Thus
$$ h = \left(\begin{array}{cc} e^{-f} & 0 \\ 0 & e^{f}\end{array}\right)$$
We find that 
$$ \db^{h}= \left(\begin{array}{cc} \db_{I} & e^{-2f} \beta\\ e^{2f}\gamma & \db_{II}\end{array}\right), $$
and 
$$  \hatF(u)= \hatF_{0}- \left( \begin{array}{cc} P & Q\\-Q^{*} &- P\end{array}\right) $$
where 
$$ P=  -\Delta f + (e^{4f}-1)\vert\gamma\vert^{2}+ (1-e^{-4f})\vert
\beta\vert^{2}. $$
So we have
\begin{equation}
 -\Delta f + (e^{4f}-1)\vert \gamma\vert^{2}+ (1-e^{-4f})\vert\beta\vert^{2} +  t f = p(\hatF_{0}),\end{equation}
 where $p(\hatF_{0})$ is  the component of $\hatF_{0}$ in ${\rm End} E_{I}$.

The maximum principle applied to this equation (53) implies the first item in Proposition 8.  
We have $\vert p(\hatF_{0})\vert\leq C$ for some $C$
 depending only on $\nabla_{0}$
and we get ${\rm max} f\leq Ct^{-1}$.  
Now we can feed this back into (53) to get the 
 the differential inequality
\begin{equation}   -\Delta f + (e^{4f}-1)\vert \gamma\vert^{2}+
 (1-e^{-4f})\vert\beta\vert^{2}\leq 2C . \end{equation}
 In particular we have $-\Delta f\leq 2C$ and while we have derived this over the open set where $u\neq 0$ it is straightforward to see that the inequality holds in a weak sense  over the whole manifold. If $f$ attains its a maximum at a point $p\in X$ then, by considering a comparison function and applying the maximum principle, we see that it is close to the maximum over a ball of fixed size about $p$, so the $L^{1}$ norm of
 $f$ is comparable to the $L^{\infty}$ norm,  which gives the second item of Proposition 8. Next, taking the $ L^{2}$ inner product of (54) with $f$ we have
$$  \int \vert \nabla f\vert^{2} \leq  c\int f\leq c \Vert f \Vert_{L^{\infty}} , $$
which is the third item.

Over $\Omega$ we have 
$$  \nabla_{0} u = \left(\begin{array}{cc}- \nabla f & f (\beta-\gamma^{*})\\ f (\beta^{*}-\gamma)
& \nabla f \end{array}\right), $$
so $$ \vert \nabla_{0} u \vert^{2}= 2( \vert \nabla f\vert^{2}+ f^{2} (\vert \beta\vert^{2}+\vert\gamma\vert^{2})). $$
This formula holds over the whole of $X$, since $f$ vanishes at the points where $\beta,\gamma$ are undefined. Similarly for the inequality (54).
Write $N=\Vert f \Vert_{L^{\infty}}$. Then it is elementary that for a suitable $\kappa$
$$  f^{2}\leq \kappa (1+ N)^{2} (e^{4f}-1)\ \ , \  \ f^{2}\leq \kappa (1+N)^{2} (1-e^{-4f}), $$
so $$ \vert \nabla_{0} u\vert^{2} \leq 2\vert \nabla f\vert^{2} + 2\kappa (1+N)^{2} ( (e^{4f}-1)\vert \gamma\vert^{2}+
 (1-e^{-4f})\vert\beta\vert^{2}). $$
 Comparing with (54) and integrating  over $X$ (so the term $\Delta f$ in (54) integrates to zero) we get
$$  \int \vert \nabla_{0}u \vert^{2}\leq 2 \int
 \vert \nabla f\vert^{2} + 2 C \kappa (1+N^{2}){\rm Vol}\ X,
 $$
 which gives item (4) of Proposition 8.

 Define $v_{t}= N(t)^{-1} u_{t} $,  where $N(t)= 
 \Vert f_{t}\Vert_{L^{\infty}}= 2^{-1/2} \Vert u \Vert_{L^{\infty}}$, as above. So the $v_{t}$ are bounded in $L^{\infty}$  and their $L^{1}$ norm has a strictly positive lower bound by item (2) of Proposition 8. By item (4) of Proposition 8 the $v_{t}$ are bounded in $L^{2}_{1}$. Suppose that there is a sequence $t_{i}$ such that $N(t_{i})\rightarrow \infty$. Passing to a subsequence, we can suppose that the $v_{t_{i}}$ have a weak $L^{2}_{1}$ limit $v_{\infty}$. This is not zero because of the lower bound on the $L^{1}$ norms. By item (3) of Proposition 8 the derivatives of $\vert v_{i}\vert$ tend to zero in $L^{2}$ and it follows that the derivative of $\vert v_{\infty}\vert$ is zero and $\vert v_{\infty}\vert$ is a non-zero constant, $\rho$ say.
  (Most likely   $\rho=\sqrt{2}$,  by our normalisation.) It follows then that the $L^{2}_{1}$ bundle endomorphism $\pi$ defined by
 $\pi= 2 \sqrt{2}\rho^{-1} v_{0}- 1$ is a rank $1$ orthogonal projection with $\pi^{*}=\pi$ and $\pi^{2}=\pi$.

 The image of $\pi$ is the limit, in a suitable sense,  of the  small-eigenvalue eigenspaces of the $u_{t_{i}}$. To visualise what is going on here, recall that the space ${\cal H}$ of positive self-adjoint $2\times 2$ matrices with determinant $1$ is a model for hyperbolic $3$-space and has a natural compactification to a closed $3$-ball, with a $2$-sphere at infinity. More intrinsically, if ${\cal H}$ consists of Hermitian forms  on a $2$-dimensional complex vector space $V$ then the sphere at infinity is $\bP(V)$. A sequence
$H_{i}\in {\cal H}$ tends to a point $[z]\in \bP(V)$ if $H_{i}(z,z)\rightarrow 0$.  In our situation we are considering sections $h_{t}$ of a bundle ${\cal H}_{E}$ over$X$ with fibre
${\cal H}$ which is  compactified by adjoining $\bP(E)$. The conclusion above is  that if the $h_{t}$ do not have a finite limit, as a section of ${\cal H}_{E}$---which would give a Hermitian-Yang-Mills connection---they have a limit \lq\lq at infinity'' which is a section of $\bP(E)$, {\it i.e.} a subbundle of $E$.

 Suppose that  $L\subset E$ is a holomorphic subbundle. Orthogonal projection onto $L$ is a smooth self-adjoint section $\varpi$ of ${\rm End}\ E$ with
$$  \varpi^{2}= \varpi \ \ \ \ \ {\rm and}\  \  \ \ \ (1-\varpi)  \db_{0}\varpi=0.$$
From another point of view, the section $\varpi$ is equivalent to a section of the projectivised bundle $\bP(E)$ and the equation $(1-\varpi)  \db_{0}\varpi=0$ is equivalent to the Cauchy-Riemann equation for this section. More generally, if ${\cal S}$ is a rank-$1$ subsheaf of $E$ we get a meromorphic section of $\bP(E)$. The basic example, in local coordinates $z_{1}, z_{2}$ on  a complex surface $X$ and a local holomorphic trivialisation of the bundle $E$, is the sheaf ${\cal S}$ which is defined by the image of the bundle map
$  s: {\cal O}\rightarrow {\cal O}\oplus {\cal O}$ with $s(z_{1}, z_{2})= (z_{1}, z_{2})$. Then the meromorphic section is given by the standard rational map $(z_{1}, z_{2})\mapsto [z_{1}, z_{2}]$ from $\bC^{2}$ to $\bP^{1}$,
 undefined at the origin.

For the $L^{2}_{1}$ projection $\pi$ constructed above, $(1-\pi)\db_{0}\pi$ is defined in $L^{2}$.  
The proof of Uhlenbeck and Yau is completed by showing three things.
\begin{enumerate}
\item $\pi$ satisfies the equation  $(1-\pi)\db_{0}\pi=0.$ 
 \item This defines a meromorphic section of $\bP(E)$ and   a coherent subsheaf of $E$.
\item The degree of this subheaf is $\geq 0$.
\end{enumerate}

To see the idea of the proof we consider the simple situation when the $u_{t}$ do not vanish anywhere and the $v_{t}$ converge in $C^{\infty}$ to $v_{0}$. In that case we have a smooth $\pi_{t}$ defined by projection onto the small eigenspace $E_{I}$ of $u_{t}$ and $\pi_{t}$ is the $C^{\infty}$ limit of the $\pi_{t}$. Then
$$  \vert (1-\pi_{t})\db_{0}\pi_{t}\vert = 
\vert \gamma_{t}\vert $$ while
$$ \int  \vert \gamma_{t}\vert^{2} (e^{4f_{t}}-1) \leq 2 C{\rm Vol} M.
$$  
Our simplifying assumptions imply that ${\rm min}_{X} f_{t}\rightarrow \infty$ as $t\rightarrow 0$, so clearly $(1-\pi_{t})\db_{0}\pi_{t}$ tends to $0$ in $L^{2}$.  In  this simple situation item (2) is trivial so we turn to item (3).   By the Chern-Weil formula the degree of $E_{I}$ is
   $$ (2\pi)^{-1} \int_{M}  p(\hatF_{0})+\vert \gamma_{t}\vert^{2} -\vert \beta_{t}\vert^{2}, $$
   while
   $$  p(\hatF_{0})+ (e^{4f_{t}}-1)\vert \gamma_{t}\vert^{2} + (1-e^{-4f_{t}})\vert\beta_{t}\vert^{2} = t f. $$
So the degree of $E_{I}$ is
$$ (2\pi)^{-1} \int_{M} t f_{t} + e^{4f_{t}} \vert \gamma_{t}\vert^{2}- e^{-4f_{t}}\vert \beta_{t}\vert^{2} \geq (2\pi)^{-1}\left( \int_{M} t f  - \int_{M} e^{-4f_{t}} \vert \beta_{t}\vert^{2}\right). $$
Under our simplifying assumptions the last term tends to zero as $t\rightarrow 0$ and we see that ${\rm deg} E_{I}\geq 0$. 
    
The proofs of items (1) and (3)  in the general case requires some more careful analysis but no fundamental difficulties. Item (2) is of a different order. Uhlenbeck and Yau write about their whole proof: {\it The technical part of the proof is quite straightforward except for one point. We obtain the [subbundles] as \lq\lq holomorphic'' in a very weak sense...obtaining enough regularity to  describe them as sheaves is more difficult.}

\

Uhlenbeck and Yau gave two largely independent treatments of this crucial difficulty. One uses complex analysis techniques. It is equivalent to show that a $L^{2}_{1}$ map into a complex Grassmannian which is a weak solution of the Cauchy-Riemann equations is meromorphic. Thus it fits into the same general realm of the regularity of weak harmonic maps we discussed in Section 4. The other uses gauge theory techniques which we will postpone
to  subsection 6.2 below.

\

The circle of ideas around this correspondence between the existence of solutions to the Hermitian-Yang-Mills equations and the algebro-geometric notion of stability has been extremely fruitful and influential in  developments in  complex differential geometry over the past four decades and we only mention a few aspects of this. In the case of a complex projective manifold, with integral K\"ahler class, the author gave alternative proofs in
 \cite{kn:SKD1}, \cite{kn:SKD2} exploiting a variational point of view. Soon after, Simpson   gave another proof which combined the variational point of view with the Uhlenbeck-Yau techniques \cite{kn:Sim}. Simpson considered a more general problem, involving a holomorphic bundle and additional fields, and there have been huge developments in that direction, one pioneer being Uhlenbeck's  student Bradlow \cite{kn:Bradlow}.
 Li and Yau extended the correspondence to the case of a general Hermitian base manifold in \cite{kn:LY}.  In place of the Uhlenbeck-Yau continity method, most subsequent work has involved the natural nonlinear heat equation---the Yang-Mills flow---associated to the existence question. (As Uhlenbeck and Yau state in \cite{kn:UY} the two methods are closely related.) Very 
complete results have been obtained, by Sibley and Wentworth \cite{kn:SW} and other authors, on the limiting behaviour of this flow  and connections with the algebro-geometric Harder-Narasimhan filtrations. 

 In a different direction
 the questions of existence of K\"ahler-Einstein metrics (in the Fano case)
 and more generally of \lq\lq extremal'' metrics turn out to fit into the same conceptual picture as the Hermitian-Yang-Mills theory and this has been the scene for much activity.  Meanwhile, on the algebraic geometry side, there have been vast extensions of the notion of stability, starting with the work of Bridgeland, and there is also much activity relating these developments to other  equations in complex differential geometry.

\subsection{Connections with small normalised energy}
 
We first discuss higher-dimensional generalisations of the small energy results Theorem 10---the gauge theory analogue of part of subsection 4.1. Then we go back to explain their relevance to the Uhlenbeck-Yau proof of Theorem 11.

The basic small-energy result for Yang-Mills connections in any dimension is just the same as in dimension $4$. For a Yang-Mills  connection over the unit ball with sufficiently small energy that energy controls all derivative of the connection, in a suitable gauge, over  an interior ball. 
 For the application to the Hermitian-Yang-Mills problem we want to consider a more general situation. 
\begin{thm}
Let $ B\subset \bC^{m}$ be the unit ball and $B'$an  interior ball and let $q$ be an exponent with $m<q<2m$. There are 
$\epsilon,\eta, C >0$ such that if   $A$ is a unitary connection over $B$  with
curvature of type $(1,1)$ and such that
\begin{equation} \Vert F(A)\Vert^{2}_{L^{2}} \leq \epsilon\ \ \ \ ,   \ \ \ \Vert \hatF\Vert_{L^{q}}\leq \eta \end{equation}
 then there is an $L^{q}$ bound on the curvature over $B'$
 
 $$  \Vert F \Vert_{L^{q}(B')}\leq C ( \Vert f\Vert_{L^{2}(B)} + \Vert \hatF \Vert_{L^{q}(B)}). $$
 \end{thm}
 
 Given this, we can apply Uhlenbeck's gauge fixing Theorem 7,
  once $F$ is sufficiently small in $L^{2}$ and $\hat{F}$ in $L^{q}$,  to get a connection form $\uA$ over $B'$ with a bound on
 $\Vert \uA\Vert_{L^{q}_{1}(B')}$.

Uhlenbeck's proof of such a result was written in the unpublished manuscript \cite{kn:UHYM}. A description of the proof, and generalisations, can be found in the recent paper \cite{kn:ChenW}.  For Yang-Mills connections (including the case of Hermitian-Yang-Mills connections, where $\hat{F}$ is a constant multiple of the identity) a proof was given by  Nakajma \cite{kn:Nak}, following the same lines as Schoen's proof  of the corresponding result for harmonic maps (Proposition 4 above). It uses a monotonicity formula going back to Price \cite{kn:Price} for the normalised Yang-Mills energy in real dimension $n$
$$  \widehat{\cE}= r^{4-n} \int_{B_{r}} \vert F\vert^{2}. $$
The proof of monotonicity, for smooth Yang-Mills connections, goes exactly as in 4.1. (In this section the reader should keep in mind that we are are working in complex dimension $m$ so the real dimension is $n=2m$, to fit in with our previous notations.)
We give a proof of Theorem 12 on the same lines as Nakajima's proof here.

\

 The first thing is to obtain a monotonicity-type property for the 
 normalised energy of connections with curvature of type $(1,1)$ and with an $L^{q}$ bound on $\hatF$.  For a connection with curvature of type $(1,1)$ we have an identity
   \begin{equation}  \vert F \vert^{2}\ {\rm vol} = \frac{1}{(m-2)!}{\rm Tr} (F^{2})\wedge \omega^{m-2} +
 \frac{1}{m} \vert \hatF\vert^{2}\ {\rm vol}.\end{equation}
   
Write the flat K\"ahler metric on the   unit ball as $\omega= d\lambda$, where the $1$-form $\lambda$ is half the contraction of $\omega$ with the radial vector field $r\partial_{r}$. Then we have
$$  \int_{B}{\rm Tr}( F^{2})\wedge\omega^{m-2}= \int_{\partial B}{\rm Tr}( F^{2})\wedge \lambda \wedge \omega^{m-3}. $$
Calculation shows that there is a pointwise bound on $\partial B$:
   $$ \frac{2}{(m-3)!} {\rm Tr}(
F^{2})\wedge \lambda \wedge \omega^{m-3} \leq \vert F\vert^{2} {\rm vol}_{\partial B}. $$
So we get
    $$  2(m-2) \int_{B}\vert F\vert^{2}\leq \int_{\partial B} \vert F\vert^{2} + \frac{2(m-2)}{m} \int_{B} \vert \hatF\vert^{2}. $$
 Applying the same argument to balls of radius $r<1$ we get the inequality for the derivative of normalised energy
   $$  \frac{d\hat{\cE}}{dr}\geq -\frac{2(m-2)}{m} r^{3-2m}\int_{B_{r}}\vert \hatF\vert^{2}. $$ 
   This gives $$\frac{d\widehat{\cE}}{dr}\geq -c r^{3-4m/q}\left(\int_{B_{r}} \vert \hatF\vert^{q}\right)^{2/q}, $$ for a suitable constant $c$. Since $q>m$ we have
$3-4m/q>-1$ and we can integrate to get, for all $r\leq 1$,
$$  \widehat{\cE}(r)\leq \widehat{\cE}(1) + c \Vert \hatF\Vert^{2}_{L^{q}}. $$

It follows from this that, in the setting of Theorem 12,  we can suppose that the normalised energy on all interior balls is as small as we please, by making $\epsilon$ and $ \eta$ suitably small.

\begin{lem} 
There are $\theta, c >0$ such that if
 $\uA$ is a connection form over the unit ball $B\subset \bC^{m}$ with curvature of type $(1,1)$, satisfying $d^{*}\uA=0$ and the boundary condition of Theorem 7, and with $\Vert \uA\Vert_{L^{m}_{1}}\leq \theta$, then the restriction of $\uA$ to the ball $B_{\tthreequart}$ satisfies
\begin{equation}  \Vert \uA \Vert_{L^{q}_{1}(B_{\tthreequart})}
\leq c \left( \Vert F \Vert_{L^{2}(B)} + \Vert \hatF\Vert_{L^{q}(B)}\right).
  \end{equation}

\end{lem}

The proof is similar to part of the proof of Theorem 7. Just as in there we have  an estimate (once $\uA$ is small in $L^{m}_{1}$):
$$  \Vert \uA\Vert_{L^{2}_{1}}\leq c \Vert F\Vert_{L^{2}}, $$
and the Sobolev embedding maps $L^{2}_{1}$ to $L^{\nu}$, with $\nu=2m/(m-1)$,
so $$   \Vert \uA\Vert_{L^{\nu}}\leq c \Vert F\Vert_{L^{2}}. $$

For a real $2$-form $\Omega$  let $\pi(\Omega)= (\Omega^{0,2}, \hat{\Omega})$.
 So the conditions on $\uA$ give
$$  \vert \pi d \uA\vert \leq \vert \hatF\vert + \vert \uA\wedge \uA\vert. $$
and we also have $d^{*}\uA=0$. The point now is that the operator
$D=d^{*}\oplus \pi d$ on $1$-forms is (overdetermined) elliptic. In fact it can be identified with
$$  \db^{*}\oplus \db: \Omega^{0,1}\rightarrow \Omega^{0,2}\oplus \Omega^{0}. $$
Let $\chi$ be a cut-off function equal to $1$ on $B_{\tthreequart}$. We have
$$  \vert D( \chi \uA)\vert \leq \vert \nabla \chi \vert \vert \uA\vert
+ \vert (\chi \uA)\wedge \uA\vert  + \vert \hatF\vert. $$
Elliptic theory gives 
$$\Vert \chi \uA\Vert_{{\nu}_{1}}\leq c \Vert D(\chi \uA)\Vert_{L^{\nu}}.$$

So $$ \Vert \chi \uA\Vert_{\nu,1}\leq c \left(\Vert\chi \uA\wedge\uA\Vert_{L^{\nu}}
+ \Vert \hatF\Vert_{L^{\nu}}+ \Vert \uA\Vert_{L^{\nu}}\right). $$
Then we can employ critical quadratic rearrangement, once $\uA$ is sufficiently small in $L^{m}_{1}$ (which implies that it is small in $L^{2m}$).
Assuming that $m\geq 3$ we have $\nu\leq m< q$ and the $L^{\nu}$ norm of $\hatF$ is dominated by a multiple of the $L^{q}$ norm. Using also the bound we have on the $L^{\nu}$ norm of $\uA$ we get 
$$ \Vert \chi \uA\Vert_{\nu,1}\leq
 c \left( \Vert F \Vert_{L^{2}}+ \Vert \hatF\Vert_{L^{q}}\right).
 $$
 Over the ball where $\chi=1$ this is an improvement on the $L^{2}_{1}$ bound we had before, since $\nu>2$. 
 Now we can repeat this process with another cut-off function, supported on the region where $\chi=1$ and equal to $1$ on $B_{\tthreequart}$. After a finite number of such steps we get a bound on the $L^{q}_{1}$ norm of $\uA$ over the $\tthreequart$ ball.

\

Note that the difference in this proof, compared with that in Theorem 7, is that the boundary condition does not combine well with the elliptic operator $D$, which is why we have to introduce cut-off functions and we only get an interior estimate.

\

 Combining this Lemma with Uhlenbeck's Theorem 7 we have  a small constant $\theta'$  such that any connection over $B$ with curvature of type $(1,1)$
 and $\Vert F\Vert_{L^{m}}\leq \theta'$ has a connection form satisfying  (57).

 \

 To prove Theorem 12, for $x$ in the unit ball $B$ let $D(x)$ be the distance to the boundary, as before. Given a connection $A$ with curvature $F$ of type $(1,1)$ over $B\subset \bC^{m}$, define $r(x)$ to be the supremum of the $r<D(x)$ such that the $L^{m}$ norm of $F$ over the $r$-ball $B_{x,r}$ is less than or equal to $ \theta'$. Define
$$ M(A)= {\rm max}_{x\in B} \frac{D(x)}{r(x)}. $$

 If we have a  bound $M(A)\leq M_{0}$ for any fixed $M_{0}$ we can apply Lemma 3 to get estimates on the $L^{q}$ norm of $F$ over any interior ball, thus proving Theorem 12. We will show that if 
 $\Vert \hatF\Vert_{L^{q}}\leq \eta$ and the normalised energy of $A$ over any interior ball is $\leq \hat{\epsilon}$,  for sufficiently small $\hat{\epsilon}$ and $\eta$ we have  $M(A)\leq 3$.

Suppose, arguing for a contradiction, that $M=M(A)>3$ and let $x_{0}$ be a point where the maximum is attained, so $D(x_{0})> 3 r(x_{0})$. Write $r=r(x_{0})$ and let $x_{1}$ be any point on the boundary of $B_{x_{0}, r}$. From the definition, $r(x_{1})\geq \ttwothree r$. By the scale invariance of the $L^{m}$ norm on $2$-forms in dimension
$2m$ we can apply Lemma 3 to the restriction of $A$ to the ball $B_{x_{0} r}$. That is, if $A'$ is the rescaled connection over the unit ball, Lemma 3 gives an $L^{q}$ bound on the curvature of $A'$ over the $\tthreequart$-ball. Since $q>m$ this implies an $L^{m}$, and by scale invariance this is an $L^{m}$ bound on the curvature of $A$ over $B_{x_{0}, 3r/4}$. This can be made as small as we please by choosing $\hat{\epsilon}, \eta$ small. In just the same way we get an $L^{m}$ bound on the curvature over $B_{x_{1}, r/2}$.
 Now take a finite number of such boundary points like $x_{1}$, whose union covers the annulus $B_{x_{0}, 5r/4}\setminus B_{x_{0}, 3r/4}$. Then by making $\hat{\epsilon}$ and $\eta $ small we can arrange that the $L^{m}$ norm of $F$ over $B_{x_{0}, 5r/4}$ is less than
$\theta'$ which is a contradiction to the definition of $r=r(x_{0})$.

\

\

  Theorem 12 leads to a global consequence for a sequence $A_{i}$ of unitary connections with curvature of type $(1,1)$
    on a fixed bundle $E$ over a compact K\"ahler manifold $X$,  with a bound on  $\Vert \hatF\Vert_{L^{q}}$ (for some $q>m$). After possibly passing to a subsequence $i'$
there is a closed set $S\subset X$ of finite
$(2m-4)$-dimensional Hausdorff measure such that $A_{i'}$ converge in $L^{q}_{1,
{\rm loc}}$ over the complement $M\setminus S $.
To see this we first  go back to equation (56). By Chern-Weil theory, the integral of ${\rm Tr}\ (F^{2})\wedge\omega^{m-2}$ is a topological invariant of the bundle, so an $L^{2}$ bound on
 $\hatF$ is equivalent to one on $F$. So since $q>2$ we have an $L^{2}$ bound on the curvatures of the $A_{i}$.
Then the proof is just as for the harmonic maps case discussed in 4.1 above , with gauge transformations constructed as  in 5.2. The only additional point is the constraint $\Vert \hat{F}\Vert_{L^{q}}\leq \eta$ in Theorem 12, for the connection over the unit ball. But this will be true for connections obtained by rescaling $A_{i}$ over sufficiently small balls in $X$, since $q>m$.

\

We return now  to the Uhlenbeck-Yau proof of Theorem 11. Item (1) of Proposition 8 shows that on the continuity path (48) the $\hatF$ satisfies a uniform
$L^{\infty}$ bound, so also and $L^{q}$ bound. So for our sequence
 $t_{i}\rightarrow 0$ the discussion above applies to the connections 
 $A_{i}=\nabla^{e^{u_{i}}}$. Thus without loss of generality there is a   weak $L^{q}_{1,{\rm loc}}$ limit convergence outside a codimension-$4$ set $S$. The limiting connection is in $L^{q}_{1}$ and has curvature of type $(1,1)$. The usual integrability theorem extends to such connections, so the connection defines a holomorphic bundle $E_{\infty}$
over $X\setminus S$. We now regard the connections $A_{i}$, over $V\setminus S$  as a convergent sequence of connections on the bundle $E_{\infty}$ and the $h_{i}$ as holomorphic bundle maps from $(E,\db_{A})$ to $(E_{\infty}, \db_{A_{i}})$. Standard arguments show that, after suitable scalings, these maps converge to a non-trivial holomorphic bundle map  from $(E,\db_{A})$ to $(E_{\infty},
\db_{A_{\infty}})$. The kernel of this map is a coherent subsheaf of $E$ over $X\setminus S$ and one sees that this is the same as the 
weak $L^{2}_{1}$ subbundle ${\rm Im} \pi$ discussed before. 

\

To sum up, Uhlenbeck and Yau use  this gauge theory argument to show that the weak subbundle is a coherent sheaf at least outside a set of real codimension $4$,  and the proof of the regularity theorem for the weak $L^{2}_{1}$ solution $\pi$
  with this extra information is  much simpler.

\

\

\subsection{Removal of codimension 4 singularities}

    In his paper \cite{kn:Tian}, Tian developed a theory for Yang-Mills connections analogous to that of Schoen and Uhlenbeck for harmonic maps. This included the notion of tangent cones at singular points. An analogue of the Schoen-Uhlenbeck small energy result in the singular case was proved by Tao and Tian in \cite{kn:TT}.
 Let $B$ be the unit ball in $\bR^{n}$ and $B'$ an interior ball. A smooth finite energy Yang-Mills connection $A$ defined on the complement 
$B\setminus \Sigma $ of a closed set $\Sigma$ of finite $(n-4)$ dimensional
Hausdorff measure is called {\it admissible}. Just as for harmonic maps, the connection is called {\it stationary} if for any vector field $v$ with compact support in the interior of $B$ generating
a $1$-parameter group of diffeomorphisms $\Phi_{t}$ we have
$$  \frac{d}{dt} \cE( \Phi_{t}^{*}(A))\vert_{t=0} = 0. $$
The monotonicity of normalised energy holds for stationary Yang-Mills connections.
 Tao and Tian's result is
\begin{thm}
There is an $\epsilon>0$ such that any stationary, admissible connection $A$ over $B$ with energy less than $\epsilon$ extends smoothly over $B'$.
\end{thm}

When $n=4$ this is equivalent to Uhlenbeck's Theorem 8 on the removal of point singularities. In the 2018 paper \cite{kn:SmU}, Smith and Uhlenbeck give a different proof of Theorem 13 which is in many ways simpler than Tao and Tian's. In fact, Smith and Uhlenbeck prove a more  general result for solutions of Yang-Mills-Higgs equations, with additional fields. In this subsection 6.3 we will discuss the Smith and Uhlenbeck proof of Theorem 13.

A consequence of Theorem 13 is that the singular set of a finite energy  stationary Yang-Mills connection has dimension strictly less than $(n-4)$. Such singular sets do arise in interesting examples. In particular, an extension by Bando and Siu of the Uhlenbeck-Yau theorem gives the existence of Hermitian -Yang-Mills connections on stable reflexive sheaves over a
 K\"ahler manifold \cite{kn:BandoSiu}. Such sheaves are vector bundles outside a singular set of complex codimension at least three and the Hermitian-Yang-Mills connection has singularities at this set.
Recent work of Chen and Sun \cite{kn:ChenSun} gives an algebro-geometric description of the tangent cone, in Tian's sense, of the Hermitian-Yang-Mills connection at a singular point.

  \

  A key component in Smith and Uhlenneck's proof is a refiment of the differential inequality (41)  for $\vert F\vert$,  where $F$ is the curvature of a Yang-Mills connection. The derivation of this used what is sometimes called the Kato inequality:  $\vert \nabla_{A} \vert F\vert\vert\leq \vert \nabla_{A} F\vert$. Similar to the calculations in 2.2 one has
$$  \vert F\vert \Delta \vert F\vert = ( \nabla_{A}^{*}\nabla_{A} F, F) +  \vert \nabla_{A} F\vert^{2} - \vert \nabla \vert F\vert\ \vert^{2}, $$ 
 and the Kato inequality gives (41).
 This can be improved via a \lq\lq refined Kato inequality'', which we now review. (See also the survey \cite{kn:kato}.)

 In general, suppose that $V$ is a Euclidean vector bundle of rank greater than $1$ over a Riemannian manifold $M$ with a compatible  covariant derivative $\nabla$ and that $s$ is a section of $V$ which does not vanish at a point $p\in M$. Then $\vert \nabla \vert s\vert\vert= \vert \nabla
s\vert$ at $p$ if and only if the image of $\nabla s$, regarded as a linear map from $TM_{p}$ to $V_{p}$, lies in the $1$-dimensional subspace spanned by
$s(p)$. In other words $\nabla s = \theta\otimes s$ for some $\theta\in T^{*}M_{p}$.
   Let $D:\Gamma(V)\rightarrow \Gamma(W)$ be a first-order
linear differential operator which is the composite of the covariant derivative with a bundle map $\sigma: T^{*}M\otimes V\rightarrow W$, so $\sigma$ is the symbol of $D$.  Alternatively, for each $\theta\in T^{*}M$ we have a $\sigma_{\theta}:V\rightarrow W$. Suppose that $D$ is overdetermined elliptic in the sense
that $\sigma_{\theta}$ is injective for all non-zero $\theta$. Then if $s$ is a non-trivial  solution of the linear equation $Ds=0$ we cannot have $\nabla s=\theta \otimes s$, since this would imply that $s$ lies in the kernel of $\sigma_{\theta}$.
It follows then from general considerations that there is some $k<1$
such that solutions of the equation $Ds=0$ satisfy the refined inequality
\begin{equation}  \vert \nabla \vert s\vert \vert \leq k \vert \nabla s\vert \end{equation} for some $k<1$.
In the case at hand, $V$ is the bundle $\Lambda^{2}\otimes \rad_{E}$, the operator  $D$ is 
 $$d_{A}\oplus d^{*}_{A}: \Omega^{2}(\rad) \rightarrow \Omega^{3}(\rad)\oplus \Omega^{1}(\rad)$$
and $s$ is the curvature $F=F_{A}$. Uhlenbeck and Smith show the constant $k$ is $\sqrt{(n-1)/n}$.  It follows then that away from the zeros of $F$ we have an improvement of the differential inequality (41) which it is convenient to  write for  $f=\vert F\vert/4$ as ,
\begin{equation} \Delta f \geq \alpha f^{-1}  \vert \nabla f\vert^{2}- f^{2}.  \end{equation}
with  $\alpha=1/(n-1)$. 
We can also write this as 
$$ \Delta f^{1-\alpha} \geq (\alpha-1) f^{2-\alpha}.  $$
Setting $\of=f^{1-\alpha}$, this is
\begin{equation}
\Delta \of\geq (\alpha-1) \  f\ \of \end{equation}
As in our previous discussion in 6.3,  these inequalities hold in a  weak sense over the zeros of $F$. The  value of $\alpha$ is not of fundamental importance, but the proof
will be simplified a bit by knowing that $\alpha<\thalf$.

  To outline the Smith and Uhlenbeck proof we begin with the 4-dimensional case, so we suppose that $A$ is a smooth Yang-Mills connection over the punctured ball $B^{4}\setminus\{0\}$ with small energy
$  \Vert F\Vert^{2}_{L^{2}}\leq \epsilon$. Our goal is to show that $F$ is in $L^{p}$ for all $p>0$. As in 6.3 (and below), once  we have this it is relatively straightforward to show that the connection extends smoothly over the origin. We divide the argument into four main steps.

\

{\bf Step 1}

\

 We claim that the function $\of$ is a weak solution of the inequality (60) over the 4-ball. In other words, if $\sigma$ is a smooth positive test function of compact support in $B^{4}$ then
\begin{equation}  \int_{B^{4}} \Delta \sigma \ \of + (1-\alpha) f \ \of  \sigma \geq 0\cdot \end{equation}
If $\chi_{\delta}$ is a standard cut-off function, with
 $\chi_{\delta}(x)=1$ for $\vert x\vert> 2\delta$ and vanishing for $\vert x\vert<\delta$,  then multiplying the inequality by $\chi_{\delta}$ and integration by parts gives
$$   \int_{B^{4}} \Delta (\chi_{\delta} \sigma)  \of +(1-\alpha) f \of \chi_{\delta}\sigma \geq 0. $$
The $k$th. derivatives of $\chi_{\delta}$ are $O(\delta^{-k})$, so we have
$$    \vert \Delta (\chi_{\delta} \sigma) - \chi_{\delta} \Delta \sigma \vert \leq K_{\sigma} \delta^{-2}, $$
for some $K_{\sigma}$ depending on $\sigma$. Thus we get
\begin{equation}     \int_{B^{4}}\chi_{\delta}\left(  \Delta \sigma\  \of +(1-\alpha) f\ \of  \sigma\right) \geq 
-K_{\sigma} \delta^{-2} \int_{B_{2\delta}} \of. \end{equation}
We know that $f$ is in $L^{2}$ so $\of$ is in $L^{2/(1-\alpha)}$ and this implies that the integral of $\of$ over the ball $B_{2\delta}$ is $O(\delta^{2+2\alpha})$.
So  the right hand side of (62) tends to $0$ with $\delta$, which establishes the claim. (Note that this step also works with $\alpha=0$.)

\

{\bf Step 2}

\

Consider the linear operator 
$$T_{f}(u)= -\Delta u - (1-\alpha) f u .  $$ 

We consider this as an operator on $L^{p}_{2}$ for $p<2$. Sobolev embedding gives  $L^{p}_{2}\rightarrow L^{r}$ for $r=2/(2-p)$. This is the same exponent as that given by H\"older's inequality for the multiplication
$$  L^{2}\times L^{r}\rightarrow L^{p},$$
so, since $f\in L^{2}$, the operator $T_{f}$ is bounded from $L^{p}_{2}$ to $L^{p}$.
Furthermore if $f$ is sufficiently small in $L^{2}$, which we can suppose, the operator $T_{f}$ is a small perturbation of $-\Delta$. Linear elliptic theory then shows that we can solve the Dirichlet problem. That is, for
$\rho\in L^{p}$ of compact support in $B^{4}$ there is a unique solution of the equation
$T_{f}(u)=\rho$ vanishing on the boundary,  and $\Vert u\Vert_{L^{p}_{2}}\leq C \Vert \rho\Vert_{L^{p}}$. 

\
\pagebreak

{\bf Step 3}

\
 
Take a cut-off function $\zeta$ of compact support in the unit ball and equal to $1$ on the half-sized ball. (This could be $\chi_{\thalf}$ but for the later discussion we prefer to use a different symbol.) Set  
\begin{equation} \rho= \Delta(\zeta \of)- \zeta \Delta \of= 2 \nabla\zeta.\nabla \of + (\Delta \zeta)  \of . \end{equation}

This is supported in an annulus on which the connection $A$ is smooth  so certainly $\rho$ is in $ L^{p}$. In fact the  standard estimates we discussed in 5.2 show that $\Vert \rho\Vert_{L^{p}}\leq C \sqrt{\epsilon}$.  Applying the linear theory, we find a function $g$ such that $T_{f}(g)=\rho$ with $\Vert g\Vert_{L^{p}_{2}}\leq C \sqrt{\epsilon}$.  By construction and the inequality (60)
$$  T_{f}( \zeta \of-g)\leq 0$$
in the weak sense. 

\

{\bf Step 4}

\

Write $h= \zeta \of-g$,  so $T_{f} h\leq 0$. 
In a case where we had $T_{f} h=0$ we would deduce from the uniqueness of the solution to the Dirichlet problem that $\zeta f=g$, so  $\of$ is in $L^{p}_{2}$ on the $\thalf$-ball and  Sobolev embedding gives $\of\in L^{r}$. Here $r$ can be made as large as as we please by taking $p$ close to $2$. Then $f$ would also in $L^{r}$ for all $r$, which was what we set out to prove.

The final, and critical, step  is a result of maximum principle type for the operator $T_{f}$, to handle the situation where we have the inequality $T_{f}h\leq 0$. 
\begin{lem}
If $f$ is sufficiently small in $L^{2}$ then for any function $h$ such that
$h\in L^{2q}$ for some $q>1$, with $h$ smooth on $B^{4}\setminus \{0\}$, vanishing on the boundary of $B^{4}$ and satisfying $T_{f}(h)\leq 0$ we have  $h\leq 0$.
\end{lem}

This gives the desired conclusion because it implies that $ \of\leq g$ over $B_{\thalf}$ and since $\of$ is positive we get $\of\in L^{r}$. 

A first step in proving  Lemma 4 is to reduce to the case when $h\geq 0$, and then show that in fact $h=0$. This can be done by replacing $h$ by ${\rm max}(h,0)$ and showing that the same differential inequality holds. This is not hard to establish if $0$ is a regular value of $h$, so the zero-set is a submanifold. In any case there are arbitrarily small positive regular values
$\tau$: replace $h$ by ${\rm max}(h,\tau)-\tau$ and take a limit as $\tau\rightarrow 0$.

For later convenience, write $q=4\beta$ and $q/2=2\beta= (1-\alpha)^{-1}=1+\gamma$.
 So $h$ is $ L^{4\beta}$.
 Let $\chi= \chi_{\delta}$ be
the cut-off function as before. The idea is to control the integral of 
$\vert\nabla (\chi h^{\beta} )\vert^{2}$.

\

 There is an identity:

$$  
  \vert \nabla (\chi h^{\beta})\vert^{2} + \frac{\beta^{2}}{\gamma} 
 \chi^{2} h^{\gamma}\  \Delta h = d^{*} J + V_{\delta} h^{2\beta}  $$
 
 where:
 \begin{itemize}
 \item  the $1$-form $J$ is a  linear combination  $$J= a_{1} h^{2\beta} \chi d\chi + a_{2} h^{2\beta-1} \chi^{2}d h; $$
\item  the function$V_{\delta}$ is a  linear combination  $$ V_{\delta}= a_{3} \vert d \chi\vert^{2} +  a_{4} \chi \ \Delta \chi. $$ 
\end{itemize}
for suitable constants $a_{i}$.
  (This is the crucial point in the proof where the fact that $\alpha>0$ is used. If we took $\alpha=0$ then $\gamma=0$ and the identity breaks down because $\gamma$ appears in the denominator.)

 \

 Integrating over $B^{4}$,  the  term from $d^{*} J$ vanishes and substituting
 the differential inequality $-\Delta h\leq (1-\alpha) f h$ gives
$$  \int \vert \nabla (\chi h^{\beta})\vert^{2}\leq \frac{(1-\alpha) \beta^{2}}{\gamma}
\int f \ ( \chi h^{\beta})^{2} + \int V_{\delta} \  h^{2\beta}. $$

In dimension 4 we have a Sobolev embedding $L^{2}_{1}\rightarrow L^{4}$ so
(using the fact that $h$ vanishes on the boundary of $B^{4}$) we get
$$  \Vert \chi h^{\beta} \Vert_{L^{4}}^{2}\leq C_{1} \Vert f \Vert_{L^{2}}
\Vert \chi h^{\beta}\Vert_{L^{4}}^{2}  + C_{2} \int V_{\delta}\  h^{2\beta}.
$$
If $\Vert f \Vert_{L^{2}}\leq (2C_{1})^{-1}$  we have
$$  \Vert \chi h^{\beta} \Vert_{L^{4}}^{2} \leq 2  C_{2}
\int V_{\delta}\  h^{2\beta}.  $$
The same kind of estimates as in Step 1 show that the right hand side tends
to $0$ as $\delta\rightarrow 0$. This uses the fact that $h\in L^{4\beta}$.

We conclude that $h=0$, as desired. Hence proving Lemma 4. 

\

Now we have shown $F\in L^{r}$ for all $r$ and this completes our discussion of  the case $n=4$.

\

\

In dimensions bigger than four, Smith and Uhlenbeck follow the same scheme but using {\it Morrey
spaces} in place of Sobolev spaces. For a function $\psi$ on $\bR^{n}$,  define the Morrey norm $\Vert\ \Vert_{\cM^{p}}$
by

$$   \Vert \psi \Vert_{\cM^{p}}^{p}= {\rm sup}_{x,r}\  r^{4-n} \int_{B^{n}_{x,r}}
\vert \psi\vert^{p}. $$

(The notation here is not standard. Smith and Uhlenbeck write $\Vert\ \Vert_{X_{k}}$ for our $\Vert\ \Vert_{\cM^{p}}$,  with $k=4/p$.)

Thus a bound on the $\cM^{2}$ Morrey norm of the curvature of a connection is the same as a bound on the normalised energy over all balls. Similar to the normalised energy, we have two properties of these norms:
\begin{enumerate}
\item For functions pulled back from $\bR^{4}$ by orthogonal  projection
$\bR^{n-4}\times \bR^{4}\rightarrow \bR^{4}$ the $\cM^{p}$ norm agrees with
the $L^{p}$ norm on $\bR^{4}$, up to a factor.
\item For $q=np/4$,  H\"older's inequality shows that the $L^{q}$ norm
  controls the $\cM^{p}$ norm,   
\end{enumerate}
 In the vein of (1), in the case when $\Sigma=\bR^{n-4}\cap B^{n}$ and the connection is pulled back from a connection over $B^{4}\setminus \{0\}$ by orthogonal projection the $n$-dimensional proof essentially reduces to that in four dimensions above. In the vein of item (2), the Morrey norm $\cM^{p}$ can be viewed for many purposes  as a slightly weakened version of the $L^{np/4}$ norm. There is an elliptic theory so that, for example, for compactly supported functions on $\bR^{n}$,
$$  \Vert \nabla^{2} \psi \Vert_{\cM^{p}}\leq C \Vert \Delta \psi\Vert_{\cM^{p}}$$ and also analogues of the Sobolev embeddings so that, for $p<2$,
\begin{equation}  \Vert \psi\Vert_{\cM^{r}} \leq C \Vert \Delta \psi\Vert_{\cM^{p}}\end{equation}
for $r=2p/(2-p)$. For $p>2$
\begin{equation}  \Vert \psi\Vert_{C^{,\mu}}\leq C \Vert \leq C \Vert \Delta
\psi\Vert_{\cM^{p}},\end{equation}
where the H\"older exponent $\mu= 2-4/p$.

The foundation of the Smith and Uhlenbeck argument is the fact, derived from monotonicity, that we can suppose the curvature $F$ is small in $\cM^{2}$-norm. 
In Steps (1), (4) we need a suitable family of cut-off functions
 $\chi_{\delta}$, equal to $1$ outside the $2\delta$-neighbourhood of the singular set $\Sigma$,  vanishing in the $\delta$-neighbourhood and with $\vert \nabla^{k} \chi_{\delta}\vert\leq c \delta^{-k}$. We need to know that the volumes of these tubular neighbourhoods are $O(\delta^{4})$. This follows from the assumption of the dimension of $\Sigma$. (In fact it appears to the author that what is required here is the assumption that $\Sigma$ has finite $(n-4)$-dimensional \lq\lq Minkowski content'' see the remark preceding Proposition 6 in 4.2 above.) With these cut-off functions in hand Step 1 works in a similar way so $\of$ is a weak solution of the inequality (60). 

For Step 2, Smith and Uhlenbeck consider the operator, for $p<2$, 
$$  T_{f}: {\cal M}^{p}_{2}\rightarrow {\cal M}^{p}$$
where ${\cal M}^{p}_{2}$ can be defined as  the space of functions $\psi$ on the ball, vanishing on the boundary,  with $\Delta \psi\in \cM^{p}$.
 More precisely, Smith and Uhlenbeck work over cubes rather than balls, which means that the boundary condition can be handled by reflection, but we will ignore this technicality. We have a Sobolev embedding $\cM^{p}_{2}\rightarrow
\cM^{r}$ with $r=2p/(2-p)$ and multiplication is defined
$$  \cM^{r}\times \cM^{2}\rightarrow \cM^{p}. $$
The upshot is that, when $f$ is small in $\cM^{2}$, the operator $T_{f}$ can be regarded as a small perturbation of the Laplacian and is invertible on these spaces. So for any $\rho\in \cM^{p}$ there is a solution $g\in \cM^{p}_{2}$ of the equation $ - \Delta g +(\alpha-1) f g =\rho$, vanishing on the boundary.

One new feature occurs in Step 3. With our cut-off function $\zeta$ we have
$T_{f}(\zeta \of) \leq \rho$ where 
$$\rho= \Delta \zeta\  \of + 2 \nabla \zeta.\nabla  \of . $$

But now the singular set $\Sigma$ can intersect the annulus on which  $\nabla \zeta$ is supported, so it is not obvious that $\rho$ is in $\cM^{p}$. That requires an estimate on the $\cM^{p}$ norm of $\nabla \of$.  Because of this the proof goes by iteration on $p$. For the first iteration we take $p=\tfourthree$, so $2p/(2-p)=4$. We begin knowing that $\vert F\vert\in {\cal M}^{2}$ (and
 is small in that norm) and after following through Steps 1-4 we get $\vert F\vert^{1-\alpha}\in {\cal M}^{\tfourthree}_{2}\subset {\cal M}^{4}$ so $\vert F\vert\in {\cal M}^{4(1-\alpha)}$ (and is small in that norm). Since $\alpha<\thalf$ this is an improvement.
Repeating the process sufficiently many times gives $\vert F\vert\in \cM^{p}$ for all $p$, or equivalently $F$ is in $L^{p}$ for all $p$.  

\

We now return to the problem of estimating $\rho$ for the first iteration, where we want to bound $\Vert \nabla \of\Vert_{\cM^{\tfourthree}}$. 
For this, Smith and Uhlenbeck go back to the inequality (58):
$$  - \Delta f  +\alpha \frac{\vert \nabla f \vert^{2}}{f}  \leq   \ f^{2}. $$
 Let $\tzeta$ be another cut off function supported in an interior $r$-ball $B_{r,x}$ and equal to $1$ on the half-sized ball. Multiplying the inequality by $\zeta$ and integrating gives
$$  \alpha \int_{B_{r,x}} \tzeta \frac{\vert \nabla f\vert^{2}}{f} \leq \int_{B_{r,x}}
(\Delta \tzeta) f + \tzeta f^{2}. $$
Knowing that $f\in \cM^{2}$ the right hand side is easily shown to be
$O(r^{n-4})$, so we get an $\cM^{2}$ bound on 
 $ f^{-\thalf}\vert \nabla f\vert$. Now write
 $$ \vert \nabla \of\vert = \vert \nabla f^{1-\alpha}\vert  = (1-\alpha) f^{\thalf-\alpha} \left( \frac{\vert \nabla f\vert}{f^{\thalf}}\right).
$$ The $\cM^{2}$ bound on $f$ gives a $\cM^{4/(1-2\alpha)}$ bound on
$ f^{\thalf-\alpha}$. Then the multiplication:
  $$  \cM^{4/(1-2\alpha)} \times \cM^{2}\rightarrow \cM^{4/(3-2\alpha)}$$ gives a bound on $\nabla \of$ in $\cM^{4/(3-2\alpha)}$ and therefore in
$\cM^{\tfourthree}$.

For the crucial step 4, Smith and Uhlenbeck establish an analogue of Lemma 4 in Morrey spaces but we will pass over that and move on to outline their argument for the removal of singularities, given a bound on the $L^{p}$ norm of the curvature for large $p$.

Let $x_{0}$ be a point in the complement of the singular set $\Sigma$. Define the {\it shadow} of $\Sigma$ to be the set of $x$ such that for some
$t\in (0,1]$ the  point $t x+(1-t)x_{0}$ lies in $\Sigma$. Then for a connection defined on the complement of $\Sigma$ the exponential gauge construction that we discussed in   5.2,  using rays emanating from $x_{0}$, defines a connection form $\uA$ over the complement of the shadow. The fact that $\Sigma$ has codimension at least $4$ implies that the shadow has codimension at least $3$. Uhlenbeck and Smith show that $\uA$ and $d\uA$ are in $L^{p}$ for all $p$, where the latter is interpreted a a distribution. By working over a small ball centred at $x_{0}$ and rescaling one can assume that the
$L^{p}$ norms are small. Then the implicit function theorem can be applied to choose a new gauge in which the Coulomb condition is satisfied, and in this gauge elliptic regularity shows that connection is smooth over the small ball. The fact that $d\uA$ is in $L^{p}$, as a distribution, uses crucially the codimension condition. By contrast suppose we had a singular set of codimension $2$, so we can have a non-trivial flat connection on the complement. Then
the shadow would have codimension $1$, the connection form $\uA$ would have a discontinuity across the shadow and the distribution $d\uA$ would have a singular component supported on the shadow. Thus the same argument would not work in that case, in agreement with the fact that the singularity is not removable.

\section{Harmonic maps to Lie groups}
\subsection{Harmonic maps, flat connections and loop groups}
In this section we discuss harmonic maps from  surfaces to unitary groups (many of the constructions extend to other compact Lie groups).
 The main topic is an important paper of Uhlenbeck \cite{kn:ULoop} which, among other things, gives semi-explicit constructions for the 
general solution in terms of geometric data when the surface is the Riemann sphere. There are many connections with gauge theory, related to the Coulomb gauge condition. The results fit into many large circles of ideas. The fundamental observation of the zero-curvature interpretation of the harmonic map equations is attributed by Uhlenbeck to  Pohlmeyer \cite{kn:Pint}, within the integrable systems literature. The equations are, as we explain below in this subsection 7.1 and in 7.2, related to Hitchin's equations on Riemann surfaces and to the Yang-Mills instanton equations in four dimensions.
Thence the integrable nature of the harmonic map equations can be related to \lq\lq twistor'' constructions. We will only  touch on a small fraction of all these ideas here and there are important parts of Uhlenbeck's paper that we do not discuss.
 In subsection 7.3 we go back to analysis and  results of H\'elein on regularity questions,  which use somewhat related ideas.

\

For any simply connected manifold $M$ the harmonic map equations for a map
from $M$ to the group $U(r)$ can be formulated as a system of three equations for a pair $(A,\psi)$
where $A$ is a $U(r)$ connection on a bundle $E$ over $M$ (in fact the trivial bundle) and $\psi \in \Omega^{1}(\rad_{E})$. These equations are
\begin{equation}  d^{*}_{A} \psi=0\ \ \ \ \ \ d_{A} \psi=0 \ \ \ \ \  \ F(A)+ \psi\wedge \psi = 0. \end{equation}
Given a solution $(A,\psi)$ to (66) the last two equations state that
$F(A\pm \psi)=0$, so $A_{1}= A+\psi$ and $A_{-1}= A-\psi$ are flat connections and hence gauge equivalent, since $M$ is simply connected. We always
 have $d^{*}_{A+\psi}\psi= d^{*}_{A}\psi$ and so the first equation states that $d^{*}_{A_{1}}\psi=0$. Choose a trivialisation in which $A_{1}$ is the product connections: {\it i.e.} the connection $1$-form $\uA_{1}$ is zero. If $\uA_{-1}$ is the connection $1$-form of  $A_{-1}$ in this trivialisation, the first equation states that $d^{*}\uA_{-1}=0$. The fact that $A_{-1}$ is flat means that $\uA_{-1}= -dg g^{-1} $ for some $g:M\rightarrow U(r)$ and then $f$ satisfies the harmonic map equation $d^{*} ( dg g^{-1})=0$.
 Conversely, given such a harmonic map $g$, let $A$ be the connection on the trivial bundle defined by the connection $1$-form $-\thalf  dg g^{-1} $ and let  $\psi= \thalf dg g^{-1}$. Then $A\pm \psi$ are flat connections and  we get a solution of the equations (66).

If we change the sign in the third equation in (66) to $F(A)-a\wedge a=0$ we have a similar discussion, with flat $GL(n,\bC)$-connections $A\pm i \psi$ and we get a correspondence with harmonic maps to the dual symmetric space $GL(r,\bC)/U(r)$.

In the case when $M$ is a Riemann surface there is an additional symmetry between the first two equations in (66). Recall that $\rad_{E}$ is the real subbbundle of ${\rm End}\  E$ consisting of the skew adjoint endomorphisms.
Over a Riemann surface $M$ we can write $\psi=\Phi-\Phi^{*}$ where
$\Phi$ is a $(1,0)$ form with values in ${\rm End} E$. Written in terms of $\Phi$, the equations (66) are 
\begin{equation}    \db_{A}\Phi=0 \ \  ,\ \  F(A) -( \Phi \wedge \Phi^{*}+\Phi^{*}\wedge \Phi) =0, \end{equation}
where $\db_{A}$ is the coupled $\db$-operator, which is complex linear.
Clearly if $(A,\Phi)$ is a solution to (67) and $\lambda$ is a complex number of modulus $1$ then $(A,\lambda^{-1} \Phi)$ is also a solution. The equations with the opposite sign of the quadratic term in $\Phi$ are {\it Hitchin's equations} \cite{kn:Hitchin}. In the  Riemann surface case, one solution of the equations (66) gives a circle of solutions corresponding to $\lambda^{-1} \Phi$ where $\lambda$ is a complex number with $\vert \lambda\vert=1$; so we have we a flat connection $A_{\lambda}$. More generally, for any $\lambda\in \bC^{*}$ we can define
a flat $GL(r,\bC)$ connection
$$A_{\lambda}= A+  \lambda^{-1} \Phi - \lambda \Phi^{*}. $$

 Fix a trivialisation of the bundle in which $A_{1}$ is the product connection. For each $\lambda\in\bC^{*}$ we have
a map $G_{\lambda}:M\rightarrow GL(n,\bC)$ 
such that $A_{\lambda}= -dG_{\lambda} G_{\lambda}^{-1}$.
If we fix a basepoint $p\in M$ then $G_{\lambda}$ is determined uniquely by the condition that $G_{\lambda}(p)=1$. Define antiholomorphic involutions of $\bC$ and $GL(r,\bC)$ by 
$$  \sigma(\lambda) = \overline{\lambda}^{-1}  \ \ \ \ \ ,\ \ \ \  \sigma(g)=\left( g^{*}\right)^{-1}. $$
Then the family $G_{\lambda}$ has the equivariance property
$$     G_{\sigma\lambda}= \sigma(G_{\lambda})$$
and in particular $G_{\lambda}$ maps into $ U(r)$when $\vert \lambda\vert=1$.
The map $G_{-1}$ is the harmonic map $g$ which we associated to the solution of the equations (66).

 We can also regard this family as a single map $G:M\times \bC^{*}\rightarrow GL(n,\bC)$. So any  harmonic map $g:M\rightarrow U(r)$ defines  an  \lq\lq extended map'' $G$.

 Uhlenbeck's first result is a characterisation of these extended maps.
\begin{prop} For  a simply connected Riemann surface $M$ there is a (1-1) correspondence between harmonic maps $g:M\rightarrow U(r)$ and families 
$G_{\lambda}:M\rightarrow GL(r,\bC)$, for $\lambda\in \bC^{*}$ such that
\begin{itemize}
\item $G_{\lambda}$ is holomorphic in $\lambda$
\item $G_{\sigma \lambda}= \sigma (G_{\lambda})$,
\item 
$$  \frac{1}{1-\lambda^{-1}} G_{\lambda}^{-1} \partial G_{\lambda} $$
is independent of $\lambda$.
\end{itemize}
\end{prop}

Given such a family $G_{\lambda}$, we define the matrix-valued $(1,0)$ form on $M$ by $\Phi= -(1-\lambda^{-1})^{-1}  G_{\lambda}^{-1} \partial G_{\lambda}$ and a connection on the trivial bundle with connection 1-form $\Phi^{*}-\Phi$
to get back to a solution of (67).

\

Let $\Omega U(r)$ be the based loop group, of smooth maps 
$\gamma: S^{1}\rightarrow U(r)$ with $\gamma(1)= 1$ and let $\epsilon:\Omega U(r) \rightarrow U(r)$ be evaluation at $-1$.  The restriction of  the family $G_{\lambda}$ to the unit circle can be regarded as a map
$\tG: M\rightarrow \Omega U(n)$ so we have a canonical lift of a harmonic map $g$ over the evaluation map $\epsilon$:
  \begin{equation}   M\stackrel{\tG}{\rightarrow} \Omega U(r) \stackrel{\epsilon}{\rightarrow} U(r). \end{equation}
   
   In \cite{kn:Segal} Segal gives an interpretation of the extended map conditions in Proposition 9 in terms of the geometry of the loop group $\Omega U(r)$. This loop space has an infinite-dimensional K\"ahler structure, preserved by left multiplication of the group. So the complex structure is determined by a complex structure on the tangent space $T$ at the identity. This tangent space $T$ consists of maps $\xi$ from the circle $\vert \lambda\vert =1$ to skew adjoint matrices with $\xi(1)=0$. Such a map has a  Fourier series
  $$\xi(\lambda)= \sum_{k=-\infty}^{\infty} a_{k} \lambda^{k} $$
(for matrix-valued coefficients $a_{k}$) with $\sum a_{k}=0$ and $a_{k}= - a_{-k}^{*}$. Thus the vector space $T$  is identified with the set of  rapidly-decreasing sequences $(a_{-1}, a_{-2}, \dots)$ of complex matrices $a_{k}$. The complex structure on $T$ is the obvious one defined by the usual complex structure on these matrix coefficents. (At the group level, this corresponds to the identification of $\Omega U(r)$ with the quotient of the loops in $ GL(r,\bC)$ by the subgroup of loops which extend holomorphically over the disc, which has a visible complex structure.)
\

Let $V$ be the subspace of $T$, of complex dimension $r^{2}$,
 corresponding to sequences $(a_{-1},0,0\dots)$. Extend this by left translation to a subbundle $\underline{V}$ of the tangent space of $\Omega U(r)$, Segal's formulation of Uhlenbeck's correspondence is
\begin{thm}  For  a simply connected Riemann surface $M$ there is a (1-1)
correspondence between harmonic maps from $M$ to $U(r)$ and holomorphic maps
$\tG:M\rightarrow \Omega U(r)$  whose derivative at each point maps into $\underline{V}$. The harmonic map corresponding to $\tG$ is $\epsilon\circ \tG$.
\end{thm}

 Given $\tG:M\rightarrow \Omega U(r)$, write $G_{\lambda}:M\rightarrow U(r)$ for the corresponding family of maps, with $\lambda\in S^{1}$. The condition that $\tG$ is holomorphic is equivalent to saying that,  at each point in $M$, the $(0,1)$ form  $G_{\lambda}^{-1} \db G_{\lambda}$ extends
 holomorphically over the unit disc, say
  $$  G_{\lambda}^{-1} \db G_{\lambda}= \sum_{k\geq 0} b_{k} \lambda^{k}. $$ Similarly, the condition that the derivative maps into $\underline{V}$ is equivalent to saying that 
$$G_{\lambda}^{-1} \partial G_{\lambda}= \alpha (1-\lambda^{-1}) +Z(\lambda)$$ for some constant matrix $\alpha$ and  $Z(\lambda)$ holomorphic over  the disc with $Z(0)=0$. The condition that the $G_{\lambda}$ map into $U(r)$  implies that $$G_{\lambda}^{-1}\partial G_{\lambda}= - (G_{\lambda}^{-1} \db G_{ \lambda})^{*}= \sum_{k\geq 0} b_{k}^{*} \lambda^{-k}. $$ So  $Z$ is zero and $G_{\lambda}^{-1} \partial G_{\lambda}= \alpha (\lambda^{-1}-1)$.
 Thus $ (\lambda^{-1}-1)^{-1}G_{\lambda}^{-1} \partial G_{\lambda}$ is independent of $\lambda$ and we can reconstruct $(A,\Phi)$ just as before. (This also shows that the family $G_{\lambda}$ extends to $\lambda\in \bC^{*}$.)

\

  The preceding discussion is essentially local in the Riemann surface $M$. and we now focus on the case when $M$ is the Riemann sphere. In that case Uhlenbeck proves an important finiteness result, that the extended maps $G_{\lambda}$ have a finite Laurent series. Here it become convenient
 to drop the normalisation of the extended map $G$ using a base point $p\in M$ (but keeping the other conditions in Proposition 9). Then the extended map can be chosen of the form
\begin{equation}   G_{\lambda}= \sum_{k=0}^{n} T_{k}\lambda^{k}, \end{equation}
where the $T_{k}$ are matrix-valued functions on $M$. 
 Uhlenbeck calls the least possible number $n$  the {\it uniton number}. (If $n=0$ the map is a constant.)  One of the most important ideas in \cite{kn:ULoop} is a construction called \lq\lq uniton addition''---a form of B\"acklund transformation, taking one solution to another---which we will describe in the next subsection.
Uhlenbeck showed that all solutions can be obtained by iterating this construction,  and one of her main results is:
\begin{thm} Let $g:S^{2}\rightarrow U(r)$ be a harmonic map.
\begin{itemize}\item
The uniton number $n$ of $g$  is  strictly less than $r$.
\item There is a unique harmonic map $\underline{g}:S^{2}\rightarrow U(r)$  of uniton number $n-1$ such that $g$ is obtained from $\underline{g}$ by the operation of uniton addition.
\end{itemize}
\end{thm}
In other words any harmonic map from $S^{2}$ to $U(r)$ is constructed in a canonical way by repeating the uniton addition construction at most $(r-1)$ times. The maps with uniton number $1$ are  the holomorphic maps from $S^{2}$ to the Grassmann manifolds $Gr_{k}(\bC^{r})$ of $k$-dimensional subspaces of $\bC^{r}$ (and, of course, translates by left and right multiplication in $U(r)$).  These Grassmann manifolds are isometrically embedded in $U(r)$ by the map which takes a subspace $W\subset \bC^{r}$ to the unitary map $R_{W}=\pi- \pi^{\perp}$, where $\pi$ is projection to $W$ and $\pi^{\perp}$ to the orthogonal complement.

This theorem puts the problem of describing the holomorphic maps into the realm of geometry and Uhlenbeck (and subsequent authors) obtained a variety of specific results. These include a derivation of the description by Eells and Wood in \cite{kn:EW} of all harmonic maps  from $S^{2}$ to complex projective space $\bC\bP^{n}$. These maps have uniton number $2$. In Segal's treatment \cite{kn:Segal} he shows that when $M=S^{2}$ the holomorphic maps $\tG$ map into explicit finite-dimensional complex submanifolds---generalised flag manifolds---of the loop group. The problem then comes down to understanding the horizontality condition for these complex curves.

   \subsection{Uniton addition and instantons on $\bR^{2,2}$}

   Uniton addition can be defined over any Riemann surface $M$. It changes a   harmonic map $g:M\rightarrow U(r)$ to a new one of the form
$g'=g \rho $,  defined using the group structure on $U(r)$, where $\rho:M\rightarrow U(r)$ maps into the space of reflections. That is, we have a map $W:M\rightarrow {\rm Gr}_{k}(\bC^{r})$ for some $k$ and $\rho(z)= R_{W(z)}$ in the notation of the previous subsection. The remarkable thing is that finding the maps $W$ for which $g \rho$ is harmonic essentially involves solving only linear  PDE. Clearly if we perform the construction again, starting with $g'$ and using the map $-\rho$, we recover $g$.

We go back to the data $(A,\Phi)$ where $A$ is a connection on a bundle $E$ and $\Phi\in \Omega^{1,0}({\rm End E})$, satisfying equation (67). Using the flat trivialisation for the connection $A_{1}= A+\Phi-\Phi^{*}$ we can regard a map $W:M\rightarrow {\rm Gr}_{k}(\bC^{r})$ as a subbundle $E_{I}\subset E$.
\begin{lem}
If $E_{I}$ is a holomorphic subbundle of $E$ with respect to the holomorphic structure defined by $\db_{A}$ which is preserved by $\Phi$ in that
$\Phi E_{I}\subset \Omega^{1,0}(E_{I})$ then the corresponding map $g \rho$ is harmonic.
\end{lem}

The equation $\db_{A}\Phi=0$ is the integrability condition for being able to find such a subbundle locally. The lemma can be proved by direct calculation  but we will take a more roundabout route which puts the construction in a wider context. However we postpone that and first discuss an interesting point of view from the work of Valli \cite{kn:Valli}, involving the harmonic maps energy.

The homotopy group $\pi_{2}(U(r))$ vanishes so there is no obvious topological
invariant of a harmonic map from $S^{2}$ to $U(r)$. However $\pi_{2}(\Omega U(r))= \pi_{3}(U(r))=\bZ$
so there is an integer degree of the map $\tG:S^{2}\rightarrow \Omega U(r)$.
As Segal shows in \cite{kn:Segal}, the energy of the harmonic map is equal to twice the
degree. In fact this is true locally, in that the energy density of $g$ is
equal to twice the pull-back by $\tG$ of a standard closed $2$-form on $\Omega U(r)$ representing
the generator of $H^{2}(\Omega U(r))$.  Valli obtained a related formula in \cite{kn:Valli}. Suppose that $g'$ is obtained from $g$ by uniton addition as above, with a subbundle $E_{I}$.  Then Valli showed that
$$   E(g')= E(g) - 2 {\rm deg} \ (E_{I}). $$
where the degree ${\rm deg}\ E_{I}$  is the first Chern class, regarded as an integer. So
if ${\rm deg}\ (E_{I})>0$ then the energy of $g'$ is strictly
smaller. As Valli explained, the existence of a $\Phi$-invariant subbundle $E_{I}$ of positive degree is a simple consequence of the algebro-geometric classification of holomorphic vector bundles over the Riemann sphere. Applying this repeatedly gives a proof of a slightly weaker form of Uhlenbeck's Theorem 15,  since we can keep on performing these uniton additions until the energy is $0$ and we have a constant map. 

 There is a intriguing connection here with notions of stability, like those we encountered for the Uhlenbeck-Yau theorem in 6.1.  
Recall that Hitchin's equations for $(A,\phi)$ are obtained by changing the
sign in (67). A pair $(E,\phi)$ consisting of a  holomorphic bundle $E$ (with
$c_{1}(E)=0$) over any compact Riemann surface $M$ and a holomorphic $\phi\in
\Omega^{1,0}({\rm End}V)$ is called stable if  any $\phi$-invariant holomorphic
subbundle has strictly negative degree. Hitchin showed that if $(E,\phi)$ is stable
there is a corresponding solutions of his equations. This is analogous to the Unhlenbeck-Yau
Theorem 11 and involves solving a similar  PDE for a metric on the bundle $E$. When $M=S^{2}$,  a solution to Hitchin's equations would give a harmonic map to
the noncompact symmetric space $GL(n,\bC)/U(n)$ and it is easy to see that
these are all constant, so there are no stable pairs. If we start with an arbitrary metric and solve a natural gradient flow equation or use a variant of the Uhlenbeck-Yau continuity path the solutions will diverge in the manner we discussed in 6.1 corresponding to a subbundle of $E$ and in fact this will be a $\phi$-invariant holomorphic subbundle of strictly positive degree.
So attempting to solve the equations with the  the reversed sign  tells us how to build the general solution of (67) by repeated uniton addition.

Before returning to the proof of Lemma 5 we make a further digression to discuss an aspect of the equations (67) which were explored in  another paper \cite{kn:ULoop2} of Uhlenbeck. In 5.3 we mentioned the  Yang-Mills instanton equation in four dimensions. These are $F^{+}(\bA)=0$, where $F^{+}$ denotes the self-dual part of the curvature. Consider connections $\bA$ over $\bR^{4}$ which are translation-invariant in two directions. In terms of coordinates
$(u_{1}, u_{2}, v_{1}, v_{2})$ we can write such a connection as
$$  \bA= A + \Phi_{1} dv_{1} + \Phi_{2} dv_{2}, $$
where $A$ is a connection on a bundle $E$ over $\bR^{2}$ (with coordinates $u_{1}, u_{2}$), lifted to $\bR^{4}$ by projection, and $\Phi_{1}, \Phi_{2}$ are sections of $\rad_{E}$ over $\bR^{2}$. 
Then the instanton equation for $\bA$ becomes Hitchin's equations for $(A,\Phi)$ where $\Phi= (\Phi_{1} + i\Phi_{2})(du_{1}+ idu_{2})$. Now, as explained in \cite{kn:ULoop2}, change the signature of the metric in four dimensions to get $\bR^{2,2}$ with metric $du_{1}^{2}+ du_{2}^{2}-dv_{1}^{2}-dv_{2}^{2}$.
The notion of anti-self-duality makes sense and the instanton equation 
for translation-invariant solutions become
(67).  From this point of view, uniton addition is a 
 special case of a more
general construction for solutions of the instanton equation over $\bR^{2,2}$
(related to the PhD thesis of Uhlenbeck's student Crane \cite{kn:Crane}). It is easiest to explain this first in a complexified setting, with the coordinates
$(x_{1}, x_{2}, \xi_{1}, \xi_{2})$ on $\bC^{4}$ and the quadratic form
$dx_{1}d\xi_{1}+dx_{2}d\xi_{2}$. The anti-self-duality equations for a $2$-form

$$ F=F_{x_{1} x_{2}} dx_{1} dx_{2} + F_{\xi_{1} \xi_{2}} d\xi_{1} d\xi_{2}+\sum
F_{x_{i}\xi_{j}} dx_{i} d\xi_{j} $$
 are
\begin{equation} F_{x_{1} x_{2}}=0 \ \ \ , \ \ \ 
F_{\xi_{1} \xi_{2}}=0\ \ \ ,\ \ \ F_{x_{1}\xi_{1}}+ F_{x_{2}\xi_{2}}=0.
 \end{equation}
 Let $L,R$ be holomorphic $1$-forms on $\bC^{4}$ of the shape
 \begin{equation}  L= L_{\xi_{1}} d\xi_{1} + L_{\xi_{2}}d\xi_{2}\ \ \ ,\
\ \ R= R_{x_{1}} dx_{1} + R_{x_{2}}dx_{2}. \end{equation}
Define
$$  \tL = -L_{\xi_{2}} dx_{1} + L_{\xi_{1}} dx_{2}\ \ \ \tR=R_{x_{2}} d\xi_{1}
- R_{x_{1}} d\xi_{2}. $$
Then $L\wedge R$ and $\tL\wedge \tR$ have the same  self-dual component
\begin{equation} \left( L\wedge R\right)^{+}= 
\left( \tL\wedge \tR\right)^{+}= \frac{1}{2}\left(
 L_{\xi_{1}}R_{x_{2}} + L_{\xi_{2}}R_{x_{1}}\right) (d\xi_{1}dx_{1}+ d\xi_{2}
dx_{2}). \end{equation}
 Write $d^{+}$ for the self-dual part of the exterior derivative. The
equation $d^{+}L=0$ has two components
 $$  L_{\xi_{1}, x_{1}}+ L_{\xi_{2},x_{2}}=0\ \ \ ,\ \ \ L_{\xi_{1}, \xi_{2}}-
L_{\xi_{2},\xi_{1}},$$
where commas denote partial derivatives.
This is the same as the equation $d^{+} \tL=0$ and similarly for $R,\tR$.

Now consider a holomorphic connection $\nabla$ on a holomorphic bundle $V$
over a domain in $\bC^{4}$ with $V=V_{I}\oplus V_{II}$. So the connection
is given by connections on $V_{I}, V_{II}$ and second fundamental forms $L,R$
which are holomorphic $1$-forms with values in ${\rm Hom}(V_{I}, V_{II}),
{\rm Hom}(V_{II},V_{I})$ respectively. In the notation we used in 6.1

\begin{equation} \nabla= \left(\begin{array}{cc} \nabla_{I} & R\\ L & \nabla_{II}\end{array}\right)
, \end{equation}
and the curvature is

\begin{equation} \left(\begin{array}{cc} F_{I} + R\wedge L & d_{I, II}
 R\\d_{I,
II}L & F_{II}+L\wedge R \end{array}\right)
\end{equation}

 The anti-self duality equation for
the connection on $V$ has four components
$$   F^{+}_{I} + ( R\wedge L )^{+}=0\ ,\   d_{I, II}^{+}L=0\ ,\  d^{+}_{I, II} R=0\ ,\  F^{+}_{II}+
( L\wedge R )^{+}=0. $$
Define $\tL,\tR$ by the same formulae as above, extended to bundle-valued
forms, and use $\tL,\tR$ as second fundamental forms defining a new connection
$\tilde{\nabla}$ on $V$, with the same connections on $V_{I}, V_{II}$. That is:
\begin{equation} \tilde{\nabla}= \left(\begin{array}{cc} \nabla_{I} & \tR\\ \tL & \nabla_{II}\end{array}\right).
 \end{equation} The
formulae above---extended to bundle valued forms in an obvious way---show
that $\tilde{\nabla}$ is an anti-self-dual connection if and only if $\nabla$
is.  

Now we look at the real forms of this construction, dealing with  unitary connections.
We would like an anti-linear involution $\sigma$ of $\bC^{4}$ such that restricted
to the fixed points of $\sigma$ the forms satisfy $R=-L^{*}, \tilde{R}=-\tilde{L}^{*}$.
Change notation to write $x_{1}=z, \xi_{1}=\overline{z}, x_{2}=w,\xi_{2}=
-\overline{w}$ and let $\sigma$ be defined by complex conjugation as indicated,
so the fixed points of $\sigma$ are the set where $z+\overline{z}, w+\overline{w}$
are real. On this set $z,w$ become complex coordinates and the metric is
$dzd\overline{z}- dwd\overline{w}$. So the metric has signature $(2,2)$ and matches up with our previous discussion when we take $z=u_{1}+ iu_{2}, w=v_{1}+i v_{2}$.
In the complex coordinates $z,w$ we have
$$  L=L_{\overline{z}} d\oz + L_{\ow} d\ow\ \ ,\ \ R=R_{z} dz+ R_{w} dw,
$$
and $$\tL=R_{\ow} dz + R_{\oz} dw\ \ ,\ \  \tR= R_{w} d\oz + R_{z} d\ow.$$
Thus $R=-L^{*}$ if and only if $\tR=-\tL^{*}$. The upshot is the following.
Suppose that $\nabla$ is a unitary anti-self-dual connection on a bundle $V$ over
 a domain $\Omega$ in $\bR^{2,2}$ and
  identify $\bR^{2,2}$ with $\bC^{2}$ as above. Then $\nabla$ defines a holomorphic
structure on $V$. Suppose that $V_{I}$ is a holomorphic
subbundle of $V\rightarrow \Omega$. The construction above produces a new
unitary anti-self-dual connection $\tilde{\nabla}$ on $V$, but now the orthogonal
complement $V_{I}^{\perp}$ is a holomorphic subbundle, with respect to the
homolomorphic structure defined by $\tilde{\nabla}$.

 When restricted to translation-invariant
connections this becomes Uhlenbeck's uniton-addition construction. In that case the second fundamental form of a subbundle $E_{I}$ over $\bR^{2}$ has a component
$\beta\in \Omega^{0,1}({\rm Hom}(E_{II}, E_{I}))$ and $\Phi$ has a component
$\Phi_{I , II}\in \Omega^{1,0}( {\rm Hom}(E_{II}, E_{I}))$. The construction takes $\Phi_{I , II}^{*}\in \Omega^{0,1}({\rm Hom}(E_{I}, E_{II}))$ to build the new second fundamental form and takes $\beta^{*}\in  \Omega^{1,0}({\rm Hom}(E_{I}, E_{II}))$ to build the new $\Phi$-field.

There is also a \lq\lq twistor'' description of this construction. The Ward correspondence relates instantons on a domain $U\subset \bR^{2,2}$ to holomorphic bundles on a three dimensional twistor space $Z$. The choice of a compatible complex structure on $\bR^{2,2}$ gives a complex surface $D\subset Z$ and the opposite structure another $D'\subset Z$. In complex geometry there is a general construction, sometimes called the Hecke transform. Let $\cE$ be a holomorphic bundle over a complex manifold $\cZ$ and $\cD$ be a hypersurface in $\cZ$. Suppose given a holomorphic subbundle $\cF$ of the restriction $\cE\vert_{\cD}$. Then we define a new holomorphic bundle $\cE'\rightarrow \cZ$ whose local holomorphic sections correspond to sections of $\cE$ which lie in $\cF$ when restricted to $\cD$. The twistor description of the construction is to apply this
 transform to the bundle over the twistor space $Z$, using subbundles over $D, D'$.

   \subsection{Weak solutions to the harmonic map equation on surfaces}
   
   A question left open in our discussion of the regularity theory for
    harmonic maps in Sections 3 and 4 is whether a weakly harmonic map from a surface---without additional minimising or stationary hypotheses---is smooth. In the case when the target space $N$ has a large isometry group this smoothness was established by H\'elein in \cite{kn:Hel1}. The proof was later extended to general target spaces \cite{kn:Hel2} but here we just consider the case of a unitary group, where the proof is particularly simple. The proof depends on the particular structure of the equations, in a similar vein to the preceding discussion in this section.

 \begin{prop} Let $g:D\rightarrow U(r) $ be a weak solution to the harmonic map equations on  the unit disc $D$ in $\bC$ with
derivative in $L^{2}$. Then $g$ is continuous.\end{prop}
 Once the continuity is established smoothness follows from general theory as in  \cite{kn:Hild}.

The proof of Proposition 10 depends on the following result of H. Wente \cite{kn:Wente}.
\begin{prop}
Suppose that $u_{1},\dots, u_{k}, v_{1},\dots, v_{k} $ are functions on the disc $D\subset \bC$ with
 derivatives in $L^{2}$. If $\phi$ satisfies $\Delta \phi= *\sum(du_{i}\wedge dv_{i})$ then $\phi$ is continuous.
\end{prop}
Of course $*\sum (du_{i}\wedge dv_{i})$ is in $L^{1}$, but in dimension $2$ it is not the case that any function
$\phi$ with $\Delta\phi$ in $L^{1}$ is  continuous. That is  the point of the result. This is
the borderline situation; if $\Delta\phi$ is in $L^{p}$ for some $p>1$
then $\phi$  is continuous.

\

Let us now see how Wente's result implies Proposition 10. The 
 harmonic map equation is $d^{*}(dg g^{-1})=0$. So $d* (dg g^{-1})=0$ and we can write $*dg g^{-1}= d \lambda$ for a matrix-valued function $\lambda$.
 Now go back to the equation $d^{*}(dg g^{-1})=0$ which is
 $$  \sum_{i=1}^{2} \frac{\partial}{\partial x_{i}}\left( \frac{\partial g}{\partial x_{i}} g^{-1}\right)=0. $$
Expanding out and multiplying on the right by $g$ we get
$$  \Delta g = \sum_{i=1}^{2}\frac{\partial g}{\partial x_{i}} g ^{-1} \frac{\partial
g}{\partial x_{i}}g = \frac{\partial \lambda}{\partial x_{2}}
\frac{\partial
g}{\partial x_{1}} -\frac{\partial \lambda}{\partial x_{1}}\frac{\partial
g}{\partial x_{2}}.$$
In other words
$$  \Delta g= * ( d\lambda\wedge dg). $$
 Thus each matrix entry $g_{ab}$ of $g$ satisfies
 $$  \Delta g_{ab}= \sum_{c}  d\lambda_{ac}\wedge dg_{cb} . $$
 Since $g$ takes values in the unitary group the operation of multiplication by $g$ preserves the standard norm on matrices. So the
 matrix-valued $1$-forms $d\lambda$ and $dg$ are in $L^{2}$ and hence the same for each entry. So we are in the situation considered in Proposition 11 and the $g_{ab}$ are continuous.

\

The proof of Wente's result, Proposition 11, has some common features with the Sacks-Uhlenbeck proof of removal of singularities in Section 3. Recall that the Green's function on $\bR^{2}$ is $(2\pi)^{-1}\log r$. It suffices to prove that there is a constant $C$ such that
\begin{equation}  \vert \int_{D} \log r *(du\wedge dv) \vert \leq C \Vert du\Vert_{L^{2}}\Vert dv\Vert_{L^{2}}\end{equation}
for all smooth compactly supported functions $u, v$ on $D$. Transfer to the cylindrical picture with coordinates $(s,\theta)$ where $r=e^{s}$.  So $u,v$ are now defined on the half-cylinder $(-\infty,0]\times S^{1}$, vanishing on $\{0\}\times S^{1}$,  and the
  left hand side of (76) is the modulus of
$$  I= \int_{(-\infty,0]\times S^{1}}  s\  du\wedge dv. $$
The $L^{2}$ norms of $du$ and $dv$ are the same computed in the disc or the cylinder. The fact that $u,v$ are smooth on the disc implies that they are bounded and their derivatives decay exponentially as $s\rightarrow -\infty$, measured in the cylinder metric.
We have $$s \ du \wedge dv= d(s\  u\  dv) - u\  ds\wedge dv$$ and the exponential decay means that it is valid to apply Stokes' Theorem, so that
$$  I=  -  \int_{[-\infty,0]\times S^{1}}u\  v_{\theta}\  d\theta ds . $$
Let $U(s)$ be the average value of $u$ over the circle, for fixed $s$. Then for each fixed $s$
$$  \int u\ v_{\theta}\  d\theta = \int (u-U)\  v_{\theta}\  d\theta, $$
and $$\int (u-U)^{2} d\theta \leq \int u_{\theta}^{2} d\theta. $$
So $$ \vert \int  u\  v_{\theta}\  d\theta \vert^{2} \leq \int u_{\theta}^{2}\ d\theta \ \int v_{\theta}^{2}\  d\theta. $$
Now, applying the Cauchy-Schwartz inequality to the $s$ integral, we get
$$  \vert I\vert \leq \Vert du\Vert_{L^{2}} \Vert dv\Vert_{L^{2}}, $$
which is (76), with $C=1$.



\end{document}